\documentclass[12pt]{article}
\usepackage{amsthm,amsmath,amsfonts,amssymb,stackengine,comment}
\usepackage{amsfonts,mathrsfs}
\usepackage{color}
\usepackage{enumerate}
\usepackage[textwidth=8em,textsize=scriptsize]{todonotes}
\usepackage[authoryear, round]{natbib} 
\RequirePackage[colorlinks,citecolor=blue,urlcolor=blue]{hyperref}
\usepackage{graphicx}
\usepackage{float}
\usepackage{caption,subcaption}
\usepackage{booktabs,multirow}
\usepackage{lscape}

\usepackage{geometry}
\theoremstyle{plain}
\newtheorem{theorem}{Theorem}
\newtheorem{lemma}{Lemma}
\newtheorem{corollary}{Corollary}

\newtheorem{definition}{Definition}
\newtheorem{assumption}{Assumption}
\newtheorem{remark}{Remark}

\numberwithin{equation}{section}

\title{Deep Quantile Regression: Mitigating the Curse of Dimensionality Through Composition}
\author{Guohao Shen\thanks{Equal contribution.  Department of Statistics, The Chinese University of Hong Kong,  Hong Kong, China. Email: ghshen@link.cuhk.edu.hk}\ \ \
Yuling Jiao\thanks{Equal contribution. School of Mathematics and Statistics,
				Wuhan University, Wuhan, Hubei Province, China 430072. Email: yulingjiaomath@whu.edu.cn}\ \ \
Yuanyuan Lin\thanks{Department of Statistics, The Chinese University of Hong Kong,  Hong Kong, China. Email: ylin@sta.cuhk.edu.hk}\ \ \
Joel L. Horowitz\thanks{Department of Economics, Northwestern University, Evanston, IL 60208, USA. Email: joel-horowitz@northwestern.edu}\ \ \
Jian Huang\thanks{Department of Statistics and Actuarial Science,
				University of Iowa, IA 52242, USA. Email: jian-huang@uiowa.edu}
}

\begin{document}
\maketitle

\begin{abstract}
This paper considers the problem of nonparametric quantile regression under the assumption that the target conditional quantile function is a composition of a sequence of low-dimensional functions. We study the nonparametric quantile regression estimator using deep neural networks to approximate the target conditional quantile function. For convenience, we shall refer to such an estimator as a deep quantile regression  (DQR) estimator. We show that the DQR estimator achieves the nonparametric optimal convergence rate up to a logarithmic factor determined by the intrinsic dimension of the underlying compositional structure of the conditional quantile function, not the ambient dimension of the predictor. Therefore, DQR is able to mitigate the curse of dimensionality under the assumption that the conditional quantile function has a compositional structure. To establish these results, we analyze the approximation error of a composite function by neural networks and show that the error rate only depends on the dimensions of the component functions. We apply our general results to several important statistical models often used in mitigating the curse of dimensionality, including the single index, the additive, the projection pursuit, the univariate composite, and the generalized hierarchical interaction models. We explicitly describe the prefactors in the error bounds in terms of the dimensionality of the data and show that the prefactors depends on the dimensionality linearly or quadratically in these models. We also conduct extensive numerical experiments to evaluate the effectiveness of DQR and demonstrate that it outperforms a kernel-based method for nonparametric quantile regression.
\end{abstract}

{\it Keywords:}  Approximation error; composite function; deep neural networks; nonparametric regression; non-asymptotic error bound.

\section{Introduction}
\label{intro}
Consider a nonparametric regression model
\begin{equation}\label{model}
		Y=f_0(X)+\eta,
\end{equation}
where $Y \in \mathbb{R}$ is a response variable, $X \in\mathcal{X} \subset\mathbb{R}^d$ is a $d$-dimensional vector of predictors,
$f_0:\mathcal{X}
\to \mathbb{R}$ is an unknown regression function, and $\eta$ is an error term that may depend on $X$.
We consider the problem of nonparametric quantile regression under the assumption that the underlying regression function is a composition of a sequence of low-dimensional functions.
We study the nonparametric quantile regression estimator using deep neural networks to approximate the target regression function. For convenience, we shall refer to such an estimator as a deep quantile regression  (\textit{DQR}) estimator.

Quantile regression \citep{koenker1978, koenker2005} is an important method in the toolkit for analyzing the relationship between a response $Y$ and a predictor $X$.
Unlike the least squares regression that models the conditional mean of $Y$ given $X$, quantile regression estimates the conditional quantiles of $Y$ given $X$.
Thus quantile regression is able to describe the conditional distribution of $Y$  given $X$.
There is a rich literature on quantile regression, much of the work focus on the parametric case when the conditional quantile function is assumed to be a linear function of the predictor. The linear quantile regression has also been studied extensively in the context of regularized estimation and variable selection in the high-dimensional settings \citep{li2008, belloni2011, belloni2019, wang2012, zheng2015globally,
zheng2018high}.
In addition, there are many important studies on nonparametric quantile regression.  Examples include the methods using smoothing splines \citep{koenker1994, he1994convergence, he1999Qsplines}  and reproducing kernels \citep{takeuchi06a, sangnier2016joint}.  These studies established the convergence rate of the nonparametric estimators and discussed related problems arising in quantile regression, including an approach to dealing with the quantile crossing problem and  a method for incorporating prior qualitative knowledge such as monotonicity constraints in the conditional quantile function estimation. An early study on nonparametric quantile regression using shallow neural networks is \citet{white1992}.
We refer to \citet{koenker2005} and the references therein for a detailed treatment of quantile regression.
More discussions on nonparametric quantile regression related to this work are given in Section \ref{related}.

\begin{figure}[H]
	\centering
	\includegraphics[width=1\textwidth]{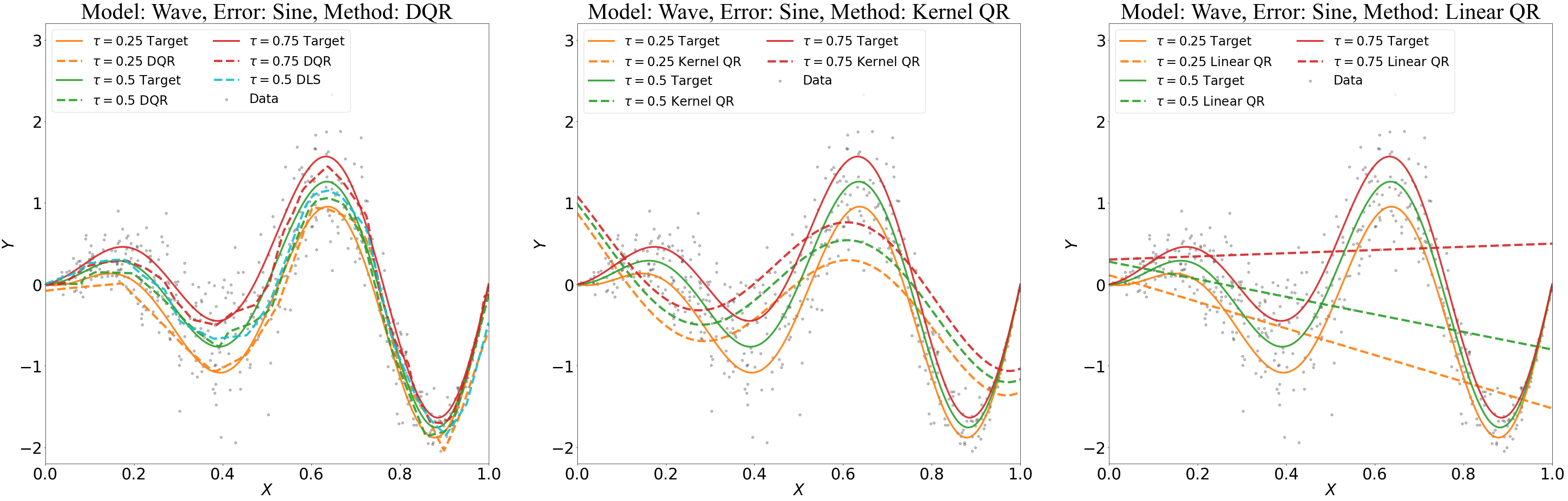}
	\caption{The fitted quantile curves by different methods under the univariate model ``Wave'' with ``Sine'' error. The training data is depicted as grey dots.The target quantile functions at the quantile levels $\tau=$0.25 (yellow), 0.5 (green), 0.75 (red)  are depicted as solid curves, and the estimated quantile functions are represented by dashed curves with the same color. From the left to right, the subfigures correspond to the methods: \textit{DQR}, \textit{kernel QR} and \textit{linear QR}. The fitted \textit{DLS}  curve (in blue)  is included in the left subfigure.}
	\label{fig:0}
\end{figure}

To give a snapshot of quantile regressions using deep neural networks compared with the traditional linear and the kernel quantile regressions, we look at the fitting of the univariate regression functions ``Wave'',  when the error term follows a ``Sine'' distribution or conditionally follows a normal distribution $(\eta\mid X=x)\sim 0.5\times\mathcal{N}(0,[\sin(\pi x)]^2).$ The functional form of the ``Wave'' function is given in Section \ref{sim}. Figure \ref{fig:0} presents the fitting results using deep quantile regression (\textit{DQR}),  quantile regression in reproducing kernel Hilbert space (\textit{kernel QR}) in \cite{sangnier2016joint} and traditional linear quantile regression (\textit{linear QR}) in \cite{koenker1978} at the $0.25$-th, the $0.50$-th and the $0.75$-th quantiles. Moreover,  least squares regression using deep neural networks (\textit{DLS}) is also compared with the above methods at the $0.50$-th quantile. We see that \textit{linear QR}  fails when the model is nonlinear, while \textit{kernel QR} and  \textit{DQR}  yield acceptable fitting curves. In particular, \textit{DQR} works best among the methods considered in this example.

In classical nonparametric statistics, including nonparametric quantile regression, the complexity of a function such as regression function and density function is measured through smoothness in terms of the order of the derivatives. The rate of convergence in estimating such functions is determined by the dimension and the smoothness index
\citep{stone1982optimal}. Specifically,
under the assumption that the target function $f_0$ is in a H\"older class with a smoothness index $\beta>0$ ($\beta$-H\"older smooth), i.e., all the partial derivatives up to order $\lfloor\beta\rfloor$ exist and the partial derivatives of order $\lfloor\beta\rfloor$ are $\beta-\lfloor\beta\rfloor$ H\"older continuous, where $\lfloor\beta\rfloor$ denotes the largest integer strictly smaller than $\beta$,
the optimal convergence rate of the prediction error is $C_d n^{-\beta/(2\beta+d)}$
under mild conditions \citep{stone1982optimal}, where $C_d$ is a prefactor independent of $n$ but depending on $d$ and other model parameters.  When $d$ is small, say,
$d =2$,  assuming the target function has a continuous second derivative, the optimal rate of convergence is
$ C_d n^{-1/3}$.  Therefore, in the low-dimensional settings, a sufficient degree of smoothness will overcome the adverse impact of the dimensionality on the convergence rate.
Moreover, in low-dimensional models with a small $d$, the impact of $C_d$ on the convergence rate is not significant.
However, in high-dimensional models with a large $d$, the situation is completely different. First, the rate of convergence can be painfully slow, unless the function $f_0$ is assumed to have an extremely large smoothness index $\beta$. But such an assumption is not realistic in practice.
Second,  the impact of $C_d$ can be substantial when $d$ is large.
For example, if the prefactor $C_d$ depends on $d$ exponentially, it can overwhelm the convergence rate $n^{-\beta/(2\beta+d)}$.
Therefore, it is important to clearly describe how $C_d$  depends on the dimensionality.

{\color{black}
Recently, several authors carried out important and inspiring studies on
the convergence properties of least squares nonparametric estimation using neural network approximation of the regression function \citep{
bauer2019deep,
schmidt2020nonparametric,
chen2019nonparametric,
kohler2019estimation, nakada2019adaptive, farrell2021deep}.
These studies show that deep neural network regression can achieve the minimax optimal rate of convergence up to a logarithmic factor for estimating the conditional mean regression function
established by \cite{stone1982optimal}.
However,   nonparametric estimation using deep neural networks cannot escape the well-know problem of {\it curse of dimensionality} in high-dimensions without any conditions on the underlying model.
}

It is clear that smoothness is not the right measure of the complexity of a function class in the high-dimensional settings, since smoothness does not help mitigate the curse of dimensionality. An effective approach to mitigating the curse of dimensionality is to consider functions with a compositional structure.
Deep neural network modeling has achieved impressive success and often outperformed
kernel based methods in many important applications with high-dimensional data, including speech recognition, image classification, object detection, drug discovery and genomics, among others \citep{lecun2015deeplearning}. Thus it is desirable to consider statistical models in a function class that can mitigate the curse of dimensionality and can be well approximated by deep neural networks.
It has been shown that deep ReLU networks are
solutions to regularized data fitting problems in the function space consisting of
compositions of functions from the Banach spaces of second-order bounded variation \citep{parhi2021kinds}.
Using composite functions in nonparametric regression modeling has a long history in statistics. For example, the nonparametric additive model, which can be considered a composition of a linear function with a vector function whose components depend on only one of the variables,
has been studied by many authors \citep{
breiman1985additive,  stone1985,stone1986dimensionality,  hastie1990generalized}.
Recently, more general composite functions for statistical modeling have been proposed in several interesting works \citep{horowitz2007, bauer2019deep, schmidt2020nonparametric}.
Under this assumption, the convergence rate $C_dn^{-\beta/(2\beta+d)}$  could be improved to $C_{d,d_*} n^{-\beta/(2\beta+d_*)}$
for some $d_*\ll d$, where $C_{d,d_*}$
is a constant depending on $(d_*,d)$, where $d_*$ is the intrinsic dimension of the model.
In these results, the convergence rate part is improved from $n^{-\beta/(2\beta+d)}$
to $n^{-\beta/(2\beta+d_*)}$. When $d_* \ll d$, the improvement is substantial.
However, the prefactor $C_{d, d_*}$  in the error bounds depends on $d$ exponentially or are not clearly described in the aforementioned works \citep{stone1985,stone1986dimensionality,horowitz2007, bauer2019deep, schmidt2020nonparametric}.
In a low-dimensional model with a
small $d$, the impact of the prefactor on the overall error bound is not significant. However, in a high-dimensional model with a large $d$, the impact of the prefactor can be substantial, even overwhelm the convergence rate part \citep{ghorbani2020discussion}. Therefore, it is important to describe how the prefactor depends on the dimension $d$ in the error bound.
	
In this paper, we establish non-asymptotic upper bounds for the excess risk and mean integrated squared error of the \textit{DQR} estimator under the assumption that the target regression function is a composite function. A novel aspect of our work is that we clearly describe how the prefactors in the error bounds depend on the ambient dimension $d$ and the dimensions of the low-dimensional component functions of the composite function.  Our error bounds achieve the minimax optimal rates and significantly improve over the existing ones in the sense that their prefactors depend  linearly or quadratically on the dimension $d$, instead of exponentially on $d$. This shows that  \textit{DQR} can mitigate the curse of dimensionality under the assumption that the target regression function belongs to the class of composite functions.
These results are based on new approximation error bounds of composite functions by the neural networks, which may be of independent interest.
Our main contributions are as follows.

\begin{enumerate}
\item We establish  excess risk bounds for the proposed \textit{DQR} estimator
under the assumption that the target conditional quantile function has a compositional structure with lower-dimensional component functions. With appropriately specified ReLU networks in terms of depth, width and size of the network, our \textit{DQR} estimator achieves near optimal convergence rate up to a logarithmic factor under a heavy-tailed error (finite $p$-th moment for $p\geq1$) and mild regular conditions on the joint distribution of the response and the predictor. Moreover, we show that DQR can mitigate the curse of dimensionality in the sense that the convergence rate of the error bound depends on the dimensions of the component functions, not the ambient dimension. We also show that the prefactors of the error bounds depend on the ambient dimension linearly or quadratically.

\item We derive novel approximation error results of composite functions using ReLU
activated neural networks under the assumption that the component functions are H\"older continuous. This result shows that the curse of dimensionality can be mitigated through composition in the sense the approximate error rate depends on the intrinsic dimension of a composite functions, instead of the ambient dimension of the function.  Equally importantly, the prefactor of the error bound is significantly improved in the sense that it depends on the dimensionality $d$ polynomially instead of exponentially as in the existing results. This approximation result is the key building block in establishing the bounds for excess risk and mean integrated squared error for \textit{DQR}.

\item We apply our general results to several important statistical models often used in
 mitigating the curse of dimensionality, including the single index, the additive, the projection pursuit, the univariate composite, and the generalized hierarchical interaction models. We show that \textit{DQR}
 achieves the optimal convergence rate up to a logarithmic factor under these models.
  We also present the prefactors of the error bounds for these models.

\item We bridge the gap between the excess risk and the mean integrated squared error of
the \textit{DQR} estimator under mild conditions. We do not require the bounded support condition  on the conditional distribution of the response given the predictor as in the existing literature. The mean integrated squared error of our \textit{DQR} estimator is shown to converge at the near optimal rate up to a logarithmic factor, inheriting the properties of
the corresponding excess risk. The convergence rate of the mean integrated squared error of the \textit{DQR} estimator is determined by the dimensions of the component functions and the prefactor depends polynomially on the widest layer of the composite functions.
\end{enumerate}

The remainder of this paper is organized as follows. In Section \ref{sec2} we describe the deep quantile regression problem,  the deep neural networks used in the estimation and the assumption on the compositional structure of the conditional quantile function. In Section \ref{summary} we provide a high level description of our main results and the overall approach we take to establish these results. In Section \ref{sec4} we
present non-asymptotic bounds on the excess risk and mean integrated squared error of the \textit{DQR} estimator. Section \ref{compositionR} includes applications of our general error bounds to  several important models in nonparametric statistics.  In Section \ref{sec3} we present a result on the approximation error of composite functions using deep neural networks. In Section \ref{sim} we present simulation results demonstrating that \textit{DQR} outperforms a kernel nonparametric quantile regression method
based on vector-valued reproducing kernel Hilbert space (RKHS) \citep{sangnier2016joint}.
Section \ref{related} contains discussions on the related work.  Concluding remarks are given in Section \ref{conclusion}. Proofs and additional simulation results are given in the appendix.

\section{Deep quantile regression}
\label{sec2}
In this section, we present the basic setup of nonparametric regression.
We describe the structure of the feedforward neural networks to be used in the
estimation and define the compositional structure for the target conditional quantile function.

For a given quantile level $\tau\in(0,1)$,  the quantile check loss function is defined by $$\rho_\tau(x)=x\{\tau-I(x\leq0)\}, \ x \in \mathbb{R}.$$
For a possibly random function $f:\mathbb{R}^d\to\mathbb{R}$, let $Z\equiv (X,Y)$ be a random vector independent of $f$.  We define the risk of $f$ under the loss function $\rho_\tau(\cdot)$ by
\begin{equation*}
	\mathcal{R}^\tau(f)=\mathbb{E}_Z\{\rho_\tau(Y-f(X))\}.
\end{equation*}
At the population level, the nonparametric quantile estimation is to find a measurable function
$f^*: \mathbb{R}^d\to \mathbb{R}$ satisfying
	$$f^* :=\arg\min_{f} \mathcal{R}^\tau(f) =\arg\min_{f}\mathbb{E}_Z\{\rho_\tau(Y-f(X))\},$$
where $\mathbb{E}_Z$ means that the expectation is taken with respect to the distribution of $Z.$
If the conditional $\tau$-th quantile of $\eta$ given $X$ is 0 and
$\mathbb{E}(|\eta|\vert X=x)<\infty$ for all $x\in\mathcal{X}$, then the true regression function $f_0$  is the optimal solution $f^*$ on $\mathcal{X}$.	

In  applications, when only a random sample $S \equiv \{(X_i,Y_i)\}_{i=1}^n$ is available, we consider the empirical risk
\begin{equation}
\label{er1}
\mathcal{R}^\tau_n(f)=\frac{1}{n} \sum_{i=1}^{n} \rho_\tau(Y_i-f(X_i)).
\end{equation}
Our goal is to construct an estimator of $f_0$ within a certain class of functions $\mathcal{F}_n$ by minimizing the empirical risk, that is,
	\begin{equation}\label{erm}
		\hat{f}_n\in\arg\min_{f\in\mathcal{F}_n}\mathcal{R}^\tau_n(f),
	\end{equation}
where $\hat{f}_n$ is called the empirical risk minimizer (ERM).
We choose $\mathcal{F}_n$  to be a function class consisting of
deep neural networks (DNN). We will also refer to $\hat{f}_n$ as a deep quantile regression (DQR) estimator below.

\subsection{Deep neural networks}

We set the function class $\mathcal{F}_n$
to be $\mathcal{F}_{\mathcal{D},\mathcal{W}, \mathcal{U},\mathcal{S},\mathcal{B}}$, a class of  feedforward neural networks $f_\phi: \mathbb{R}^d \to \mathbb{R} $ with parameter $\phi$, depth $\mathcal{D}$, width $\mathcal{W}$, size $\mathcal{S}$, number of  neurons $\mathcal{U}$ and $f_{\phi}$ satisfying $\Vert f_\phi\Vert_\infty\leq\mathcal{B}$ for some $0 <B < \infty$, where
$\Vert f \Vert_\infty$ is the supreme norm of a function $f:\mathbb{R}^d\to\mathbb{R}$.
Note that the network parameters may depend on the sample size $n$, but  the dependence is omitted in the notation for simplicity.
A brief description of multilayer perceptrons (MLPs), the commonly used feedforward neural networks, are given below. The architecture of a MLP can be expressed as a composition of a series of functions
\[
f_\phi(x)=\mathcal{L}_\mathcal{D}\circ\sigma\circ\mathcal{L}_{\mathcal{D}-1}
\circ\sigma\circ\cdots\circ\sigma\circ\mathcal{L}_{1}\circ\sigma\circ\mathcal{L}_0(x),\  x\in \mathbb{R}^d,
\]
where $\sigma(x)=\max(0, x)$ is the rectified linear unit (ReLU) activation function (defined for each component of $x$ if $x$ is a vector) and
$$\mathcal{L}_{i}(x)=W_ix+b_i,\quad i=0,1,\ldots,\mathcal{D},$$
where $W_i\in\mathbb{R}^{d_{i+1}\times d_i}$  is a weight matrix,  $d_i$ is the width (the number of neurons or computational units) of the $i$-th layer, and $b_i\in\mathbb{R}^{d_{i+1}}$ is the bias vector in the $i$-th linear transformation $\mathcal{L}_i$.

Such a network $f_\phi$ has $\mathcal{D}$ hidden layers and $(\mathcal{D}+1)$ layers in total.
We use a $(\mathcal{D}+1)$-vector $(w_0,w_1,\ldots,w_\mathcal{D})^\top$ to describe the width of each layer; particularly in nonparametric regression problems,  $w_0=d$ is the dimension of the input and $w_\mathcal{D}=1$ is the dimension of the response . The width $\mathcal{W}$ is defined as the maximum width of hidden layers, i.e.,
$\mathcal{W}=\max\{w_1,\ldots,w_\mathcal{D}\}$; 	the size $\mathcal{S}$ is defined as the total number of parameters in the network $f_\phi$, i.e., $\mathcal{S}=\sum_{i=0}^\mathcal{D}\{w_{i+1}\times(w_i+1)\}$; 	
the number of neurons $\mathcal{U}$ is defined as the number of computational units in hidden layers, i.e.,
$\mathcal{U}=\sum_{i=1}^\mathcal{D} w_i$.  For an MLP  $\mathcal{F}_{\mathcal{D},\mathcal{U},\mathcal{W},\mathcal{S},\mathcal{B}}$,
its parameters satisfy the simple
relationship
\[
\max\{\mathcal{W},\mathcal{D}\}\leq\mathcal{S}\leq\mathcal{W}(\mathcal{D}+1)
+(\mathcal{W}^2+\mathcal{W})(\mathcal{D}-1)+\mathcal{W}+1
=O(\mathcal{W}^2\mathcal{D}).
\]

\subsection{Structured composite functions}
\label{composition}

Let the target quantile regression function $f_0: \mathbb{R}^d \to \mathbb{R}$ be a $d$-dimensional function.
We assume that $f_0$ is a composition of a series of functions $h_i, i=0\ldots,q$, i.e.,
$$f_0=h_q\circ\cdots\circ h_0,$$
where $h_i:[a_i,b_i]^{d_i}\to[a_{i+1},b_{i+1}]^{d_{i+1}}$. Here $d_0=d$ and $d_{q+1}=1$.
For each $h_i$, denote by $h_i=(h_{ij})^\top_{j=1,\ldots,d_{i+1}}$ the components of $h_i$ and let $t_i$ be the maximal number of variables on which each of $h_{ij}$ the depends on. Note that $t_i\leq d_i$ and each $h_{ij}$ is a $t_i$-variate function for $j=1,\ldots,d_i$.


Many well-known important models in  semiparametric  and nonparametric statistics
have a compositional structure. Examples include the single index model \citep{hardle1993optimal, 
horowitz1996},
 the additive model  \citep{stone1985, stone1986dimensionality, hastie1990generalized},
 the projection pursuit model  \citep{friedman1981projection},
the interaction model  \citep{stone1994use},  the composite regression model \citep{horowitz2007},
 and  the {generalized hierarchical interaction model}
\citep{bauer2019deep}.
We consider the bounds for the excess risk of \textit{DQR} under these models in Section \ref{compositionR}.

In this work, we focus on the quantile regression models in which the conditional quantile function has a compositional structure. This is the key condition we use to mitigate the curse of dimensionality.
We will only assume the H\"older continuity on the component functions of the composite conditional quantile function.
A function $h: [a_1, b_1]^{d_1} \to [a_2, b_2]^{d_2}$ is said to be
 H\"older continuous with order $\alpha$ and H\"older constant $\lambda$ if there exist $\alpha\in(0,1]$ and $\lambda \ge 0$ such that
 \begin{equation}
 \label{holder}
 \Vert h(x)-h(y)\Vert_2 \leq \lambda\Vert x-y\Vert_2^\alpha
 \end{equation}
  for any $x,y\in[a_1, b_1]^{d_1}$.

We now describe the assumptions on the target regression function $f_0$ in detail below.

\begin{assumption}[Structured target regression function with continuous components]
	\label{structure}
	The target quantile regression function $f_0=h_q\circ\cdots\circ h_0$ is a composition of a series of functions $h_i, i=0\ldots,q$,  where $h_i:[a_i,b_i]^{d_i}\to[a_{i+1},b_{i+1}]^{d_{i+1}}$ with $d_0=d$ and $d_{q+1}=1$. For each $h_i=(h_{ij})^\top_{j=1,\ldots,d_{i+1}}$ ($i=0,\ldots,q$), its components $h_{ij}:[a_i,b_i]^{t_{i}}\to[a_{i+1},b_{i+1}]$ ($j=1,\ldots,d_{i+1}$) are H\"older continuous functions with order $\alpha_i\in[0,1]$ and constant $\lambda_i\ge0$, where $t_i$ is the maximal number of variables on which each of $h_{ij}$ depends on ($t_i\leq d_i$). Let $J\subset\{0,\ldots,q\}$ be a set consisting of the indices of linear transformation layers of $f_0$ (if any) and $J^c:=\{0,\ldots,q\}\backslash J$ denote the complement of $J$.
\end{assumption}
We will show that,
if the target regression function $f_0$ satisfies Assumption \ref{structure},
the \textit{DQR} estimator can automatically adapt to the compositional structure and circumvent the curse of dimensionality.

\section{A high-level description of the results}
\label{summary}

In this section, we present a high-level description of our approach, the non-asymptotic bounds for the excess risk and the mean integrated squared error of the \textit{DQR} estimator. Detailed statements of the results and the assumptions are given in the Sections \ref{sec4}-\ref{sec3} below.

For a \textit{DQR} estimator $\hat{f}_n \in \mathcal{F}_n$ defined in (\ref{erm}),
we evaluate its quality via the \textit{excess risk}, defined as the difference between the risks of $\hat{f}_n$ and $f_0$,
\begin{align*}
	\mathcal{R}^\tau(\hat{f}_n)-\mathcal{R}^\tau(f_0) = \mathbb{E}_{Z}\rho_\tau(\hat{f}_n(X)-Y)-\mathbb{E}_{Z}\rho_\tau(f_0(X)-Y). \nonumber
\end{align*}
{\color{black}
We first establish an upper bound on the excess risk, which is the starting point of our error analysis.

\begin{lemma}\label{lemma1}
For any random sample $S=\{(X_i, Y_i)_{i=1}^n\}$, the excess risk of the
\textit{DQR} estimator $\hat{f}_n$
satisfies
		\begin{align}
\label{lem1}
			\mathcal{R}^\tau(\hat{f}_n)-\mathcal{R}^\tau(f_0)\leq  2\sup_{f\in\mathcal{F}_n}\vert \mathcal{R}^\tau(f)-\mathcal{R}^\tau_{n}(f)\vert+\inf_{f\in\mathcal{F}_n} \mathcal{R}^\tau({f})-\mathcal{R}^\tau(f_0),
		\end{align}
where $\mathcal{R}^\tau_{n}$ is defined in (\ref{er1}).
	\end{lemma}	
	
The excess risk of the \textit{DQR} estimator is bounded above by the sum of two terms: the stochastic error $2\sup_{f\in\mathcal{F}_n}\vert \mathcal{R}^\tau(f)-\mathcal{R}^\tau_n(f)\vert$ and the approximation error $\inf_{f\in\mathcal{F}_n}\mathcal{R}^\tau(f)- \mathcal{R}(f_0)$. It is interesting to note that the upper bound no longer depends on the \textit{DQR} estimator itself, but the function class $\mathcal{F}_n$, the loss function $\rho_\tau$ and the random sample $S$.}

The stochastic error  $2\sup_{f\in\mathcal{F}_n}\vert \mathcal{R}^\tau(f)-\mathcal{R}^\tau_n(f)\vert$ can be analyzed using the empirical process theory \citep{vw1996, anthony1999, bartlett2019nearly}.
A key step is to calculate the complexity measure of $\mathcal{F}_n$ in terms of its covering number.
The details are given in Section \ref{sec4}.

The approximation error term $\inf_{f\in\mathcal{F}_n} \mathcal{R}^\tau({f})-\mathcal{R}^\tau(f_0)$ measures the approximation error of the function class $\mathcal{F}_n$ for $f_0$ under
the loss function $\rho_\tau$. To utilize the approximation theories of neural networks, we need to relate $\inf_{f\in\mathcal{F}_n} \mathcal{R}^\tau({f})-\mathcal{R}^\tau(f_0)$ to the quantity $\inf_{f\in\mathcal{F}_n} \Vert f-f_0\Vert$ for some functional norm $\Vert\cdot\Vert$.
The power of neural network functions approximating high-dimensional functions have been studied by many authors, some recent works include \citet{yarotsky2017error, yarotsky2018optimal,shen2019nonlinear, shen2019deep}, among others. For a composite function $f_0$ under Assumption  \ref{structure}, we derive new approximation results in Section \ref{sec3}.

To clearly describe how the error bounds depend on various parameters, including the network parameters such as depth, width and size of the network,  as well as the model parameters such as the intrinsic and ambient dimensions of the model, we present general expressions of the stochastic errors and the approximation errors, which constitute
the upper bounds for the excess risk and the mean integrated squared error (MISE),  in Theorems \ref{thm2} and \ref{thm3} in Section \ref{sec4} below. The network parameters, similar to the bandwidth in kernel nonparametric regression or density estimation, can be tuned as a function of the sample size and the model dimension to obtain the best trade-off between the stochastic error and the approximation error, and therefore achieve the best overall error rate.
An appealing aspect of our results is that they clearly and explicitly describe how the prefactors in the error bounds depend on the network parameters and the dimensionality of the model.
Explicit expressions of the bounds for the excess risk and the MISE are presented in Corollaries \ref{thm2d} and \ref{thm3c} in Section \ref{sec4}.

In Section \ref{compositionR}, we consider several well-known semiparametric and nonparametric models that are widely used to mitigate the curse of dimensionality,  including the single index model, the additive model, the projection pursuit model,  the interaction model,  the univariate composite regression model, and the generalized hierarchical interaction model. We derive explicit expressions of the error bounds when the underlying conditional quantile function takes the form of these well-known models

As can be seen in Corollary \ref{thm2d} for the excess risk of the \textit{DQR} estimator and the error bounds for the models considered in Section \ref{compositionR},  based on appropriately specified network parameters (depth, width and size of the network), we have the following upper bound for the excess risk,
	\begin{equation}
		\label{erb0}
		\mathbb{E}\big\{\mathcal{R}^\tau(\hat{f}_\phi)-\mathcal{R}^\tau(f_0)\big\}\leq 	C_0 C_{d, d^*} (\log n)^2 n^{-\left(1-\frac{1}{p}\right)\frac{2\alpha^*}{2\alpha^*+t^*}},
	\end{equation}
where $C_0$ is a constant only depending on the model parameters such as the smoothness index of the underlying conditional quantile function,  $C_{d,d^*}$ is the prefactor depending on $d$,  the dimension of the predictor;  and $d^*$, determined by the dimensions of the component functions in the composite function. The convergence rate part of the error bound (\ref{erb0}), $n^{-(1-1/p)2\alpha^*/(2\alpha^*+t^*)}$, is determined by the number of moments $p$ of the response $Y$ (see Assumption \ref{moment} below), the smoothness index of the composite function $\alpha^*$, and the intrinsic dimension of the model $t^*$. If $Y$ has sub-exponential tail probabilities, we can set $p=\infty$.
The bound for the mean integrated squared error of the \textit{DQR} estimator has a form similar to
(\ref{erb0}), see Corollary \ref{thm3c}.

Explicit expressions for $C_{d, d^*}$ in (\ref{erb0})
 are given in Corollaries \ref{thm2d} and \ref{thm3c}, as well as for the examples in Section \ref{compositionR}. For example, for the single index model (\ref{siM}), the additive model (\ref{additiveM})  and the additive model with an unknown link function (\ref{additiveUM}),
$C_{d, d^*} = d^2 \log d$. For the interaction model  (\ref{interactionM}),
$C_{d, d^*}=(Kdd^*)^2\log(Kdd^*),$ where $K$ is the number of component functions  and $d^*$ is the dimension of the component functions in the model. For the projection pursuit model (\ref{ppM}),  $C_{d, d^*}=(\max\{K, d\})^2 \log (\max\{K, d\}),$ where $K$
is the number of component functions in the model. For the univariate composite model (\ref{ucM}) and the generalized hierarchical interaction model (\ref{ghiMe}), the forms of $C_{d, d^*}$ are more complicated, they are given in Section \ref{compositionR}.

These results  demonstrate that \textit{DQR} with deep neural networks can significantly attenuate the curse of dimensionality when the underlying conditional quantile function takes the form of one of these models, even though the construction of the \textit{DQR} estimator does not
use the specific structure of these models.

\section{Non-asymptotic error bounds}
\label{sec4}
	
In this section,
we  present non-asymptotic error bounds for the \textit{DQR} estimator, including bounds for the excess risk upper bounds in section \ref{sec4.1} and bounds for mean integrated squared error in \ref{sec4.2}. The bounds are determined by a trade-off between the stochastic error and the approximation error.

\subsection{Excess risk bounds}
\label{sec4.1}



For analyzing the stochastic error of the \textit{DQR} estimator, we make the following assumption.

{\color{black}
\begin{assumption}
\label{moment} (i)
The conditional $\tau$-th quantile of $\eta$ given $X=x$ is 0 and $\mathbb{E}(\vert \eta\vert |X=x)<\infty$ for almost every $x\in\mathcal{X}$. (ii) The support of covariates $\mathcal{X}$ is a bounded compact set in $\mathbb{R}^d$, and without loss of generality $\mathcal{X}=[0,1]^d$. (iii) The response variable $Y$ has a finite $p$-th moment for some  $p>1$, i.e., there exists a finite constant $M>0$ such that $\mathbb{E}\vert Y\vert^p\leq M$.
\end{assumption}
Note that throughout the paper, we focus on the case when $\mathcal{X}=[0,1]^d$.
In the nonparametric regression problems, we can always first transform the predictors to a bounded region.
}

For a class $\mathcal{F}$ of functions: $\mathcal{X}\to \mathbb{R}$,
its pseudo dimension, denoted by $\text{Pdim}(\mathcal{F}),$  is defined to be the largest integer $m$ for which there exists $(x_1,\ldots,x_m,y_1,\ldots,y_m)\in\mathcal{X}^m\times\mathbb{R}^m$ such that for any $(b_1,\ldots,b_m)\in\{0,1\}^m$ there exists $f\in\mathcal{F}$ such that $\forall i:f(x_i)>y_i\iff b_i=1$ \citep{anthony1999, bartlett2019nearly}.
For a class of real-valued functions generated by neural networks, pseudo dimension is a natural measure of its complexity. In particular, if $\mathcal{F}$ is the class of functions generated by a neural network with a fixed architecture and fixed activation functions, we have $\text{Pdim}(\mathcal{F})=\text{VCdim}(\mathcal{F})$ (Theorem 14.1 in \cite{anthony1999}), where $\text{VCdim}(\mathcal{F})$ is the VC dimension of $\mathcal{F}$.
In our results, we require the sample size $n$ to be greater than the pseudo dimension of the class of neural networks considered.
	
For a given sequence $x=(x_1,\ldots,x_n)\in\mathcal{X}^n,$ let  $\mathcal{F}_\phi|_x=\{(f(x_1),\ldots,f(x_n):f\in\mathcal{F}_\phi\} \subset \mathbb{R}^{n}$. For a positive number $\delta$, let $\mathcal{N}(\delta,\Vert\cdot\Vert_\infty,\mathcal{F}_\phi|_x)$ be the covering number of $\mathcal{F}_\phi|_x$ under the norm $\Vert\cdot\Vert_\infty$ with radius $\delta$.
Define the uniform covering number
$\mathcal{N}_n(\delta,\Vert\cdot\Vert_\infty,\mathcal{F}_\phi)$ to be the maximum
	over all $x\in\mathcal{X}$ of the covering number $\mathcal{N}(\delta,\Vert\cdot\Vert_\infty,\mathcal{F}_\phi|_x)$, i.e.,
\begin{equation}
\label{ucover}
\mathcal{N}_n(\delta,\Vert\cdot\Vert_\infty,\mathcal{F}_\phi)=
\max\{\mathcal{N}(\delta,\Vert\cdot\Vert_\infty,\mathcal{F}_\phi|_x):x\in\mathcal{X}\}.
\end{equation}
We give an upper bound of the stochastic error in the following lemma.

\begin{lemma}\label{lemma2}
Consider the $d$-variate nonparametric regression model in (\ref{model}) with an unknown regression function $f_0$. Let $\mathcal{F}_\phi=\mathcal{F}_{\mathcal{D},\mathcal{W},\mathcal{U},\mathcal{S},\mathcal{B}}$ be a class of feedforward neural networks with a continuous piecewise-linear activation function of finite pieces and $\hat{f}_\phi \in\arg\min_{f\in\mathcal{F}_\phi}R^\tau_n(f) $ be the empirical risk minimizer over $\mathcal{F}_\phi$. Assume that Assumption \ref{moment} holds and $ \Vert f_0\Vert_\infty\leq \mathcal{B}$ for $\mathcal{B}\geq 1$.
Then,  for $2n \ge \text{Pdim}(\mathcal{F}_\phi)$ and any $\tau\in(0,1)$,
\begin{equation} \label{entropy}
	\sup_{f\in\mathcal{F}_\phi}\vert \mathcal{R}^\tau(f)-\mathcal{R}^\tau_{n}(f)\vert \leq c_0\frac{\max\{\tau,1-\tau\}\mathcal{B}}{n^{1-1/p}} \log\mathcal{N}_{2n}(n^{-1},\Vert \cdot\Vert_\infty,\mathcal{F}_\phi),
\end{equation}
where $c_0>0$ is a constant independent of $n,d,\tau,\mathcal{B},\mathcal{S},\mathcal{W}$ and $\mathcal{D}$. Moreover,
\begin{equation} \label{oracle}
	\mathbb{E}\big\{\mathcal{R}^\tau(\hat{f}_\phi)-\mathcal{R}^\tau(f_0)\big\}\leq C_0\frac{\max\{\tau,1-\tau\}\mathcal{B}\mathcal{S}\mathcal{D}\log(\mathcal{S})\log(n)}{n^{1-1/p}}+ 2\inf_{f\in\mathcal{F}_\phi}\big\{\mathcal{R}^\tau(f)-\mathcal{R}^\tau(f_0)\big\},
\end{equation}
where $C_0>0$ is a constant independent of $n,d,\tau,\mathcal{B},\mathcal{S},\mathcal{W}$ and $\mathcal{D}$.
	\end{lemma}


\begin{remark}
	The denominator $n^{1-1/p}$ in (\ref{entropy}) and (\ref{oracle}) can be improved to $n$ if the response $Y$ is assumed to be sub-exponentially distributed, i.e., there exists a constant $\sigma_Y>0$ such that $\mathbb{E}\exp(\sigma_Y\vert Y\vert)<\infty$. This corresponds to the case that $p=+\infty$.
\end{remark}

The stochastic error is bounded by a term determined by the metric entropy of $\mathcal{F}_\phi$
in (\ref{entropy}), which is measured by the covering number  of $\mathcal{F}_\phi$.
To obtain (\ref{oracle}), we further bound the covering number of $\mathcal{F}_\phi$ by its pseudo dimension (VC dimension).
According to \cite{bartlett2019nearly}, the pseudo dimension (VC dimension) of $\mathcal{F}_\phi$ with piecewise-linear activation function can be further contained and expressed in terms of its parameters $\mathcal{D}$ and $\mathcal{S}$, i.e.,  ${\rm Pdim}(\mathcal{F}_\phi)=O(\mathcal{S}\mathcal{D}\log(\mathcal{S}))$. This leads to the upper bound for the prediction error by the sum of the stochastic error and  the approximation error of $\mathcal{F}_\phi$ to $f_0$ in (\ref{oracle}).
	
To derive an upper bound for the approximation error $\inf_{f\in\mathcal{F}_\phi}\{\mathcal{R}^\tau(f) -\mathcal{R}^\tau(f_0)\}$, we first bound it in terms of $\inf_{f\in\mathcal{F}_{\phi}}\Vert f-f_0\Vert$ for some functional norm $\Vert\cdot\Vert$. In the following, we let $\nu$ denote the marginal distribution of $X$ and define $\Vert f-f_0\Vert_{L^p(\nu)}:=\{\mathbb{E} \vert f(X)-f_0(X) \vert^p\}^{1/p}$ for $p\in(0,\infty)$.

\begin{lemma}\label{lemma3}
	Assume that Assumption \ref{moment} (i) holds. Let $f_0$ be the target function defined in (\ref{model}) and $\mathcal{R}^\tau(f_0)$ be its risk. Then, we have
	$$\inf_{f\in\mathcal{F}_\phi}\{\mathcal{R}^\tau(f) -\mathcal{R}^\tau(f_0)\}\leq \max\{\tau,1-\tau\}\inf_{f\in\mathcal{F}_{\phi}}\mathbb{E}\vert f(X)-f_0(X)\vert=\max\{\tau,1-\tau\}\inf_{f\in\mathcal{F}_{\phi}}\Vert f-f_0\Vert_{L^1(\nu)},$$
where $\nu$ denotes the marginal distribution of $X$.
\end{lemma}
As a consequence of Lemma \ref{lemma3}, we only need to give upper bounds on the approximation error $\inf_{f\in\mathcal{F}_{\phi}}\Vert f-f_0\Vert_{L^1(\nu)}$ to give the overall bounds on the excess risk of the ERM $\hat{f}_\phi$ defined in (\ref{erm}).
Furthermore, if the conditional distributions of error given covariates satisfy proper  conditions and the risk function $\mathcal{R}(\cdot)$ has a local quadratic approximation around $f_0$, the convergence rate results can be further improved.

\begin{assumption}[Local quadratic bound of the excess risk]
	\label{quadratic}
	There exist some constants $c^0_\tau=c^0_\tau(\tau,X,\eta,f_0)>0$ and $\delta^0_\tau=\delta^0_\tau(\tau,X,\eta,f_0)>0$ which may depend on $\tau$, $X$, $\eta$ and $f_0$ such that
	$$\mathcal{R}^\tau(f) -\mathcal{R}^\tau(f_0)\leq c^0_\tau \Vert f-f_0 \Vert^2_{L^2(\nu)},$$
	for any $f$ satisfying $\Vert f-f_0 \Vert_{L^\infty(\mathcal{X}^0)}\leq \delta^0_\tau$, where $\mathcal{X}^0$ is any subset of $\mathcal{X}$ such that $P(X\in\mathcal{X}^0)=P(X\in\mathcal{X})$.
\end{assumption}
\begin{remark}
 Assumption \ref{quadratic} is generally satisfied when the conditional density of $\eta$ given $X=x$ is positive
 in a neighborhood of its $\tau$-th conditional quantile.
\end{remark}

By  Lemma \ref{lemma3} and Assumption \ref{quadratic},  a sharper bound for the approximation error improves over that of Lemma \ref{lemma3}
can be obtained and presented in the next lemma.

\begin{lemma}\label{lemma4}
	Assume that Assumption \ref{moment} (i) and \ref{quadratic} hold, let $f_0$ be the target function defined in (\ref{model}) and $\mathcal{R}^\tau(f_0)$ be its risk, then we have
	$$\inf_{f\in\mathcal{F}_\phi}\{\mathcal{R}^\tau(f) -\mathcal{R}^\tau(f_0)\}\leq c_\tau\inf_{f\in\mathcal{F}_{\phi}}\Vert f-f_0\Vert^2_{L^2(\nu)},$$
	where $c_\tau\geq \max\big\{c^0_\tau,\max\{\tau,1-\tau\}/\delta^0_\tau\big\}>0$ is a constant, $\nu$ denotes the marginal probability measure of $X$ and $\mathcal{F}_\phi=\mathcal{F}_{\mathcal{D},\mathcal{W},\mathcal{U},\mathcal{S},\mathcal{B}}$ denotes the class of feedforward neural networks with parameters $\mathcal{D},\mathcal{W},\mathcal{U},\mathcal{S}$ and $\mathcal{B}$.
\end{lemma}

\begin{remark}
We establish the error bounds for approximating a composite function using deep neural networks in Theorem \ref{thm1} in Section \ref{sec3}. Theorem \ref{thm1} can be used
to bound the approximation error term $\inf_{f\in\mathcal{F}_{\phi}}\Vert f-f_0\Vert_{L^2(\nu)}$ in Lemmas \ref{lemma3} and \ref{lemma4}, which leads to the bound for the approximation error in Theorem \ref{thm2} below.
\end{remark}

	

Before stating the results for the excess risk bounds, we specify the network parameters.
For any given $N_i,L_i\in\mathbb{N}^+, i\in J^c$, we set the function class $\mathcal{F}_{\phi}=\mathcal{F}_{\mathcal{D},\mathcal{W},\mathcal{U},\mathcal{S},\mathcal{B}}$
consisting of ReLU multi-layer perceptrons
with width no more than $\mathcal{W}$ and depth $\mathcal{D}$, where
\begin{align}
\label{widtha}
\mathcal{W}&=\max_{i=0,\ldots,q}d_{i}\max\{4t_i\lfloor N_i^{1/t_i}\rfloor+3t_i,12N_i+8\}, \\
 \label{deptha}
 \mathcal{D}&=\sum_{i\in J^c}(12L_i+15)+2\vert J\vert.
 \end{align}
 Here recall $J\subset\{0,\ldots,q\}$ is a set collecting the indices of linear layers of $f_0$ (if any) and $J^c:=\{0,\ldots,q\}\backslash J$ denotes the complement of $J$.

\begin{theorem}[Non-asymptotic excess risk bound] \label{thm2}
Under model (\ref{model}), suppose that Assumptions \ref{structure} and \ref{moment} hold, $\nu$ is absolutely continuous with respect to the Lebesgue measure, and $ \Vert f_0\Vert_\infty\leq\mathcal{B}$ for some $\mathcal{B}\geq1$. Suppose the network parameters of the function class $\mathcal{F}_{\phi}$ are specified as in (\ref{widtha}) and (\ref{deptha}).
Then,  for $2n \ge \text{Pdim}(\mathcal{F}_\phi)$,
the excess risk
of the \textit{DQR} estimator $\hat{f}_\phi$ satisfies
\begin{equation*}
	\mathbb{E}\big\{\mathcal{R}^\tau(\hat{f}_\phi)-\mathcal{R}^\tau(f_0)\big\}\leq C\frac{\lambda_\tau\mathcal{B}\mathcal{S}\mathcal{D}\log(\mathcal{S})\log(n)}{n^{1-1/p}}+ 2\lambda_\tau\sum_{i\in J^c} C_i^*\lambda_i^* t_i^*(N_iL_i)^{-2\alpha_i^*/t_i},
\end{equation*}
where $\lambda_\tau=\max\{\tau,1-\tau\}$ and $C>0$ is a constant which does not depend on $n,d,\tau,\mathcal{B},$ $\mathcal{S},$ $\mathcal{D},$ $C_i^*,$ $\lambda_i^*,$ $\alpha_i^*,$ $N_i$ or $L_i$, and $C_i^*=18^{\Pi_{j=i+1}^{q}\alpha_{j}}$, $\lambda_i^*=\Pi_{j=i}^{q}\lambda_{j}^{\Pi_{k=j+1}^{q} \alpha_{k}}$, $\alpha_i^*=\Pi_{j=i}^q \alpha_j$ and $t_i^*={(\Pi_{j=i}^{q}\sqrt{t_{j}}^{\Pi_{k=j}^{q} \alpha_{k}})}/{\sqrt{t_i}^{\alpha_i}}$.

Additionally if Assumption \ref{quadratic} also holds, we have
\begin{equation*}
	\mathbb{E}\big\{\mathcal{R}^\tau(\hat{f}_\phi)-\mathcal{R}^\tau(f_0)\big\}\leq C\frac{\lambda_\tau\mathcal{B}\mathcal{S}\mathcal{D}\log(\mathcal{S})\log(n)}{n^{1-1/p}}+ 2c_\tau\big[\sum_{i\in J^c} C_i^*\lambda_i^* t_i^*(N_iL_i)^{-2\alpha_i^*/t_i}\big]^2,
\end{equation*}
where $c_\tau>0$ is a constant defined in
Lemma \ref{lemma4}
and $C>0$ is a constant not depending on $n,d,\tau,\mathcal{B},$ $\mathcal{S},$ $\mathcal{D},$ $C_i^*,$ $\lambda_i^*,$ $\alpha_i^*,$ $N_i$ or $L_i$.
\end{theorem}

\begin{remark}
In Theorem \ref{thm2}, the bounds for the excess risk are explicitly expressed in terms of
the network parameters $\mathcal{D}$ and $\mathcal{S}$ and the parameters  $N_i$ and $L_i$.
, which determine the width and the depth of the network as specified in (\ref{widtha}) and
(\ref{deptha}).
The dependence of the bounds on the dimensions of the functions $(d, t_j)$ and the H\"older constants
$(\alpha_j, \lambda_j)$  for the functions is also explicitly described. These constants are given and determined by the underlying model, so we cannot change them.
The constants $C$ and $c_{\tau}$ are independent of all the above parameters, in particular, they do not
depend on the dimensions $(d, t_j)$.
\end{remark}

Theorem \ref{thm2} gives a general expression of the upper bound for the excess risk. This bound clearly describes how the bounds depend on various parameters.
The parameters that can be changed or tuned are the network parameters given in terms of
$N_i$ and $L_i$. We note that the stochastic error term increases with $(N_i, L_i)$, while the approximation error term decreases with $(N_i, L_i)$. Thus we can select $(N_i, L_i)$ to balance these
two error terms, which lead to the best error bound. We will present an explicit expression of the risk bound in Corollary \ref{thm2d} below. First, we state a simpler bound assuming that all the component functions in the composition are  Lipschitz continuous with $\alpha_i=1, i=0, 1, \ldots, q.$

\begin{corollary}\label{thm2c}
Under model (\ref{model}),  suppose Assumptions \ref{structure} and \ref{moment} hold and  all $h_{ij}:D_{ij}\to\mathbb{R}$ in Theorem \ref{thm1} are Lipschitz continuous functions ($\alpha_i=1$ for $i=0,\ldots,q$) with Lipschitz constants $\lambda_i\ge0$.
 Given any $N,L\in\mathbb{N}^+$, for $i\in J^c$, we set the same shape for each subnetwork with $N_i=N\in\mathbb{N}^+$ and $L_i=L\in\mathbb{N}^+$, and for $j\in J$, we set the 3-layer subnetwork with width $(d_j,2d_j,d_{j+1})$ according to Lemma \ref{lemma8}.
Suppose the network parameters of the function class $\mathcal{F}_{\phi}$ are specified as in (\ref{widtha}) and (\ref{deptha}). Then,
for $2n \ge \text{Pdim}(\mathcal{F}_\phi)$,
the excess risk of the \textit{DQR} estimator $\hat{f}_\phi$ satisfies
\begin{equation*}
	\mathbb{E}\big\{\mathcal{R}^\tau(\hat{f}_\phi)-\mathcal{R}^\tau(f_0)\big\}\leq C\frac{\lambda_\tau\mathcal{B}\mathcal{S}\mathcal{D}\log(\mathcal{S})\log(n)}{n^{1-1/p}}+ 36\lambda_\tau\sum_{i\in J^c}  \Pi_{k=i+1}\sqrt{t_k}(N_iL_i)^{-2/t_i},
\end{equation*}
where $\lambda_\tau=\max\{\tau,1-\tau\}$ and $C>0$ is a constant independent of $n,d,\tau,\mathcal{B},\mathcal{S},\mathcal{D},N$ or $L$. Additionally if Assumption \ref{quadratic} also holds, we have
\begin{equation*}
	\mathbb{E}\big\{\mathcal{R}^\tau(\hat{f}_\phi)-\mathcal{R}^\tau(f_0)\big\}\leq C\frac{\lambda_\tau\mathcal{B}\mathcal{S}\mathcal{D}\log(\mathcal{S})\log(n)}{n^{1-1/p}}+ 648c_\tau\big[\sum_{i\in J^c}  \Pi_{k=i+1}\sqrt{t_k}(N_iL_i)^{-2/t_i}\big]^2,
\end{equation*}
where $c_\tau>0$ is a constant defined in Assumption \ref{quadratic} and $C>0$ is a constant independent of $n,d,\tau,\mathcal{B},\mathcal{S},\mathcal{D},N$ or $L$.

\end{corollary}

\begin{remark}
	The $\log(n)$ factor  in the stochastic error of the upper bound in Theorem \ref{thm2} and Corollary \ref{thm2c} is due to the truncation technique used in the proof. Power of log factors, $ (log n)^k$ for some $k\in\mathbb{N}^+$, are commonly seen in the results of related work, e.g.,  \cite{bauer2019deep,schmidt2020nonparametric}  and \cite{farrell2021deep}. By properly setting the network size $\mathcal{S}$ or depth $\mathcal{D}$ to have order $O(n^c/(\log n)^k)$ for some constant $c>0$ and $k\in\mathbb{N}^+$, the final convergence rate of the excess risk could be made optimal. However, this will make the selection of the network parameters more complicated.
Therefore, we will not do so in this paper.  The rate of convergence is (nearly) optimal up to a logarithmic factor $(\log n)^2$.
\end{remark}

We now present an explicit risk bound for the \textit{DQR} estimators with three sets of network parameters with different depth and width. All these three different specifications of the network parameters lead to the same risk bound.

\begin{corollary}
\label{thm2d}
	Under model (\ref{model}), suppose that Assumptions \ref{structure}-\ref{quadratic} hold, $\nu$ is absolutely continuous with respect to the Lebesgue measure, $ \Vert f_0\Vert_\infty\leq\mathcal{B}$ for some $\mathcal{B}\geq1$ and $2n \ge \text{Pdim}(\mathcal{F}_\phi)$. Let $(\alpha^*,t^*)=\arg\min_{(\alpha_i^*,t_i),i\in J^c} \{\alpha^*_i/t_i\}$, $\lambda^*=\max_{i=0,\ldots,q} \lambda_i^*$ and $d^*=\max_{i=0,\ldots,q}t_i^*$, where $\alpha_i^*,\lambda_i^*$ and $t_i^*$ are defined in Theorem \ref*{thm2}. Suppose the network parameters of the function class $\mathcal{F}_{\phi}$ are specified as follows:
	\begin{itemize}
		\item [1.] (Deep and fixed width MLP) Let $N_i=1$ and $L_i=\lfloor n^{(1-1/p)t^*/(4\alpha^*+2t^*)}\rfloor$. The corresponding width, depth and size of the networks satisfy:
		\begin{align*}
			&\mathcal{W}_1=\max_{i=0,\ldots,q}d_i\max\{7t_i,20\},\\
			&\mathcal{D}_1=(12\lfloor n^{(1-1/p)t^*/(4\alpha^*+2t^*)}\rfloor+15)\vert J^c\vert+2\vert J\vert,\\
			&\mathcal{S}_1\leq \mathcal{W}_1^2\mathcal{D}_1\leq \max_{i=0,\ldots,q} (20d_it_i)^2\times29q\lfloor n^{(1-1/p)t^*/(4\alpha^*+2t^*)}\rfloor.
		\end{align*}
	\item [2.] (Deep and wide MLP) Let $N_i=\lfloor n^{(1-1/p)t^*/(8\alpha^*+4t^*)}\rfloor$ and $L_i=\lfloor n^{(1-1/p)t^*/(8\alpha^*+4t^*)}\rfloor$. The corresponding width, depth and size of the networks satisfy:
	\begin{align*}
	&\mathcal{W}_2=\max_{i=0,\ldots,q}d_{i}\max\{4t_i\lfloor \lfloor n^{(1-1/p)t^*/(8\alpha^*+4t^*)}\rfloor^{1/t_i}\rfloor+3t_i,12\lfloor n^{(1-1/p)t^*/(8\alpha^*+4t^*)}\rfloor+8\},\\
	&\mathcal{D}_2=(12\lfloor n^{(1-1/p)t^*/(8\alpha^*+4t^*)}\rfloor+15)\vert J^c\vert+2\vert J\vert,\\
	&\mathcal{S}_2\leq \mathcal{W}_2^2\mathcal{D}_2\leq \max_{i=0,\ldots,q} (20d_it_i)^2\times 29q\lfloor n^{(1-1/p)t^*/(4\alpha^*+2t^*)}\rfloor^{3/2}.
	\end{align*}

		\item [3.] (Fixed depth and wide MLP) Let $N_i=\lfloor n^{(1-1/p)t^*/(4\alpha^*+2t^*)}\rfloor$ and $L_i=1$. The corresponding width, depth and size of the networks satisfy:
	\begin{align*}
		&\mathcal{W}_3=\max_{i=0,\ldots,q}d_{i}\max\{4t_i\lfloor \lfloor n^{(1-1/p)t^*/(4\alpha^*+2t^*)}\rfloor^{1/t_i}\rfloor+3t_i,12\lfloor n^{(1-1/p)t^*/(4\alpha^*+2t^*)}\rfloor+8\},\\
		&\mathcal{D}_3=27\vert J^c\vert+2\vert J\vert,\\
		&\mathcal{S}_3\leq \mathcal{W}_3^2\mathcal{D}_3\leq \max_{i=0,\ldots,q} (20d_it_i)^2\times 29q \lfloor n^{(1-1/p)t^*/(4\alpha^*+2t^*)}\rfloor^2.
	\end{align*}
	\end{itemize}	

	Then,  the  excess risk satisfies
	\begin{equation}
		\label{erb}
\mathbb{E}\big\{\mathcal{R}^\tau(\hat{f}_\phi)-\mathcal{R}^\tau(f_0)\big\}\leq C_0 C_{d, d^*} (\log n)^2 n^{-\left(1-\frac{1}{p}\right)\frac{2\alpha^*}{2\alpha^*+t^*}},
	\end{equation}
where $C_{d, d^*}= (d^*)^2(\max_{i=0,\ldots,q}d_it_i)^2\log(\max_{i=0,\ldots,q}d_it_i)$, $C_0=c\lambda_\tau c_\tau\mathcal{B} q^2\log(q)(\lambda^*)^2$. Here $c$ is a universal constant not depending on any parameters.
\end{corollary}

In Corollary \ref{thm2d}, three sets of different network parameters lead to the same risk bound. Therefore, generally the choice of network parameters is not unique to achieve a desired risk bound.
Although the three sets of network parameters given in  Corollary \ref{thm2d} yield the same risk bound,
the sizes of the networks are different. As can be seen from the expressions of the network sizes $\mathcal{S}_1$,
$\mathcal{S}_2$ and $\mathcal{S}_3$, we have, on the logarithmic scale,
\[
\log \mathcal{S}_1 : \log \mathcal{S}_2 : \log \mathcal{S}_3 = 1: \frac{3}{2}: 2.
\]
Therefore, the deep and fixed width network in the  first network specification with width $\mathcal{W}_1$ and depth $\mathcal{D}_1$ is the most efficient design among the three network structures in the sense that it has the smallest network size.
Corollary \ref{thm2d} shows that deep networks have advantages over shallow ones in the sense that deep networks achieve the same risk bound with a smaller network size.
More detailed discussions on the relationship between convergence rate and  network structure can be found in \citet{jiao2021deep}.

\subsection{Mean integrated squared error} 
\label{sec4.2}
The empirical risk minimization quantile estimator typically
 results in an estimator $\hat{f}_n$ for which its risk $\mathcal{R}^\tau(\hat{f}_n)$ is close to optimal risk $\mathcal{R}^\tau(f_0)$ in expectation or with high probability. However, small excess risk in general only implies in a weak sense that the ERM $\hat{f}_n$ is
close to $f_0$ (Remark 3.18, \cite{steinwart2007compare}). Hence, in this subsection, we bridge the gap between the excess risk and {\color{black} the mean integrated squared error (MISE)} of the estimated conditional quantile function. 
To this end, we need the following condition on the conditional distribution of $Y$ given $X$.

\begin{assumption}\label{calib}
	There exist constants $\gamma>0$ and $\kappa>0$ such that for any $\vert\delta\vert\leq \gamma$,
	\begin{equation*}
		\big\vert P_{Y|X}(f_0(x)+\delta\mid x)-P_{Y|X}(f_0(x)\mid x)\big\vert\geq \kappa \vert \delta\vert,
	\end{equation*}
for all $x\in\mathcal{X}$ up to a $\nu$-negligible set, where $P_{Y|X}(\cdot|x)$ denotes the conditional distribution function of $Y$ given $X=x$.
\end{assumption}

\begin{remark}
A similar condition is assumed by \cite{padilla2021adaptive} in studying nonparametric quantile trend filtering. This condition is weaker than Condition 2.1 in \cite{he1994convergence} and condition D.1 in \cite{belloni2011}, which require the conditional density of $Y$ given $X=x$ to be bounded below near its $\tau$-th quantile.
\end{remark}

Under Assumption \ref{calib}, the self-calibration condition can be established as stated below.
This will lead to a bound on the MISE of the estimated quantile function based on a bound for the excess risk.

{\color{black}
\begin{lemma}[Self-calibration]	\label{lemma9}
	Suppose  that Assumption \ref{moment} (i) and Assumption \ref{calib} hold. For any $f:\mathcal{X}\to\mathbb{R}$, denote $\Delta^2 (f,f_0)=\mathbb{E}\big[\min\{\vert f(X)-f_0(X)\vert^2,\vert f(X)-f_0(X)\vert\}\big]$ where $\kappa$ and $\gamma>0$ are defined in Assumption \ref{calib}. Then we have
	\begin{equation*}
		\Delta^2 (f,f_0)\leq c_{\kappa,\gamma}\big\{\mathcal{R}^\tau(f)-\mathcal{R}^\tau(f_0)\big\},
	\end{equation*}
for any $f:\mathcal{X}\to\mathbb{R}$ where $c_{\kappa,\gamma}=\max\{2/\kappa,4/(\kappa\gamma)\}$. More exactly, for $f:\mathcal{X}\to\mathbb{R}$ satisfying $\vert f(x)-f_0(x)\vert\leq\gamma$ for $x\in\mathcal{X}$ up to a $\nu$-negligible set, we have
\begin{equation*}
	\Vert f-f_0\Vert^2_{L^2(\nu)}\leq \frac{2}{\kappa}\big\{\mathcal{R}^\tau(f)-\mathcal{R}^\tau(f_0)\big\},
\end{equation*}
otherwise we have
\begin{equation*}
	\Vert f-f_0\Vert_{L^1(\nu)}\leq \frac{4}{\kappa\gamma}\big\{\mathcal{R}^\tau(f)-\mathcal{R}^\tau(f_0)\big\}.
\end{equation*}
\end{lemma}
}

\begin{remark}
	Similar self-calibration conditions can be found in \cite{christmann2007svms,steinwart2011estimating,lv2018oracle} and \cite{padilla2020quantile}. A general result is obtained in \cite{steinwart2011estimating} under the so-called $\tau$-quantile of $t$-average type assumption on the joint distribution $P$, where $\Vert f-f_0\Vert_{L^r(\nu)}$ is upper bounded by the $q$-th root of excess risk $\mathcal{R}^\tau(f)-\mathcal{R}^\tau(f_0)$ for $t\in(0,\infty]$, $q\in[1,\infty)$ and $r=tq/(t+1)$. However, those assumptions on the joint distribution $P$ generally require that the conditional distribution of $Y$ given $X$ is bounded, which may not be applicable to  models with  heavy-tailed response as in our setting, see, e.g.,  Assumption \ref{moment}.
\end{remark}

\begin{theorem}[Non-asymptotic bound for {\color{black} mean integrated squared error}]
\label{thm3}
Under model (\ref{model}), suppose that Assumptions \ref{structure}, \ref{moment} and \ref{calib} hold, $\nu$ is absolutely continuous with respect to the Lebesgue measure, and $ \Vert f_0\Vert_\infty\leq\mathcal{B}$ for some $\mathcal{B}\geq1$.
Then, given any $N_i,L_i\in\mathbb{N}^+, i\in J^c$, for the function class of ReLU multi-layer perceptrons $\mathcal{F}_{\phi}=\mathcal{F}_{\mathcal{D},\mathcal{W},\mathcal{U},\mathcal{S},\mathcal{B}}$
with width no larger than $\mathcal{W}=\max_{i=0,\ldots,q}d_{i}\max\{4t_i\lfloor N_i^{1/t_i}\rfloor+3t_i,12N_i+8\}$ and depth $\mathcal{D}=\sum_{i\in J^c}(12L_i+15)+2\vert J\vert$,
for $2n \ge \text{Pdim}(\mathcal{F}_\phi)$,
{\color{black}
the MISE 
 of the \textit{DQR} estimator $\hat{f}_\phi$ satisfies
\begin{equation*}
	\mathbb{E} \big\{\Delta^2 (\hat{f}_\phi,f_0)\big\} \leq c_{\kappa,\gamma}\lambda_\tau\Big[ C\frac{\mathcal{B}\mathcal{S}\mathcal{D}\log(\mathcal{S})\log(n)}{n^{1-1/p}}+ 2\sum_{i\in J^c} C_i^*\lambda_i^* t_i^*(N_iL_i)^{-2\alpha_i^*/t_i}\Big],
\end{equation*}
where $c_{\kappa,\gamma}=\max\{4/(\kappa\gamma),2/\kappa\}$ and $\Delta^2(\cdot,\cdot)$ are defined in Lemma \ref{lemma9},
}
$\lambda_\tau=\max\{\tau,1-\tau\}$ and $C>0$ is a constant not depending on $n,d,\tau,\mathcal{B},\mathcal{S},\mathcal{D},C_i^*,\lambda_i^*,\alpha_i^*,N_i$ or $L_i$, and $C_i^*=18^{\Pi_{j=i+1}^{q}\alpha_{j}}$, $\lambda_i^*=\Pi_{j=i}^{q}\lambda_{j}^{\Pi_{k=j+1}^{q} \alpha_{k}}$, $\alpha_i^*=\Pi_{j=i}^q \alpha_j$ and $t_i^*={(\Pi_{j=i}^{q}\sqrt{t_{j}}^{\Pi_{k=j}^{q} \alpha_{k}})}/{\sqrt{t_i}^{\alpha_i}}$. Additionally if Assumption \ref{quadratic} also holds, we have
\begin{equation*}
		\mathbb{E}\Vert \hat{f}_\phi-f_0\Vert_{L^*(\nu)}\leq c_{\kappa,\gamma}\Big[ C\frac{\lambda_\tau\mathcal{B}\mathcal{S}\mathcal{D}\log(\mathcal{S})\log(n)}{n^{1-1/p}}+ 2c_\tau\big\{\sum_{i\in J^c} C_i^*\lambda_i^* t_i^*(N_iL_i)^{-2\alpha_i^*/t_i}\big\}^2\Big],
\end{equation*}
where $c_\tau>0$ is a constant defined in Assumption \ref{quadratic} and $C>0$ is a constant independent of $n,d,\tau,\mathcal{B},\mathcal{S},\mathcal{D},C_i^*,\lambda_i^*,\alpha_i^*,N_i$ or $L_i$.
\end{theorem}

Similar to Corollary \ref{thm2d}, we have the following corollary for {\color{black}  the MISE} of the \textit{DQR} estimator.

\begin{corollary}
\label{thm3c}
	Under model (\ref{model}), suppose that Assumptions \ref{structure}-\ref{quadratic} hold, $\nu$ is absolutely continuous with respect to the Lebesgue measure, $ \Vert f_0\Vert_\infty\leq\mathcal{B}$ for some $\mathcal{B}\geq1$ and $2n \ge \text{Pdim}(\mathcal{F}_\phi)$. Let $(\alpha^*,t^*)=\arg\min_{(\alpha_i^*,t_i),i\in J^c} \{\alpha^*_i/t_i\}$, $\lambda^*=\max_{i=0,\ldots,q} \lambda_i^*$ and $d^*=\max_{i=0,\ldots,q}t_i^*$, where $\alpha_i^*,\lambda_i^*$ and $t_i^*$ are defined  in Theorem \ref*{thm2}. Suppose that the network parameters of the function class $\mathcal{F}_{\phi}$ are specified as follows:
	\begin{itemize}
		\item [1.] (Deep and fixed width MLP) Let $N_i=1$ and $L_i=\lfloor n^{(1-1/p)t^*/(4\alpha^*+2t^*)}\rfloor$. The corresponding width, depth and size of the networks satisfy:
		\begin{align*}
			&\mathcal{W}_1=\max_{i=0,\ldots,q}d_i\max\{7t_i,20\},\\
			&\mathcal{D}_1=(12\lfloor n^{(1-1/p)t^*/(4\alpha^*+2t^*)}\rfloor+15)\vert J^c\vert+2\vert J\vert,\\
			&\mathcal{S}_1\leq \mathcal{W}_1^2\mathcal{D}_1\leq \max_{i=0,\ldots,q} (20d_it_i)^2\times29q\lfloor n^{(1-1/p)t^*/(4\alpha^*+2t^*)}\rfloor.
		\end{align*}
	\item [2.] (Deep and wide MLP) Let $N_i=\lfloor n^{(1-1/p)t^*/(8\alpha^*+4t^*)}\rfloor$ and $L_i=\lfloor n^{(1-1/p)t^*/(8\alpha^*+4t^*)}\rfloor$. The corresponding width, depth and size of the networks satisfy:
	\begin{align*}
	&\mathcal{W}_2=\max_{i=0,\ldots,q}d_{i}\max\{4t_i\lfloor \lfloor n^{(1-1/p)t^*/(8\alpha^*+4t^*)}\rfloor^{1/t_i}\rfloor+3t_i,12\lfloor n^{(1-1/p)t^*/(8\alpha^*+4t^*)}\rfloor+8\},\\
	&\mathcal{D}_2=(12\lfloor n^{(1-1/p)t^*/(8\alpha^*+4t^*)}\rfloor+15)\vert J^c\vert+2\vert J\vert,\\
	&\mathcal{S}_2\leq \mathcal{W}_2^2\mathcal{D}_2\leq \max_{i=0,\ldots,q} (20d_it_i)^2\times 29q\lfloor n^{(1-1/p)t^*/(4\alpha^*+2t^*)}\rfloor^{3/2}.
	\end{align*}

		\item [3.] (Fixed depth and wide MLP) Let $N_i=\lfloor n^{(1-1/p)t^*/(4\alpha^*+2t^*)}\rfloor$ and $L_i=1$. The corresponding width, depth and size of the networks satisfy:
	\begin{align*}
		&\mathcal{W}_3=\max_{i=0,\ldots,q}d_{i}\max\{4t_i\lfloor \lfloor n^{(1-1/p)t^*/(4\alpha^*+2t^*)}\rfloor^{1/t_i}\rfloor+3t_i,12\lfloor n^{(1-1/p)t^*/(4\alpha^*+2t^*)}\rfloor+8\},\\
		&\mathcal{D}_3=27\vert J^c\vert+2\vert J\vert,\\
		&\mathcal{S}_3\leq \mathcal{W}_3^2\mathcal{D}_3\leq \max_{i=0,\ldots,q} (20d_it_i)^2\times 29q \lfloor n^{(1-1/p)t^*/(4\alpha^*+2t^*)}\rfloor^2.
	\end{align*}
	\end{itemize}	

{\color{black}
Then,  we have
	\begin{equation}
		\label{miseb}
		 \mathbb{E} \big\{\Delta^2 (\hat{f}_\phi,f_0)\big\} \leq A_0 A_{d, d^*}
(\log n)^2 n^{-\left(1-\frac{1}{p}\right)\frac{2\alpha^*}{2\alpha^*+t^*}},
	\end{equation}
}
	where $A_{d, d^*}= (d^*)^2(\max_{i=0,\ldots,q}d_it_i)^2\log(\max_{i=0,\ldots,q}d_it_i)$, $A_0=cc_{\kappa,\gamma}\lambda_\tau c_\tau\mathcal{B} q^2\log(q)(\lambda^*)^2$,
with  $c$ a universal constant independent of any parameters.
\end{corollary}
We note that, according to Corollary \ref{thm3c},  the same comments about the relationship between the network sizes and the risk bound
following Corollary \ref{thm2d} apply to the relationship between the network size and {\color{black} the MISE }of the \textit{DQR} estimator.

\section{Examples}
\label{compositionR}

In this section, we specialize the general results in Theorems \ref{thm2} and \ref{thm3} and Corollaries \ref{thm2d} and \ref{thm3c}
to several important models widely used in statistics. We explicitly describe how the prefactor depends on the ambient dimension and the intrinsic dimension of the model.  We present the results with $\mathcal{F}_n$ consisting of  deep and fixed-width network functions in constructing the \textit{DQR} estimators,
as such networks are more efficient in the sense that they require a smaller network size to achieve the optimal convergence rate compared with other shaped networks, see Corollaries \ref{thm2d} and \ref{thm3c}.

We note that, in computing the \textit{DQR} estimator as defined in (\ref{erm}), we do not use the information about the specific structure of the models considered below. This is different from the methods in literature that are designed based on the model structure. For example,
the backfitting algorithm \citep{breiman1985additive} for fitting the additive conditional mean model (\ref{additiveM}) with the least squares loss  specifically use the
additive structure of the model. See also \citet{cds1997} and \citet{horowitz2005qadditive} for methods that estimate a conditional quantile model based on the additive structure assumption.
In the single index conditional mean model, \citet{hjs2001sim} described a method for estimating the index regression coefficient $\theta$.  With their method and regularity conditions, the difference between the distribution of  their estimator $\hat{\theta}_{\text{HJS}}$ and a mean-zero multivariate normal distribution converges to zero at a rate that does not depend on the dimension $d$ of the predictor.  This suggests that a kernel estimator of the index function using $\hat{\theta}_{\text{HJS}}$ in place of  $\theta$ has the usual one-dimensional rate of convergence that does not depend on the dimension $d$.   \citet{khan2001} also developed a two-stage method for estimating a
model satisfying a monotonicity condition on the conditional quantile function of the response variable.
However, these estimators heavily depend on the single index model assumption, they may not be consistent if this model assumption is not satisfied.

Let  $c_{\kappa,\gamma}=\max\{4/(\kappa\gamma),2/\kappa\}$ in all the examples below,
where $\kappa$ and $\gamma$ are the constants defined in Assumption \ref{calib}.

\subsection{Single index model}
A  popular semiparametric model in statistics and econometrics for mitigating the curse of dimensionality is the single index model
\begin{equation}
\label{siM} f_0(x)=g(\theta^\top x),  \quad {x\in\mathbb{R}^d},
\end{equation}
where $g:\mathbb{R}\to\mathbb{R}$ is a univariate function and $\theta\in\mathbb{R}^d$ is a $d$-dimensional vector.
Such $f_0$ can be written as a composition of functions
$$f_0=h_1\circ h_0,$$
where $h_0(x)=\theta^\top x$ is a linear transformation and $h_1(x)=g(x)$. Then $d_0=t_0=d, d_1=t_1=1$ and $d_2=1$ according to the definition in Assumption \ref{structure}. Suppose  that Assumptions \ref{structure}-\ref{moment} and the conditions in Theorem \ref{thm2} are satisfied, where $g$ or $h_1$ is H\"older continuous with order $\alpha_1$ and constant $\lambda_1$. Then by Theorem \ref{thm2}, given any $N,L\in\mathbb{N}^+$, for the function class of ReLU multi-layer perceptrons $\mathcal{F}_{\phi}=\mathcal{F}_{\mathcal{D},\mathcal{W},\mathcal{U},\mathcal{S},\mathcal{B}}$
with width $\mathcal{W}=\max\{12N+8,2d\}$ and depth $\mathcal{D}=12L+17$,
for $2n \ge \text{Pdim}(\mathcal{F}_\phi)$, the excess risk of the \textit{DQR} estimator $\hat{f}_\phi$ satisfies
\begin{equation*}
	\mathbb{E}\big\{\mathcal{R}^\tau(\hat{f}_\phi)-\mathcal{R}^\tau(f_0)\big\}\leq C\frac{\lambda_\tau\mathcal{B}\mathcal{S}\mathcal{D}\log(\mathcal{S})\log(n)}{n^{1-1/p}}+ 36\lambda_\tau \lambda_1(NL)^{-2\alpha_1},
\end{equation*}
where $C>0$ is a constant not depending on $n,d,\tau,\mathcal{B},\mathcal{S},\mathcal{D},\lambda_1,\alpha_1,N,L$ and $\lambda_\tau=\max\{\tau,1-\tau\}$. If we choose $N=1$ and $L=\lfloor n^{(1-1/p)/(2\alpha_1+2)}\rfloor$, then $\mathcal{S}\leq (20^2+20)\times(12L+15)+d\times(2d)+2d\leq 8\times20\times21\times27\times d^2\times \lfloor n^{(1-1/p)/(2\alpha_1+2)}\rfloor$ and
$$\mathbb{E}\big\{\mathcal{R}^\tau(\hat{f}_\phi)-\mathcal{R}^\tau(f_0)\big\}\leq C\mathcal{B}\times d^2\log(d) \times (\log n)^2 n^{-(1-1/p)\alpha_1/(\alpha_1+1)},$$
where $C>0$ is a constant independent of $n,d,\mathcal{B}$ and $\alpha_1$.

If Assumption \ref{quadratic} also holds, we have
\begin{equation*}
	\mathbb{E}\big\{\mathcal{R}^\tau(\hat{f}_\phi)-\mathcal{R}^\tau(f_0)\big\}\leq C\frac{\lambda_\tau\mathcal{B}\mathcal{S}\mathcal{D}\log(\mathcal{S})\log(n)}{n^{1-1/p}}+ 648c_\tau \lambda_1^2(NL)^{-4\alpha_1},
\end{equation*}
where $c_\tau>0$ is a constant defined in Lemma \ref{lemma4}. Alternatively, if we choose $N=1$ and $L=\lfloor n^{(1-1/p)/(4\alpha_1+2)}\rfloor$, then
$$\mathbb{E}\big\{\mathcal{R}^\tau(\hat{f}_\phi)-\mathcal{R}^\tau(f_0)\big\}\leq C_0 \mathcal{B}\times d^2\log(d) \times (\log n)^2
n^{-\left(1-\frac{1}{p}\right) \frac{2 \alpha_1}{2\alpha_1+1}},
$$
where $C_0 >0$ is a constant not depending on $n,d,\mathcal{B}$ and $\alpha_1$.

Additionally,  if Assumption \ref{calib} holds, it follows from Theorem \ref{thm3} that
{\color{black}
$$\mathbb{E} \big\{\Delta^2 (\hat{f}_\phi,f_0)\big\} \leq  c_{\kappa,\gamma}C_0 \mathcal{B}\times d^2\log(d) \times (\log n)^2
n^{-\left(1-\frac{1}{p}\right) \frac{2 \alpha_1}{2\alpha_1+1}}.
$$
}

\subsection{Additive model}
A well-known structured model is the additive model \citep{stone1985, stone1986dimensionality, hastie1990generalized}
\begin{equation}
\label{additiveM}
f_0(x_1, \ldots,x_d)=f_{0,1}(x_1)+\cdots +f_{0,d}(x_d),  \quad {x=(x_1, \ldots, x_d)^\top \in\mathbb{R}^d},
\end{equation}
where $f_{0,j}:\mathbb{R}\to\mathbb{R}$,  $j=1,\ldots,d$,  are univariate functions.
This model is a direct nonparametric extension of the linear model. It has certain appealing computational and theoretical properties. In particular, it  can be estimated with the optimal rate of convergence of the univariate nonparametric regression \citep{stone1986dimensionality}.
The additive function  $f_0$ can be written as a simple composition of functions
$$f_0=h_1\circ h_0,$$
where $h_0(x)=(f_{0,1}(x),\ldots,f_{0,d}(x))^\top$ and $h_1(x)=\sum_{i=1}^d x_i$ where $x=(x_1,\ldots, x_d)^\top \in\mathbb{R}^d$. In this case, $d_0=d,t_0=1,d_1=t_1=d$ and $d_2=1$.  Suppose  that Assumption \ref{structure}-\ref{moment} and those conditions in Theorem \ref{thm2} are satisfied, where $f_{0,i}$ is H\"older continuous with order $\alpha_0$ and constant $\lambda_0$ for $i=1,\ldots,d$. Then by Theorem \ref{thm2}, given any $N,L\in\mathbb{N}^+$, for the function class of ReLU multi-layer perceptrons $\mathcal{F}_{\phi}=\mathcal{F}_{\mathcal{D},\mathcal{W},\mathcal{U},\mathcal{S},\mathcal{B}}$
with width $\mathcal{W}=(12N+8)d$ and depth $\mathcal{D}=12L+17$,
for $2n \ge \text{Pdim}(\mathcal{F}_\phi)$, the excess risk of the \textit{DQR} estimator $\hat{f}_\phi$ satisfies
$$\mathbb{E}\big\{\mathcal{R}^\tau(\hat{f}_\phi)-\mathcal{R}^\tau(f_0)\big\}\leq C\frac{\lambda_\tau\mathcal{B}\mathcal{S}\mathcal{D}\log(\mathcal{S})\log(n)}{n^{1-1/p}}+ 36\lambda_\tau \lambda_0\sqrt{d}(NL)^{-2\alpha_0},$$
where $C>0$ is a constant that does not depend on $n,d,\tau,\mathcal{B},\mathcal{S},\mathcal{D},\lambda_0,\alpha_0,N,L$ and $\lambda_\tau=\max\{\tau,1-\tau\}$. If we choose $N=1$ and $L=\lfloor n^{(1-1/p)/(2\alpha_0+2)}\rfloor$, then $\mathcal{S}\leq \{(20d)^2+20d\}\times(12L+15)+d\times(2d)+2d\leq 20\times21\times27\times d^2\times \lfloor n^{(1-1/p)/(2\alpha_0+2)}\rfloor$ and
$$\mathbb{E}\big\{\mathcal{R}^\tau(\hat{f}_\phi)-\mathcal{R}^\tau(f_0)\big\}\leq C \mathcal{B}\times d^2\log(d) \times (\log n)^2
n^{-\left(1-\frac{1}{p}\right) \frac{\alpha_0}{\alpha_0+1}},
$$
where $C>0$ is a constant not depending on $n,d,\mathcal{B}$ and $\alpha_0$.

If Assumption \ref{quadratic} also holds, we have
\begin{equation*}
\mathbb{E}\big\{\mathcal{R}^\tau(\hat{f}_\phi)-\mathcal{R}^\tau(f_0)\big\}\leq C\frac{\lambda_\tau\mathcal{B}\mathcal{S}\mathcal{D}\log(\mathcal{S})\log(n)}{n^{1-1/p}}+ 648c_\tau \lambda^2_0d(NL)^{-4\alpha_0},
\end{equation*}
where $c_\tau>0$ is a constant defined in Lemma \ref{lemma4}. Alternatively, if we choose $N=1$ and $L=\lfloor n^{(1-1/p)/(4\alpha_0+2)}\rfloor$, then
$$\mathbb{E}\big\{\mathcal{R}^\tau(\hat{f}_\phi)-\mathcal{R}^\tau(f_0)\big\}\leq C_0 \mathcal{B}\times d^2\log(d) \times (\log n)^2
n^{-\left(1-\frac{1}{p}\right) \frac{2\alpha_0}{2\alpha_0+1}},
$$
where $C_0>0$ is a constant not depending on $n,d,\mathcal{B}$ and $\alpha_0$.

 Additionally,  if Assumption \ref{calib} holds, it follows from Theorem \ref{thm3} that
 {\color{black}
$$ \mathbb{E} \big\{\Delta^2 (\hat{f}_\phi,f_0)\big\} \leq  c_{\kappa,\gamma}C_0\mathcal{B}\times d^2\log(d) \times (\log n)^2
n^{-\left(1-\frac{1}{p}\right) \frac{2\alpha_0}{2\alpha_0+1}}.
$$
}

\subsection{Additive model with an unknown link function}
The additive model with an
unknown link function is
\begin{equation}
\label{additiveUM}
f_0(x)= f_1 (f_{0,1}(x_1)+\cdots  + f_{0,d}(x_d)), \ x \in \mathbb{R}^d,
\end{equation}
where $f_1,f_{0,1},\ldots,f_{0,d}$ are univariate real-functions. Such $f_0$ has one more hierarchy than that of Additive model, which can be written as
$$f_0=h_2\circ h_1\circ h_0,$$
where $h_0(x)=(f_{0,1}(x),\ldots,f_{0,d}(x))^\top$, $h_1(x)=\sum_{i=1}^d x_i$ and $h_2(x)=f_1(x)$ where $x=(x_1,\ldots, x_d)^\top \in\mathbb{R}^d$. In this case, $d_0=d,t_0=1,d_1=t_1=d, d_2=t_2=1$ and $d_3=1$. Suppose that Assumptions \ref{structure}-\ref{moment} and those conditions in Theorem \ref{thm2} hold, where $f_{0,i}$ is H\"older continuous with order $\alpha_0$ and constant $\lambda_0$ for $i=1,\ldots,d$ and $f_1$  is H\"older continuous with order $\alpha_2$ and constant $\lambda_2$. By Theorem \ref{thm2}, given any $N,L\in\mathbb{N}^+$, for the function class of ReLU multi-layer perceptrons $\mathcal{F}_{\phi}=\mathcal{F}_{\mathcal{D},\mathcal{W},\mathcal{U},\mathcal{S},\mathcal{B}}$
with width $\mathcal{W}=(12N+8)d$ and depth $\mathcal{D}=24L+32$,
for $2n \ge \text{Pdim}(\mathcal{F}_\phi)$, the excess risk of the \textit{DQR} estimator $\hat{f}_\phi$ satisfies
\begin{align*}
\mathbb{E}\big\{\mathcal{R}^\tau(\hat{f}_\phi)-\mathcal{R}^\tau(f_0)\big\}\leq
& \
C\frac{\lambda_\tau\mathcal{B}\mathcal{S}\mathcal{D}\log(\mathcal{S})\log(n)}{n^{1-1/p}}\\
&+2\lambda_\tau\{18^{\alpha_2}\lambda_0^{\alpha_2}d^{\alpha_2/2}(NL)^{-2\alpha_0\alpha_2}
+18\lambda_2(NL)^{-2\alpha_2}\},
\end{align*}
where $C>0$ is a constant that does not depend on $n,d,\tau,\mathcal{B},\mathcal{S},\mathcal{D},\lambda_0,\lambda_2,\alpha_2,N,L$ and $\lambda_\tau=\max\{\tau,1-\tau\}$. If we choose $N=1$ and $L=\lfloor n^{(1-1/p)/(2\alpha_0\alpha_2+2)}\rfloor$, then $\mathcal{S}\leq \{(20d)^2+20d+20^2+20\}\times(12L+15)+d\times(2d)+2d\leq 2\times20\times21\times27\times d^2\times \lfloor n^{(1-1/p)/(2\alpha_0\alpha_2+2)}\rfloor$ and
$$\mathbb{E}\big\{\mathcal{R}^\tau(\hat{f}_\phi)-\mathcal{R}^\tau(f_0)\big\}\leq C\mathcal{B}\times d^2\log(d) \times (\log n)^2
n^{-\left(1-\frac{1}{p}\right) \frac{2\alpha_0\alpha_2}{\alpha_0\alpha_2+1}},
$$
where $C>0$ is a constant not depending on $n,d,\mathcal{B}$ and $\alpha_0,\alpha_2$.

Additionally,  if Assumption \ref{quadratic} holds, we have
\begin{align*}
	\mathbb{E}\big\{\mathcal{R}^\tau(\hat{f}_\phi)-\mathcal{R}^\tau(f_0)\big\}\leq&  C\frac{\lambda_\tau\mathcal{B}\mathcal{S}\mathcal{D}\log(\mathcal{S})\log(n)}{n^{1-1/p}}\\ &+2c_\tau\{18^{\alpha_2}\lambda_0^{\alpha_2}d^{\alpha_2/2}(NL)^{-2\alpha_0\alpha_2}+18\lambda_2(NL)^{-2\alpha_2}\}^2,
\end{align*}
where $c_\tau>0$ is a constant defined in Lemma \ref{lemma4}.
Alternatively, if we choose $N=1$ and $L=\lfloor n^{(1-1/p)/(4\alpha_0\alpha_2+2)}\rfloor$,
 then
$$\mathbb{E}\big\{\mathcal{R}^\tau(\hat{f}_\phi)-\mathcal{R}^\tau(f_0)\big\}\leq C_0\mathcal{B}\times d^2\log(d) \times (\log n)^2
n^{-\left(1-\frac{1}{p}\right) \frac{2\alpha_0\alpha_2}{\alpha_0\alpha_2+1}},
$$
where $C_0>0$ is a constant not depending on $n,d,\mathcal{B},\alpha_0$ and $\alpha_2$.

Moreover, if Assumption \ref{calib} holds,  Theorem \ref{thm3} implies that
{\color{black}
$$\mathbb{E} \big\{\Delta^2 (\hat{f}_\phi,f_0)\big\} \leq c_{\kappa,\gamma}C_0 \mathcal{B}\times d^2\log(d) \times (\log n)^2
n^{-\left(1-\frac{1}{p}\right) \frac{2\alpha_0\alpha_2}{\alpha_0\alpha_2+1}}.
$$
}
\subsection{Interaction model}

The additive model was also generalized to an interaction model  \citep{stone1994use}
\begin{equation}
\label{interactionM}
f_0(x)=\sum_{I\subseteq\{1,\ldots,d\},\vert I\vert=d^*}f_I(x_I),  \quad x=(x_1,\ldots,x_d)^\top \in\mathbb{R}^d,
\end{equation}
where $d^*\in\{1,\ldots,d\}$, $I=\{i_1,\ldots,i_{d^*}\}$,  $1\leq i_1<\ldots<i_{d^*}\leq d$, $x_I=(x_{i_1}, \ldots, x_{i_{d^*}})$ and all $f_I$ are H\"older continuous $d^*$-variate functions with order $\alpha_0$ and constant $\lambda_0$ defined on $\mathbb{R}^{\vert I\vert}$. Let $\mathcal{I}$ be the collection of index set $I$ in the summation, and let $K=\vert \mathcal{I}\vert$ be the cardinality of $\mathcal{I}$. For such $f_0$, in our notation, it can be written as a composition of two functions:
$$f_0=h_1\circ h_0,$$
where $h_0(x)=(f_1(x),\ldots,f_K(x))^\top$ and $h_1(x)=\sum_{i=1}^K x_i$ for $x=(x_1,\ldots,x_K)^\top\in\mathbb{R}^K$. Here $d_0=d,t_0=d^*,d_1=t_1=K$ and $d_2=1$. Suppose that Assumptions \ref{structure}-\ref{moment} and the conditions in Theorem \ref{thm2} are satisfied. Then by Theorem \ref{thm2}, given any $N,L\in\mathbb{N}^+$, for the function class of ReLU multi-layer perceptrons $\mathcal{F}_{\phi}=\mathcal{F}_{\mathcal{D},\mathcal{W},\mathcal{U},\mathcal{S},\mathcal{B}}$
with width $\mathcal{W}=d\max\{4d^*\lfloor N^{1/d^*}\rfloor+3d^*,12N+8\}$ and depth $\mathcal{D}=12L+17$,
for $2n \ge \text{Pdim}(\mathcal{F}_\phi)$, the excess risk of the \textit{DQR} estimator $\hat{f}_\phi$ satisfies
\begin{equation*}
\mathbb{E}\big\{\mathcal{R}^\tau(\hat{f}_\phi)-\mathcal{R}^\tau(f_0)\big\}\leq C\frac{\lambda_\tau\mathcal{B}\mathcal{S}\mathcal{D}\log(\mathcal{S})\log(n)}{n^{1-1/p}}+ 36\lambda_\tau \lambda_0\sqrt{K}(NL)^{-2\alpha_0},
\end{equation*}
where $C>0$ is a constant  not depending on $n,d,\tau,\mathcal{B},\mathcal{S},\mathcal{D},\lambda_0,\alpha_0,N,L$ and $\lambda_\tau=\max\{\tau,1-\tau\}$. If we choose $N=1$ and $L=\lfloor n^{(1-1/p)/(2\alpha_0+2)}\rfloor$, then $\mathcal{S}\leq \{d^2\max\{7d^*,20\}^2+d\max\{7d^*,20\}\}\times(12L+15)+K\times(2K)+2K\leq 2\times27^3\times (Kdd^*)^2\times\lfloor n^{(1-1/p)/(2\alpha_0+2)}\rfloor$ and
$$\mathbb{E}\big\{\mathcal{R}^\tau(\hat{f}_\phi)-\mathcal{R}^\tau(f_0)\big\}\leq C\mathcal{B}\times (Kdd^*)^2\log(Kdd^*) \times (\log n)^2
n^{-\left(1-\frac{1}{p}\right) \frac{\alpha_0}{\alpha_0+1}},
$$
where $C>0$ is a constant not depending on $n,d,d^*,K,\mathcal{B}$ and $\alpha_0$.

If Assumption \ref{quadratic} also holds, we have
\begin{equation*}
	\mathbb{E}\big\{\mathcal{R}^\tau(\hat{f}_\phi)-\mathcal{R}^\tau(f_0)\big\}\leq C\frac{\lambda_\tau\mathcal{B}\mathcal{S}\mathcal{D}\log(\mathcal{S})\log(n)}{n^{1-1/p}}+ 648c_\tau \lambda_0^2K(NL)^{-4\alpha_0},
\end{equation*}
where $c_\tau>0$ is a constant defined in Lemma \ref{lemma4}.
If we choose $N=1$ and $L=\lfloor n^{(1-1/p)/(4\alpha_0+2)}\rfloor$, then
$$\mathbb{E}\big\{\mathcal{R}^\tau(\hat{f}_\phi)-\mathcal{R}^\tau(f_0)\big\}\leq C_0 \mathcal{B}\times (Kdd^*)^2\log(Kdd^*) \times (\log n)^2
n^{-\left(1-\frac{1}{p}\right) \frac{2\alpha_0}{2\alpha_0+1}},
$$
where $C_0>0$ is a constant not depending on $n,d,d^*,K,\mathcal{B}$ and $\alpha_0$.

Furthermore,  if Assumption \ref{calib} also holds, it follows from Theorem \ref{thm3} that
{\color{black}
$$\mathbb{E} \big\{\Delta^2 (\hat{f}_\phi,f_0)\big\} \leq   c_{\kappa,\gamma}C_0 \mathcal{B}\times (Kdd^*)^2\log(Kdd^*) \times (\log n)^2
n^{-\left(1-\frac{1}{p}\right) \frac{2\alpha_0}{2\alpha_0+1}}.
$$
}

\subsection{Projection pursuit}
The projection pursuit model assumes
\begin{equation}
\label{ppM}
f_0(x)=\sum_{k=1}^K g_k(\theta_k^\top x),  \quad x\in\mathbb{R}^d,
\end{equation}
where $K\in\mathbb{N}$, $g_k:\mathbb{R}\to\mathbb{R}$ and $\theta_k\in\mathbb{R}^d$ \citep{friedman1981projection}. Such $f_0$ can be written as
$$f_0=h_2\circ h_1\circ h_0,$$
where $h_0(x)=\Theta x$ is a linear transformation from $\mathbb{R}^d$ to $\mathbb{R}^{K}$ with $\Theta=[\theta_1,\ldots,\theta_K]^\top $, $h_1(x)=(g_1(x),\ldots,g_K(x))^\top$ and $h_2(x)=\sum_{i=1}^Kx_i$ for $x=(x_1,\ldots,x_k)^\top\in\mathbb{R}^K$.  Correspondingly, $d_0=t_0=d, d_1=K,t_1=1, d_2=t_2=K$ and $d_3=1$. Suppose that Assumptions \ref{structure}-\ref{moment} and those conditions in Theorem \ref{thm2} are satisfied, where $g_i$ is H\"older continuous with order $\alpha_1$ and constant $\lambda_1$,  $i=1,\ldots,K$. By Theorem \ref{thm2}, given any $N,L\in\mathbb{N}^+$, for the function class of ReLU multi-layer perceptrons $\mathcal{F}_{\phi}=\mathcal{F}_{\mathcal{D},\mathcal{W},\mathcal{U},\mathcal{S},\mathcal{B}}$
with width $\mathcal{W}=\max\{2d,K(12N+8)\}$ and depth $\mathcal{D}=12L+19$,
for $2n \ge \text{Pdim}(\mathcal{F}_\phi)$, the excess risk of the \textit{DQR} estimator $\hat{f}_\phi$ satisfies
\begin{align*}
	\mathbb{E}\big\{\mathcal{R}^\tau(\hat{f}_\phi)-\mathcal{R}^\tau(f_0)\big\}\leq&  C\frac{\lambda_\tau\mathcal{B}\mathcal{S}\mathcal{D}\log(\mathcal{S})\log(n)}{n^{1-1/p}}+36\lambda_\tau\lambda_1\sqrt{K}(NL)^{-2\alpha_1},
\end{align*}
where  $C>0$ is a constant that does not depend on $n,d,\tau,\mathcal{B},\mathcal{S},\mathcal{D},\lambda_1,\alpha_1,N,L$ and $\lambda_\tau=\max\{\tau,1-\tau\}$.  If we choose $N=1$ and $L=\lfloor n^{(1-1/p)/(2\alpha_1+2)}\rfloor$, then $\mathcal{S}\leq \{(20K)^2+20K\}\times(12L+15)+d\times(2d)+2d+2d\times K +K\times 2K+2K\leq 20\times21\times27\times \max\{K,d\}^2\times \lfloor n^{(1-1/p)/(2\alpha_1+2)}\rfloor$ and
$$\mathbb{E}\big\{\mathcal{R}^\tau(\hat{f}_\phi)-\mathcal{R}^\tau(f_0)\big\}\leq C \mathcal{B}\times \max\{K,d\}^2\log(\max\{K,d\}) (\log n)^2
n^{-\left(1-\frac{1}{p}\right) \frac{\alpha_1}{\alpha_1+1}},
$$
where $C>0$ is a constant not depending on $n,d,\mathcal{B}$ and $\alpha_1$.

Additionally,  if Assumption \ref{quadratic}  holds, we have
\begin{align*}
	\mathbb{E}\big\{\mathcal{R}^\tau(\hat{f}_\phi)-\mathcal{R}^\tau(f_0)\big\}\leq&  C\frac{\lambda_\tau\mathcal{B}\mathcal{S}\mathcal{D}\log(\mathcal{S})\log(n)}{n^{1-1/p}}+648c_\tau\lambda_1^2{K}(NL)^{-4\alpha_1},
\end{align*}
where $c_\tau>0$ is a constant defined in Lemma \ref{lemma4}.
Alternatively, if we choose $N=1$ and $L=\lfloor n^{(1-1/p)/(4\alpha_1+2)}\rfloor$, then
$$\mathbb{E}\big\{\mathcal{R}^\tau(\hat{f}_\phi)-\mathcal{R}^\tau(f_0)\big\}\leq C_0 \mathcal{B}\times \max\{K,d\}^2\log(\max\{K,d\}) (\log n)^2
n^{-\left(1-\frac{1}{p}\right) \frac{2\alpha_1}{2\alpha_1+1}},
$$
and $C_0>0$ is a constant not depending on $n,d,\mathcal{B},K$ and $\alpha_1$.

Furthermore,  if Assumption \ref{calib} holds,  Theorem \ref{thm3} implies that
{\color{black}
$$\mathbb{E} \big\{\Delta^2 (\hat{f}_\phi,f_0)\big\} \leq c_{\kappa,\gamma}C_0 \mathcal{B}\times \max\{K,d\}^2\log(\max\{K,d\}) (\log n)^2
n^{-\left(1-\frac{1}{p}\right) \frac{2\alpha_1}{2\alpha_1+1}}.
$$
}

\subsection{The univariate composite model}

The univariate composite model \citep{horowitz2007}  takes the form
\begin{equation}
\label{ucM}
f_0(x) =m\Big\{\sum_{j_1=1}^{K_1} m_{j_1}\Big(\sum_{j_2=1}^{K_2}m_{j_1,j_2}
\Big[\cdots \sum_{j_{q-1}=1}^{K_{q-1}}m_{j_1, \ldots, j_{q-1}}\Big\{
\sum_{j_q=1}^{K_q}m_{j_1, \ldots, j_q}(x^{j_1, \ldots, j_q})
\Big\}\Big]\Big)\Big\},
\end{equation}
where $m$, $m_1, \ldots, m_{L_1, \ldots, K_q}$ are unknown univariate functions and
$x^{j_1, \ldots, j_q}$ are  one-dimensional elements of $x\in \mathbb{R}^d$,
which could be identical for two different indices $(j_1, \ldots, j_q)$. According to our notation, the target function $f_0$ can be written as
$$f_0=h_{2q}\circ\cdots\circ h_0,$$
where $h_{2q}(\cdot)=m(\cdot)$ and $h_{2i}(\cdot)=(m_{1,\cdots,1}(\cdot),\ldots,m_{j_1,\cdots,j_{q-i}}(\cdot),\cdots,m_{K_1,\cdots,K_{q-i}}(\cdot))^\top$ for $i=0,\ldots,q-1$ are all univariate functions. Correspondingly, $d_0=K_q,t_0=1,d_1=t_1=K_q,d_2=K_{q-1},t_2=1,\ldots,d_{q-2}=K_1,t_{q-2}=1,d_{2q-1}=t_{2q-1}=K_1,d_{2q}=t_{2q}=1$ and $d_{2q+1}=1$.
Suppose that Assumptions \ref{structure}-\ref{moment} and those conditions in Theorem \ref{thm2} hold, where $ m_{1,\cdots,1}(\cdot),\ldots,m_{j_1,\cdots,j_{q-i}}(\cdot),\cdots,m_{K_1,\cdots,K_{q-i}}(\cdot)$ are H\"older continuous with order $\alpha_i$ and constant $\lambda_i$ for $i=0,\ldots,q-1$, and $m$ is H\"older continuous with order $\alpha_q$ and constant $\lambda_q$. Then by Theorem \ref{thm2}, given any $N,L\in\mathbb{N}^+$, for the function class of ReLU multi-layer perceptrons $\mathcal{F}_{\phi}=\mathcal{F}_{\mathcal{D},\mathcal{W},\mathcal{U},\mathcal{S},\mathcal{B}}$
with width $\mathcal{W}=(12N+8) \Pi_{i=1}^qK_i$ and depth $\mathcal{D}=(12L+15)(q+1)+2q$,
for $2n \ge \text{Pdim}(\mathcal{F}_\phi)$, the excess risk of the \textit{DQR} estimator $\hat{f}_\phi$ satisfies
\begin{align*}
	\mathbb{E}\big\{\mathcal{R}^\tau(\hat{f}_\phi)-\mathcal{R}^\tau(f_0)\big\}\leq&  C\frac{\lambda_\tau\mathcal{B}\mathcal{S}\mathcal{D}\log(\mathcal{S})\log(n)}{n^{1-1/p}}+2\lambda_\tau\sum_{i=0}^{q} C_i^*\lambda_i^* K_i^*(NL)^{-2\alpha_i^*},
\end{align*}
where  $C>0$ is a constant not depending on $n,d,\tau,\mathcal{B},\mathcal{S},\mathcal{D},N,L,C_i^*,\lambda_i^*,\alpha_i^*$, $\lambda_\tau=\max\{\tau,1-\tau\}$ and $C_i^*=18^{\Pi_{j=i+1}^{q}\alpha_{j}}$, $\lambda_i^*=\Pi_{j=i}^{q}\lambda_{j}^{\Pi_{k=j+1}^{q} \alpha_{k}}$, $\alpha_i^*=\Pi_{j=i}^q \alpha_j$ and $K_i^*={(\Pi_{j=i}^{q}\sqrt{K_{q-j+1}}^{\Pi_{k=j}^{q} \alpha_{k}})}$.
To specify the network parameters, we set $N=1$, $L=\lfloor n^{(1-1/p)/(2\alpha_0^*+2)}\rfloor$ and let $K_0=1$.  Then $\mathcal{S}\leq (12L+15)\sum_{i=0}^{q}(20^2\Pi_{j=0}^iK_j^2+20\Pi_{j=0}^iK_j)+\sum_{i=0}^q(2K_i^2+2K_iK_{i+1})\leq 20\times21\times27\times (q+1)\Pi_{j=0}^qK_i^2\times \lfloor n^{(1-1/p)/(2\alpha_0^*+2)}\rfloor$ and
$$\mathbb{E}\big\{\mathcal{R}^\tau(\hat{f}_\phi)-\mathcal{R}^\tau(f_0)\big\}\leq C \mathcal{B}\times (\Pi_{j=0}^qK_i)^2\log(\Pi_{j=0}^qK_i) (\log n)^2
n^{-\left(1-\frac{1}{p}\right)\frac{\alpha_0}{\alpha_0+1}},
$$
where $C>0$ is a constant independent of $n,d,\mathcal{B},K_i$ and $\alpha^*_0$.

If Assumption \ref{quadratic} also holds, we have
\begin{align*}
	\mathbb{E}\big\{\mathcal{R}^\tau(\hat{f}_\phi)-\mathcal{R}^\tau(f_0)\big\}\leq&  C\frac{\lambda_\tau\mathcal{B}\mathcal{S}\mathcal{D}\log(\mathcal{S})\log(n)}{n^{1-1/p}}+2c_\tau\big[\sum_{i=0}^{q} C_i^*\lambda_i^* K_i^*(NL)^{-2\alpha_i^*}\big]^2,
\end{align*}
where $c_\tau>0$ is a constant defined in Lemma \ref{lemma4}.  If we choose $N=1$ and $L=\lfloor n^{(1-1/p)/(4\alpha_0^*+2)}\rfloor$, then
$$\mathbb{E}\big\{\mathcal{R}^\tau(\hat{f}_\phi)-\mathcal{R}^\tau(f_0)\big\}\leq C_0\mathcal{B}\times (\Pi_{j=0}^qK_i)^2\log(\Pi_{j=0}^qK_i) (\log n)^2
n^{-\left(1-\frac{1}{p}\right)\frac{\alpha_0}{\alpha_0+1}},
$$
where $C_0>0$ is a constant independent of $n,d,\mathcal{B},K_i$ and $\alpha^*_0$.

Moreover,  if Assumption \ref{calib} holds, it follows from Theorem \ref{thm3} that
{\color{black}
$$\mathbb{E} \big\{\Delta^2 (\hat{f}_\phi,f_0)\big\} \leq c_{\kappa,\gamma}C_0 \mathcal{B}\times (\Pi_{j=0}^qK_i)^2\log(\Pi_{j=0}^qK_i) (\log n)^2
n^{-\left(1-\frac{1}{p}\right)\frac{\alpha_0}{\alpha_0+1}}.
$$
}

\subsection{Generalized hierarchical interaction model}
\label{ghim}

Another general model is the {\it generalized hierarchical interaction model} of order $d^*$ and level $l$
\citep{bauer2019deep}.
For $d^*\in\{1,\ldots, d\}, l\in\mathbb{N}$ and $f_0:\mathbb{R}^d\to\mathbb{R}$, the generalized hierarchical interaction  model is defined as follows:
\begin{itemize}
	\item[(a)] The function $f_0$ satisfies a generalized hierarchical interaction model of order $d^*$ and level $0$, if there exist $\theta_1,\ldots, \theta_{d^*}\in\mathbb{R}^d$ and $f:\mathbb{R}^{d^*}\to\mathbb{R}$ such that
\begin{equation}
\label{miM}
f_0(x)=f(\theta_1^\top x,\ldots,\theta_{d^*}^\top x) \quad {\rm for\ all}\ x\in\mathbb{R}^d;
\end{equation}
	\item[(b)] The function $f_0$ satisfies a generalized hierarchical interaction model of order $d^*$ and level $l+1$, if there exist $K\in\mathbb{N}$, $g_k:\mathbb{R}^{d^*}\to\mathbb{R} \ (k=1,\ldots,K)$ and $f_{1,k},\ldots,f_{d^*,k}:\mathbb{R}^d\to\mathbb{R} \ (k=1,\ldots,K)$ such that $f_{1,k},\ldots,f_{d^*,k}   (k=1,\ldots,K)$ satisfy a generalized hierarchical interaction model of order $d^*$ and level $l$ and
\begin{equation}
\label{ghiMe}
	f_0(x)=\sum_{k=1}^K g_k(f_{1,k}(x),\ldots,f_{d^*,k}(x))\quad {\rm for \ all}\ x\in\mathbb{R}^d;
\end{equation}
	\item[(c)] the generalized hierarchical interaction model defined above is $\beta$-H\"older smooth if all the functions involve in its definition are $\beta$-H\"older smooth.
\end{itemize}
The generalized hierarchical interaction model includes the aforementioned models as special cases.  For instance, the single index model belongs to the class of generalized hierarchical interaction models of  order $1$ and level $0$; the additive model and projection pursuit correspond to order $1$ and level $1$; the interaction model is in conformity with order $d^*$ and level $1$; the univariate composite model in
 \citet{horowitz2007}  is a generalized hierarchical interaction model of order $1$ and level $q+1$. Moreover, the level zero generalized hierarchical interaction model (\ref{miM}) is the semiparametric multiple index model used in the sufficient dimension reduction
 \citep{li1991sir}.

In the  generalized hierarchical interaction models,
the target function $f_0$ is a composition of multi-index model and $d^*$-dimensional smooth functions, which resembles a multilayer feedforward neural networks in terms of the compositional structure. \citet{bauer2019deep} showed that the convergence rate of the least squares estimator based on sigmoid or bounded continuous activated deep regression networks is
$C_{d, d^*}(\log n)^3n^{-2\beta/(2\beta+d^*)}$. However, in their result, how the prefactor $C_{d,d^*}$ depends on $(d, d^*)$ is unclear.

 For the generalized hierarchical interaction model of order $d^*$ and level $l$ ($d^*\in\{1,\ldots,d\}$ and $l\in\mathbb{N}$) studied in \citet{bauer2019deep}, the target function $f_0$ is  a composition of multi-index model and $d^*$-dimensional smooth functions, which can be written as
$$f_0=h_{2l-1}\circ\cdots\circ h_0,$$  where $h_{2i}(\cdot)=(m_{1,\cdots,1}(\cdot),\ldots,m_{j_1,\cdots,j_{l-i}}(\cdot),\cdots,m_{K_1,\cdots,K_{l-i}}(\cdot))^\top$ for $i=0,\ldots,l-1$ are all $d^*$-variate functions and $h_{2i+1}(x)=\sum_{j=1}^{K_{l-i}} x_j$ for $x=(x_1,\ldots,x_{K_{l-i}})^\top\in\mathbb{R}^{K_{l-i}}$ and $i=0,\ldots,l-1$. Correspondingly, $d_0=K_l,t_0=d^*,d_1=t_1=K_l,d_2=K_{l-1},t_2=d^*,\ldots,d_{l-2}=K_1,t_{l-2}=d^*,d_{2l-1}=t_{2l-1}=K_1$ and $d_{2l}=t_{2l}=1$. Suppose that Assumptions \ref{structure}-\ref{moment} and those conditions in Theorem \ref{thm2} are satisfied, where $m_{1,\cdots,1}(\cdot),\ldots,m_{j_1,\cdots,j_{l-i}}(\cdot),\cdots,m_{K_1,\cdots,K_{l-i}}(\cdot)$ are H\"older continuous with order $\alpha_i$ and constant $\lambda_i$ for $i=0,\ldots,l-1$. Then by Theorem \ref{thm2}, given any $N,L\in\mathbb{N}^+$, for the function class of ReLU multi-layer perceptrons $\mathcal{F}_{\phi}=\mathcal{F}_{\mathcal{D},\mathcal{W},\mathcal{U},\mathcal{S},\mathcal{B}}$
 with width $\mathcal{W}=\max\{4d^*\lfloor N^{1/d^*}\rfloor+3d^*,12N+8\}\Pi_{i=1}^lK_i$ and depth $\mathcal{D}=(12L+17)l$,
 for $2n \ge \text{Pdim}(\mathcal{F}_\phi)$, the excess risk of the \textit{DQR} estimator $\hat{f}_\phi$ satisfies
 \begin{align*}
 	\mathbb{E}\big\{\mathcal{R}^\tau(\hat{f}_\phi)-\mathcal{R}^\tau(f_0)\big\}\leq&  C\frac{\lambda_\tau\mathcal{B}\mathcal{S}\mathcal{D}\log(\mathcal{S})\log(n)}{n^{1-1/p}}+2\lambda_\tau\sum_{i=0}^{q} C_i^*\lambda_i^* K_i^*(NL)^{-2\alpha_i^*/d^*},
 \end{align*}
 where  $C>0$ is a constant independent of $n,d,\tau,\mathcal{B},\mathcal{S},\mathcal{D},N,L,C_i^*,\lambda_i^*,\alpha_i^*$, $\lambda_\tau=\max\{\tau,1-\tau\},$  $C_i^*=18^{\Pi_{j=i+1}^{l}\alpha_{j}}$, $\lambda_i^*=\Pi_{j=i}^{l}\lambda_{j}^{\Pi_{k=j+1}^{l} \alpha_{k}}$, $\alpha_i^*=\Pi_{j=i}^l \alpha_j$ and $K_i^*={(\Pi_{j=i}^{l}\sqrt{K_{l-j+1}d^*}^{\Pi_{k=j}^{l} \alpha_{k}})}/{d^*}^{\alpha_i/2}$.
To specify the network parameters, we choose $N=1$ and $L=\lfloor n^{(1-1/p)d^*/(2\alpha_0^*+d^*)}\rfloor$. Then we have
 $\mathcal{S}\leq (12L+15)\sum_{i=0}^{l}(\max\{7d^*,20\}^2\Pi_{j=0}^iK_j^2+\max\{7d^*,20\}\Pi_{j=0}^iK_j)+\sum_{i=0}^l(2K_i^2+2K_iK_{i+1})\leq 7\times20\times21\times27\times d^*\times (l+1)\Pi_{i=0}^qK_i^2\times \lfloor n^{(1-1/p)d^*/(2\alpha_0^*+d^*)}\rfloor$ and
 $$\mathbb{E}\big\{\mathcal{R}^\tau(\hat{f}_\phi)-\mathcal{R}^\tau(f_0)\big\}\leq C\mathcal{B}\times (d^*)^2 (\Pi_{i=0}^lK_i)^2\log(\Pi_{i=0}^lK_i) (\log n)^2
 n^{-\left(1-\frac{1}{p}\right)\frac{\alpha^*_0}{\alpha^*_0+d^*}}
 $$
 where $C>0$ is a constant that does not depend on $n,d^*,\mathcal{B},K_i$ and $\alpha^*_0$.

 If Assumption \ref{quadratic} also holds, we have
 \begin{align*}
	\mathbb{E}\big\{\mathcal{R}^\tau(\hat{f}_\phi)-\mathcal{R}^\tau(f_0)\big\}\leq&  C\frac{\lambda_\tau\mathcal{B}\mathcal{S}\mathcal{D}\log(\mathcal{S})\log(n)}{n^{1-1/p}}+2c_\tau\big[\sum_{i=0}^{q} C_i^*\lambda_i^* K_i^*(NL)^{-2\alpha_i^*/d^*}\big]^2,
\end{align*}
 where $c_\tau>0$ is a constant defined in Lemma \ref{lemma4}. Alternatively, choosing $N=1$ and $L=\lfloor n^{(1-1/p)d^*/(4\alpha_0^*+2d^*)}\rfloor$,
 we have
$$\mathbb{E}\big\{\mathcal{R}^\tau(\hat{f}_\phi)-\mathcal{R}^\tau(f_0)\big\}\leq C_0\mathcal{B}\times (d^*)^2 (\Pi_{i=0}^lK_i)^2\log(\Pi_{i=0}^lK_i) (\log n)^2
n^{-\left(1-\frac{1}{p}\right)\frac{2\alpha^*_0}{2\alpha^*_0+d^*}},
$$
 and $C_0>0$ is a constant
not depending on $n,d^*,\mathcal{B},K_i$ and $\alpha^*_0$

Furthermore,  if Assumption \ref{calib} holds, it follows from Theorem \ref{thm3} that
{\color{black}
$$\mathbb{E} \big\{\Delta^2 (\hat{f}_\phi,f_0)\big\} \leq c_{\kappa,\gamma}C_0 \mathcal{B}\times (d^*)^2 (\Pi_{i=0}^lK_i)^2\log(\Pi_{i=0}^lK_i) (\log n)^2
n^{-\left(1-\frac{1}{p}\right)\frac{2\alpha^*_0}{2\alpha^*_0+d^*}}.
$$
}

In summary, these examples demonstrate that the \textit{DQR} estimator is able to mitigate the curse of dimensionality by taking advantage of  the compositional structure of these models. The prefactor only depends quadratically on $d$, instead of exponentially on $d$ as in the existing results for least squares conditional mean regression using deep neural networks. However, even with a quadratic dependence on the $d$, the error bounds can still be large for a large $d$. In particular, based on the risk bounds obtained above, a sample size of a polynomial order of $d$ is needed to achieve a small excess risk.

\section{Approximation of composite functions} \label{sec3}
In this section, we establish the error bound for approximating composite functions
defined in Assumption \ref{structure} using deep ReLU neural networks.
To bound the excess risk in Lemma \ref{lemma2}, we must first bound the
 approximation error due to the use of neural networks in constructing the estimator, as represented in  the second term on the right side of
(\ref{lem1}) or (\ref{oracle}). The stochastic error term can be analyzed using the empirical process theory by computing the cover number of the class of neural networks, as is given in (\ref{oracle}). So the remaining crucial task is to deal with the approximation error.

We will  express the error bounds in terms of the network parameters, the dimensionality of
the  components of $f_0$ and their continuity indices. To describe smoothness, we use the concept of the modulus of continuity.

\begin{definition}[Modulus of continuity]\label{modulus}
	For a function $f:D\to\mathbb{R}$, let $\omega_f(\cdot)$ denote its modulus of continuity, i.e.,
	\begin{equation}
		\label{mod}
		\omega_f(r) :=\sup\{\vert f(x)-f(y)\vert:x,y\in D,\Vert x-y\Vert_2\leq r\}, {\rm for\ any\ } r\geq0.
	\end{equation}
\end{definition}

For a uniformly continuous function $f$, $\lim_{r\to0}\omega_f(r)=\omega_f(0)=0$.  In addition, based on the modulus of continuity, different equicontinuous families of functions can be defined. For instance, the modulus $\omega_f(r)= \theta r$ describes the $\theta$-Lipschitz continuity; the modulus $\omega_f(r)=\lambda r^\alpha$ with $\lambda,\alpha>0$ describes the H\"older continuity.

In our problem, rather than imposing smoothness condition directly on the target function $f_0$, we make smoothness assumptions on the components of $f_0$. We assume that the functions $h_{ij}:[a_i,b_i]^{d_i}\to[a_{i+1},b_{i+1}]^{d_{i+1}}$ are H\"older continuous with order $\alpha_i$ and constant $\lambda_i$, i.e.,
\begin{align*}
	\vert h_{ij}(x)-h_{ij}(y)\vert\leq \lambda_i\Vert x-y\Vert^{\alpha_i}\quad, \forall x,y\in D_{ij},\ {\rm for\ }j=1,\ldots,d_{i+1}.
\end{align*}

For ease of reference, we first state an important result on the error bounds for approximating a general continuous function $f_0:[0,1]^d\to\mathbb{R}$ using ReLU neural networks \citep{shen2019deep}.
Our error bounds on approximating a composite function build on this result.

\begin{lemma}[Theorem 2.1 of \citet{shen2019deep}] \label{lemma5}
	Given $f\in\mathcal{C}([0,1]^d)$, for any $L\in\mathbb{N}^+$ and $N\in\mathbb{N}^+$, there exists a function $\phi$ implemented by a ReLU FNN with width $\max\{4d\lfloor N^{1/d}\rfloor+3d,12N+8\}$ and depth $12L+14$ such that $\Vert\phi\Vert_{L^\infty(\mathbb{R}^d)}\leq\vert f(\textbf{0})\vert+\omega_f(\sqrt{d})$ and,
	$$\vert f(x)-\phi(x)\vert\leq18\sqrt{d}\omega_f(N^{-2/d}L^{-2/d}), \quad{\rm for\ any\ } x\in[0,1]^d\backslash\Omega([0,1]^d,K,\delta),$$
	where $K=\lfloor N^{1/d}\rfloor^2\lfloor L^{1/d}\rfloor^2$ and $\delta$ is an arbitrary number in $(0,1/(3K)]$, and the trifling region $\Omega([0,1]^d,K,\delta)$ of $[0,1]^d$ is defined as $$\Omega([0,1]^d,K,\delta)=\cup_{i=1}^d\{x=[x_1,x_2,...,x_d]^T:x_i\in\cup_{k=1}^{K-1}(k/K-\delta,k/K)\}.$$
	Especially, if $f$ is H\"older continuous of order $\alpha>0$ with constant $\lambda$, then
		$$\vert f(x)-\phi(x)\vert\leq18\sqrt{d}\lambda N^{-2\alpha/d}L^{-2\alpha/d}, \quad{\rm for\ any\ } x\in[0,1]^d\backslash\Omega([0,1]^d,K,\delta).$$
\end{lemma}

According to Lemma \ref{lemma5}, for a function $h_i:[a_i,b_i]^{d_i}\to[a_{i+1},b_{i+1}]^{d_{i+1}}$,  each of its components $h_{ij}:[a_i,b_i]^{t_i}\to\mathbb{R}$ can be approximated by a ReLU network. Then $d_i$ such (parallel) networks can be stacked to form a new ReLU network
for approximating  $h_i$.

\begin{lemma}[Parallel networks]\label{lemma6}
	Let $h=(h_{j})_j^\top: [0,1]^{d}\to\mathbb{R}^{m}$ be a continuous function, and suppose that $(h_{j})_j^\top, j=1,\ldots,m,$ are $t$-variate functions with the same modulus of continuity $\omega(\cdot)$. Then, for any $L\in\mathbb{N}^+$ and $N\in\mathbb{N}^+$, there exists a function $\phi$ implemented by a ReLU FNN with width $d\max\{4t\lfloor N^{1/t}\rfloor+3t,12N+8\}$ and depth $12L+14$ such that $\Vert\phi\Vert_{L^\infty(\mathbb{R}^d)}\leq \max_{j=1,\ldots,m} \vert h_j(\textbf{0})\vert+\omega(\sqrt{t})$ and
	$$\vert h(x)-\phi(x)\vert\leq18\sqrt{t}\omega(N^{-2/t}L^{-2/t}), \quad{\rm for\ any\ } x\in[0,1]^d\backslash\Omega([0,1]^d,K,\delta),$$
	where $K=\lfloor N^{1/d}\rfloor^2\lfloor L^{1/d}\rfloor^2$ and $\delta$ is an arbitrary number in $(0,1/(3K)]$.
\end{lemma}

By Lemma \ref{lemma6}, for a composite function $h_q\circ\cdots\circ h_0$, each function  $h_i$ in the composition can be approximated by a
 ReLU network $\tilde{h}_i$ under the H\"older continuity assumption.   It is thus natural to consider stacking these networks $\tilde{h}_i$ in a sequence as $\tilde{h}_q\circ\ldots\tilde{h}_0$ to approximate $h_q\circ\ldots\circ h_0$.

\begin{definition}[Norms of a vector of functions]
	For a function $h=(h_{j})_j^\top: \mathbb{R}^{d_{in}}\to\mathbb{R}^{d_{out}}$ with domain $D=D_1\otimes\ldots\otimes D_{d_{out}}$, we define its supremum-norm by the sup-norm of the vectors of its outputs,
	$$\Vert h\Vert_{L_\infty(D)}:=\sup_{x\in D}\Vert h(x)\Vert_\infty,$$
	and define its $L_2$-norm by the $L_2$ of the vectors of its outputs,
	$$\Vert h\Vert_{L_2(D)}:=\sup_{x\in D}\Vert h(x)\Vert_2.$$
\end{definition}

\begin{lemma} [Approximation by composition]\label{lemma7}
	Let $h_{ij}:\mathbb{R}^{t_{i}}\to\mathbb{R}$, $i=0,\ldots,q$ and $j=1,\ldots,d_{i+1}$ be H\"older continuous functions with order $\alpha_i\in[0,1]$ and constant $\lambda_i\ge0$ and let $h_i=(h_{ij})_j^\top: \mathbb{R}^{d_{i}}\to\mathbb{R}^{d_{i+1}}$ be vectors of functions with domain $D_i$. Then any functions $\tilde{h}_i=(\tilde{h}_{ij})_j^\top: \mathbb{R}^{d_{i}}\to\mathbb{R}^{d_{i+1}}$ with $\tilde{h}_{ij}:\mathbb{R}^{t_{i}}\to\mathbb{R}$, which have the same domain as $h_i,$ will satisfy,
	\begin{align*}
\Vert h_q\circ\ldots h_0-\tilde{h}_q\circ\ldots \tilde{h}_0\Vert_{L_\infty(D_0)}
		\leq\sum_{i=0}^q \Pi_{j=i+1}^{q}
\lambda_{j}^{\Pi_{k=j+1}^{q} \alpha_{k}} \Pi_{j=i+1}^{q}\sqrt{t_{j}}^{\Pi_{k=j}^{q} \alpha_{k}} \Vert h_i-\tilde{h}_i\Vert_{L_\infty(D_i)}^{\Pi_{j=i+1}^{q}\alpha_{j}}.
	\end{align*}
\end{lemma}

\begin{remark}
	Lemma \ref{lemma7} can be generalized without further difficulty for any other continuous functions $h_i$ with different types of modulus of continuity. The generalized result is expressed in term of the modulus of continuities of $h_i$, where the expression is analytical but complicated with a nested or compositional form of modulus functions.
\end{remark}
Note that the domains of $h_i$ are generally not $[0,1]^{d_i}$ as required in Lemma \ref{lemma5} and Lemma \ref{lemma6}. Thus the domain of the constructed ReLU networks have to be aligned with the approximated functions $h_{i}$. In light of this, we add an additional invertible linear layer ${A}_i(\cdot):D_{i}\to[0,1]^{d_i}$ at the beginning of each of the subnetworks $\tilde{h}_{i}$ in Lemma \ref{lemma6} for $i=1,\ldots,q$. With a slight abuse of notation, in the following we let $\tilde{h}_{i}$ denote the networks with an additional invertible linear layer as their first layer. In this case, $\tilde{h}_{i}:D_i\to\mathbb{R}^{d_{i+1}}$.

Moreover, there are many popular statistical models containing a linear function as a layer in a composite function, i.e.,  there exists some $i\in\{0,\ldots,q\}$ such that $h_i(x)=T_ix+u_i$ for some matrix $T_i\in\mathbb{R}^{d_{i}\times d_{i+1}}$ and $u_i\in\mathbb{R}^{d_{i+1}}$. For such a linear function $h_i$, it is possible to construct ReLU neural networks to approximate it perfectly.

\begin{lemma}[Approximation of linear functions]	\label{lemma8}
	Let $h=(h_{j})_j^\top: \mathbb{R}^{d}\to\mathbb{R}^{m}$ be a linear function, i.e. $h(x)=Tx+u$ with $T\in\mathbb{R}^{m\times d}$ and $u\in\mathbb{R}^{m}$.  Then there exists
a three-layer ReLU neural network $\tilde{h}$ with width vector $(d,2d,m)$ such that $\tilde{h}(x)=h(x)$ for any $x\in\mathbb{R}^d$.
\end{lemma}

By Lemma \ref{lemma8}, the approximation of composite functions can be further improved if some of the compositions are linear functions.

\begin{theorem} [Approximation of composite functions]\label{thm1}
	Let $H_q=h_q\circ\ldots\circ h_0$ be a function from $[a,b]^d$ to $\mathbb{R}$ and $h_i=(h_{ij})_j^\top: D_i\to\mathbb{R}^{d_{i+1}}, i=0,\ldots,q$ be vectors of functions with domain $D_i\subset \mathbb{R}^{d_{i}}$ where $h_{ij}:D_{ij}\to\mathbb{R}$, $i=0,\ldots,q$ and $j=1,\ldots,d_{i+1}$ with domain $D_{ij}\subset \mathbb{R}^{t_{i}}$ are H\"older continuous functions with order $\alpha_i\in[0,1]$ and constant $\lambda_i\ge0$.Then for any $L_i\in\mathbb{N}^+$ and $N_i\in\mathbb{N}^+$, there exist functions $\tilde{h}_i$ for $i=0,\ldots,q$ implemented by ReLU FNNs with width $d_{i}\max\{4t_i\lfloor N_i^{1/t_i}\rfloor+3t_i,12N_i+8\}$ and depth $12L_i+15$ such that $\Vert\tilde{h}_i\Vert_{L_i^\infty(\mathbb{R}^{d_i})}\leq \max_{j=1,\ldots,d_i} \vert h_{ij}(\textbf{0})\vert+\omega(\sqrt{t_i})$ and
	$$\vert \tilde{h}_i(x)-h_i(x)\vert\leq18\sqrt{t_i}\lambda_i(N_iL_i)^{-2\alpha_i/t_i}, \quad{\rm for\ any\ } x\in D_i\backslash\ A^{-1}_i(\Omega([0,1]^{d_i},K,\delta)),$$
	where $A_i:D_i\to[0,1]^{d_i}$ is an invertible linear layer (the first layer of $\tilde{h}_i$), $K_i=\lfloor N_i^{1/d_i}\rfloor^2\lfloor L_i^{1/d_i}\rfloor^2$ and $\delta_i$ is an arbitrary number in $(0,1/(3K_i)]$.
	
	Furthermore,  if $h_j$ are linear functions for $j\in J\subset \{0,\ldots,q\}$ with H\"older constant $\lambda_j=1$  and order $\alpha_{j}=1$, then there exists functions $\tilde{h}_j$ implemented by ReLU FNNs with width vector $(d_j,2d_j,d_{j+1})$ and depth 3 such that,
	$$\vert\tilde{h}_j(x)-h_j(x)\vert=0,\quad {\rm for\ any\ } x\in\mathbb{R}^{d_j}.$$
	
	Let $\tilde{H}_q=\tilde{h}_q\circ\ldots\circ \tilde{h}_0$ denote the function implemented by ReLU FNN with width no more than $\max_{i=0,\ldots,q}d_{i}\max\{4t_i\lfloor N_i^{1/t_i}\rfloor+3t_i,12N_i+8\}$ and depth $\sum_{i\in J^c}(12L_i+15)+2\vert J\vert$, where $\vert J\vert$ denotes its cardinality and $J^c:=\{0,\ldots,q\}\backslash J$, then we have
	\begin{align*}
		\vert \tilde{H}_q(x)-H_q(x)\vert\leq\sum_{i\in J^c} C_i^*\lambda_i^* t_i^*(N_iL_i)^{-2\alpha_i^*/t_i}, \qquad {\rm for\ any\ }x\in [a,b]^d\backslash \Omega_0,
	\end{align*}
	where $C_i^*=18^{\Pi_{j=i+1}^{q}\alpha_{j}}$, $\lambda_i^*=\Pi_{j=i}^{q}\lambda_{j}^{\Pi_{k=j+1}^{q} \alpha_{k}}$, $\alpha_i^*=\Pi_{j=i}^q \alpha_j$, $t_i^*={(\Pi_{j=i}^{q}\sqrt{t_{j}}^{\Pi_{k=j}^{q} \alpha_{k}})}/{\sqrt{t_i}^{\alpha_i}}$ and $\Omega_0$ is a subset of $[a,b]^d$ which satisfies
	$$ \Omega([0,1]^{d_i},K_i,\delta_i)\subseteq A_i\circ\tilde{h}_{i-1}\circ\cdots\circ \tilde{h}_0(\Omega_0), \qquad {\rm for\ } i=0,\ldots,q,$$
	where $A_j$ is defined as identity map for $j\in J$.
\end{theorem}
\vspace{0.1in}
\begin{remark}
	In Theorem \ref{thm1}, since $\tilde{h}_i, A_i$ are continuous mappings, the Lebesgue measure of $\Omega_0$ can be arbitrarily small as $\delta_i\in(0,1/(3K_i)]$ can  be arbitrarily small,
thus the Lebesgue measure of $\Omega([0,1]^{d_i},K_i,\delta_i)$ can be arbitrarily small.
\end{remark}

When all the component functions $h_{ij}$ are Lipschitz continuous,
the approximation error bound in
Theorem \ref{thm1} can be simplified considerably. Because Lipschitz continuity is a reasonable assumption in practice, we state the following corollary on the approximation error for Lipschitz continuous functions.

\begin{corollary}
	Suppose all $h_{ij}:D_{ij}\to\mathbb{R}$ in Theorem \ref{thm1} are Lipschitz continuous functions ($\alpha_i=1$ for $i=0,\ldots,q$) with Lipschitz constant $\lambda_i\ge0$. We set the same shape for each subnetwork with $N_0=\ldots=N_q=N\in\mathbb{N}^+$ and $L_0=\ldots=L_q=L\in\mathbb{N}^+$,
then we have
\begin{align*}
		\vert \tilde{H}_q(x)-H_q(x)\vert& \leq18\sum_{i=0}^q \big(\Pi_{j=i}^{q}\lambda_{j}\big)\big(\Pi_{j=i+1}^{q}\sqrt{t_{j}}\big) (NL)^{-2/t_i}\\
		&=18 \sum_{i=0}^q \lambda_i^* t_i^* (NL)^{-2/t_i}, \ \text{for any }  \ x\in[a,b]^d\backslash\Omega_0,
\end{align*}
where $\lambda_i^*=\Pi_{j=i}^{q}\lambda_{j}$ and $t_i^*=\Pi_{j=i+1}^{q}\sqrt{t_{j}}$.

Furthermore , if $h_j$ are linear functions for $j\in J\subset \{0,\ldots,q\}$, then we have
	\begin{align*}
	\vert \tilde{H}_q(x)-H_q(x)\vert &\leq18 \sum_{i\in J^c} \lambda_i^* t_i^* (NL)^{-2/t_i}, \
\text{for  any } \  x \in [a,b]^d\backslash\Omega_0.
\end{align*}
\end{corollary}
This lemma shows that, if $t_i\ll d_i$, the approximation rate improves, which lessens
 the curse of dimensionality.

\section{Numerical studies}\label{sim}
In this section, we compare deep quantile regression with traditional linear quantile regression and reproducing kernel methods on simulated data.  To be specific, we compare the following  methods of quantile regressions:
\begin{itemize}
	\item The traditional linear quantile regression as described in \cite{koenker1978}, denoted by \textit{linear QR}. Without regularization, the empirical risk is minimized over the parameter space (intercept included) $\mathbb{R}^{d+1}$ to give an linear estimator. These estimation are implemented on Python via package \textit{statsmodels}.
	\item   Kernel-based nonparametric quantile regression  as described in \cite{sangnier2016joint}, denoted by \textit{kernel QR}.
	This is a joint quantile regression method based on vector-valued reproducing kernel Hilbert space (RKHS), which enjoys fewer quantile crossings and enhanced performances compared to independent estimations and hard non-crossing constraints. In our implementation, the radial basis function (RBF) kernel is chosen and a coordinate descent primal-dual algorithm \citep{fercoq2019coordinate} is used via Python package \textit{qreg}.
	\item Deep quantile regression as described in Section \ref{sec2}, denoted by \textit{DQR}. We implement it in Python via \textit{Pytorch} and use \textit{Adam} \citep{kingma2014adam} as the optimization algorithm with default learning rate 0.01 and default $\beta=(0.9,0.99)$ (coefficients used for computing running averages of gradients and their squares).
	\item Deep least squares regression, denoted by \textit{DLS}. We minimize the mean square error on the training data to get the nonparametric least square estimator using deep neural networks. Similarly we implement it on Python via \textit{Pytorch} and use \textit{Adam} as the optimization algorithm with default settings. The comparison with DLS mainly focuses on the $0.5$-th quantile curve since the conditional mean and the conditional median coincident with each other when error is symmetric.
\end{itemize}

\subsection{Estimations and Evaluations}
We consider estimating the quantile curves at 5 different levels for each simulated model, i.e., we estimate quantile curves for $\tau\in\{0.05,0.25,0.5,0.75,0.95\}$. For each model $f_0$ and each error $\eta$, according to model (\ref{model}) we generate the training data $(X_i^{train},Y_i^{train})_{i=1}^n$ with sample size $n$ to train the empirical risk minimizer at $\tau\in\{0.05,0.25,0.5,0.75,0.95\}$ by different methods, i.e.
\begin{align*}
	\hat{f}^\tau_n\in\arg\min_{f\in\mathcal{F}} \frac{1}{n}\sum_{i=1}^n\rho_\tau(Y_i^{train}-f(X_i^{train})),
\end{align*}
where $\mathcal{F}$ is the class of linear functions, RKHS
or the class of
ReLU neural network functions. For each $f_0$ and each error $\eta$, we also generate the testing data $(X_t^{test},Y_t^{test})_{t=1}^T$ with sample size $T$ from the same distribution of the training data. Then for each obtained $\hat{f}^\tau_n$, we calculate its testing risk on $(X_t^{test},Y_t^{test})_{t=1}^T$, i.e.,
 \begin{align*}
 	\mathcal{R}^\tau(\hat{f}^\tau_n)=\frac{1}{T}\sum_{t=1}^T \rho_\tau(Y_t^{test}-\hat{f}^\tau_n(X_t^{test})).
 \end{align*}
Moreover, for each obtained $\hat{f}^\tau_n$, we calculate the $L_1$ distance between $\hat{f}^\tau_n$ and the corresponding risk minimizer $f_0^\tau$, i.e.
 \begin{align*}
\Vert\hat{f}^\tau_n-f_0^\tau\Vert_{L^1(\nu)}=\frac{1}{T}\sum_{t=1}^T \vert\hat{f}^\tau_n(X_t^{test})-f_0^\tau(X_t^{test})\vert,
\end{align*}
and we also calculate the $L_2$ distance between $\hat{f}^\tau_n$ and the corresponding risk minimizer $f_0^\tau$, i.e.
 \begin{align*}
	\Vert\hat{f}^\tau_n-f_0^\tau\Vert^2_{L^2(\nu)}=\frac{1}{T}\sum_{t=1}^T \vert\hat{f}^\tau_n(X_t^{test})-f_0^\tau(X_t^{test})\vert^2.
\end{align*}
All the $L_2$ test error results are provided in the appendix.
The specific forms of $f_0^\tau$ are given in the part on the data generation models below.

In the simulation studies, we take $T=100,000$ as the sample size of testing data for each data generation model. We report the mean and standard deviation of statistics including excess risk $\mathcal{R}^\tau(\hat{f}^\tau_n)-\mathcal{R}^\tau({f}^\tau_0)$, $L_1$ distance and $L^2_2$ distance  over $R = 10$ replications under different scenarios.  For \textit{DLS}, the testing risk and the excess risk are calculated in terms of mean squares loss function other than the check loss $\rho_\tau$.

\subsection{Data generation: univariate models}
We generate data according to model (\ref{model}), i.e., $Y=f_0(X)+\eta$. We consider three basic univariate models, including ``Linear'', ``Wave'' and ``Triangle'', which corresponds to different specifications of $f_0$. The formulae are given below.
\begin{enumerate}[(a)]
\setlength\itemsep{-0.05 cm}
	\item Linear:
	\begin{align*} f_0(x)=2x.
\end{align*}
	\item Wave:
	$$f_0(x)=2x\sin(4\pi x).$$
	\item Triangle:
	$$f_0(x)=4(1-\vert x-0.5\vert).$$
\end{enumerate}
We use the linear model as a baseline model in our simulations and expect all the methods
perform well under the linear model. The ``Wave'' is a nonlinear smooth model and the ``Triangle'' is a
nonlinear continuous but non-differentiable model. These models
are chosen so that we can evaluate the
performance of \textit{DQR},  \textit{kernel QR} and \textit{linear QR} under different types of  models.

\begin{figure}[H]
\centering
\includegraphics[width=\textwidth, height=4.5 in]{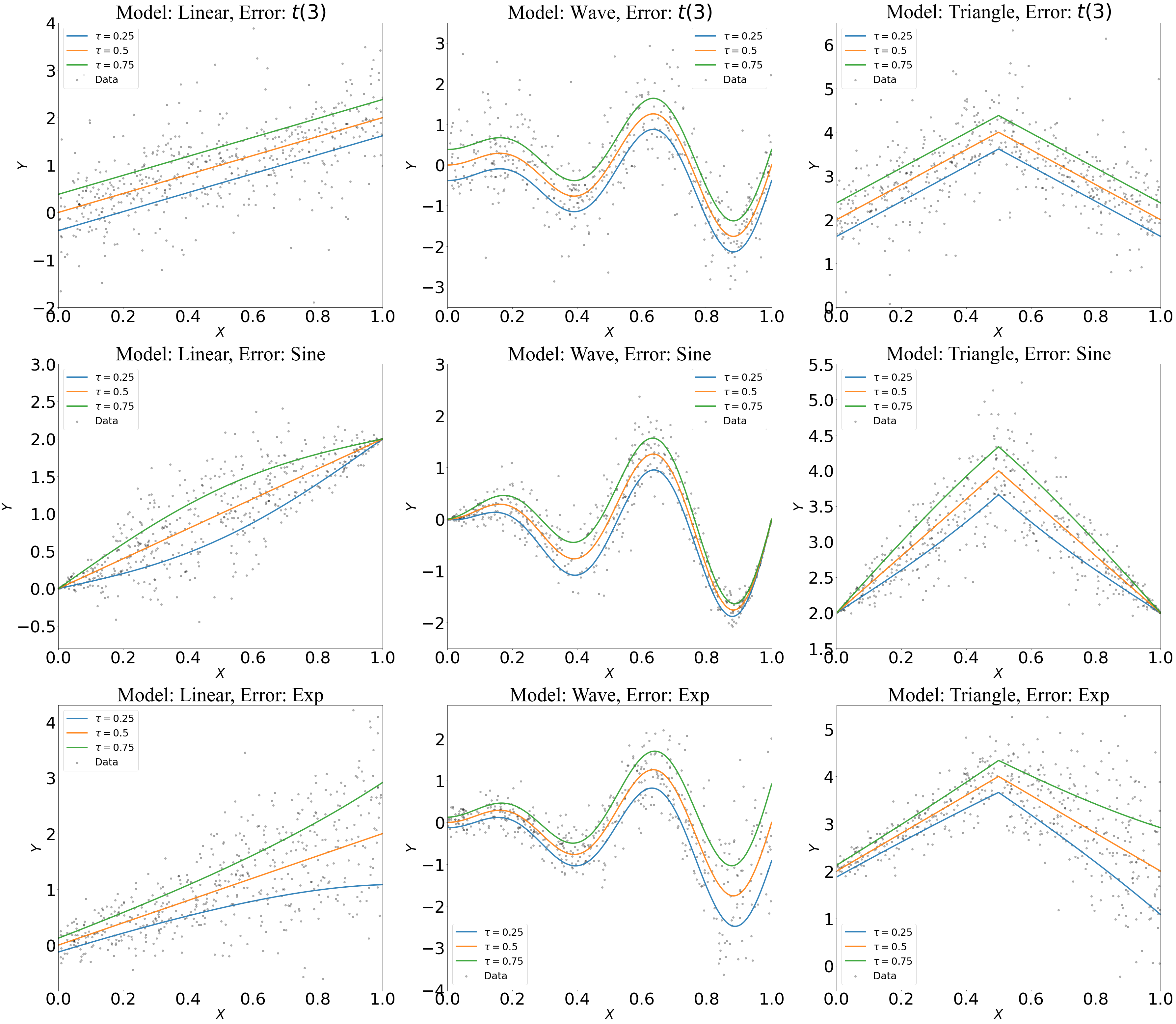}
\caption{The target quantiles curves at $\tau=0.25, 0.50. 0.75$
under different models and error distributions. From the left to the right, each column corresponds a data generation model, ``Linear'', ``Wave'' and ``Triangle''. From the top to the bottom, each row corresponds a error distribution,  $t(3)$, ``\textit{Sine}'' and ``\textit{Exp}''.}
\label{fig:target}
\end{figure}

For these models, we generate $X$ uniformly from the unit interval $[0,1]$. We generate the error $\eta$ from the following distributions.
\begin{enumerate}[(i)]
\setlength\itemsep{-0.05 cm}
	\item $\eta$ follows a scaled Student's t distribution with degrees of freedom 3, i.e., $\eta\sim 0.5\times t(3)$, denoted by $t(3)$;
	\item  Conditioning on $X=x$, the error $\eta$ follows a normal distribution of which variance depends on the covariate $X$, i.e., $\eta\mid X=x\sim 0.5\times\mathcal{N}(0,[\sin(\pi x)]^2)$, denoted by \textit{Sine};
	\item Conditioning on $X=x$, the error $\eta$ follows a normal distribution of which variance depends on the covariate $X$, i.e., $\eta\mid X=x\sim0.5\times\mathcal{N}(0,\exp(4x-2))$, denoted by \textit{Exp}.
\end{enumerate}
Note that except for $t(3)$, other two types of errors depend on the predictor $X$.
The $\tau$-th conditional quantile $f_0^\tau(x)$ of the response $Y$ given $X=x$ can be calculated by
$$f_0^\tau(x)=f_0(x)+F^{-1}_{\eta\mid X=x}(\tau),$$
where $F^{-1}_{\eta\mid X=x}(\cdot)$ is the inverse of the conditional cumulated distribution function of $\eta$ given $X=x$. For $t(3)$ error, $\eta$ is independent with $X$, then $F^{-1}_{\eta\mid X=x}(\cdot)$ is simply the inverse of distributional function of the $2 t(3)$. For the \textit{Sine} error, $F^{-1}_{\eta\mid X=x}(\tau)=0.5\times\sin(\pi x)\times\Phi^{-1}(\tau)$ where $\Phi^{-1}(\cdot)$ is the inverse of the CDF of a standard normal random variable. Similarly, for the \textit{Exp} error, $F^{-1}_{\eta\mid X=x}(\tau)=0.5\times\exp(2x-1)\times\Phi^{-1}(\tau)$.
Figure \ref{fig:target} shows all these univariate data generation models and their corresponding conditional quantiles at $\tau=0.25, 0.50, 0.75$.

\begin{figure}[H]
	\centering
	\begin{subfigure}{\textwidth}
		\includegraphics[width=1\textwidth]{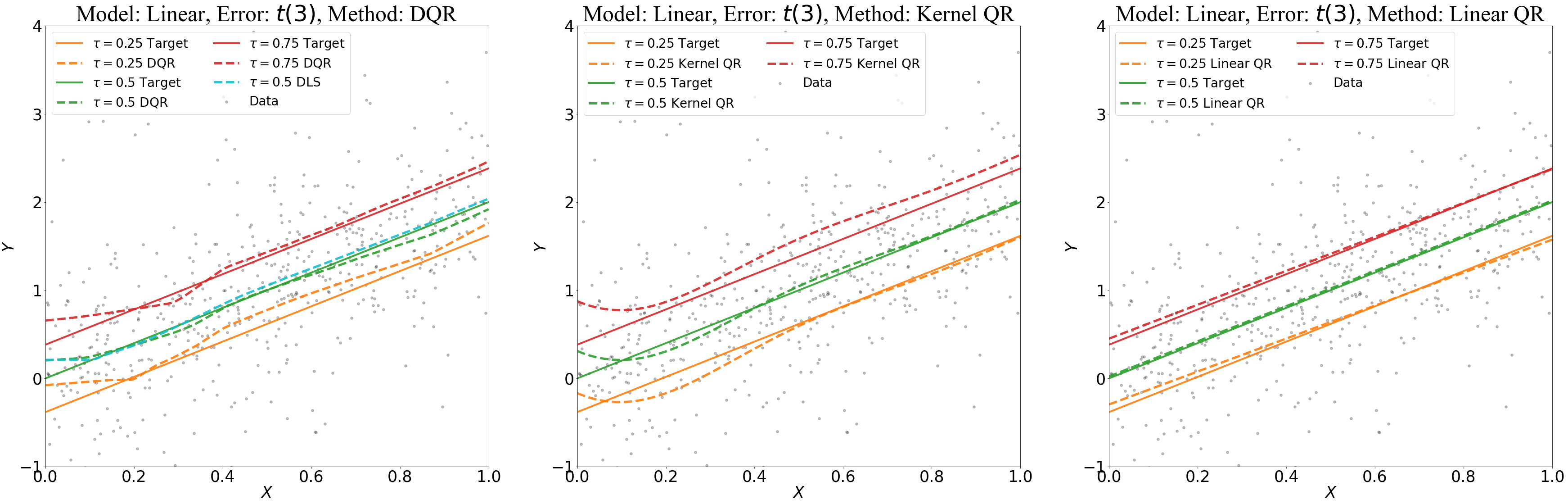}
	\end{subfigure}
	\begin{subfigure}{\textwidth}
		\includegraphics[width=1\textwidth]{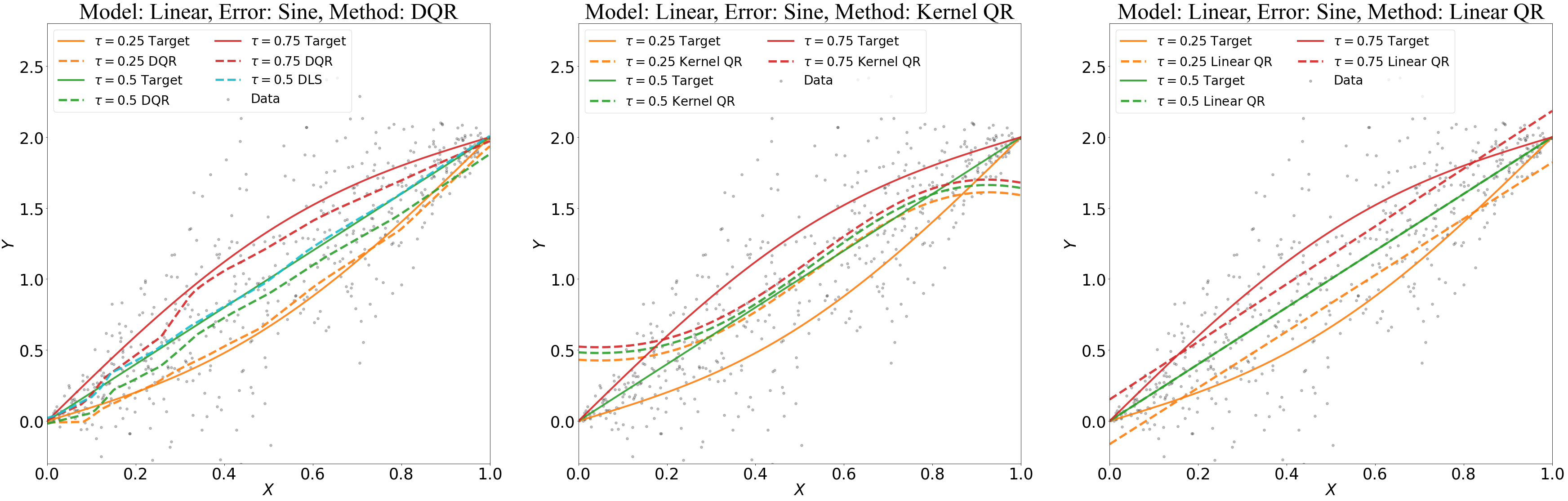}
	\end{subfigure}
	\begin{subfigure}{\textwidth}
		\includegraphics[width=\textwidth]{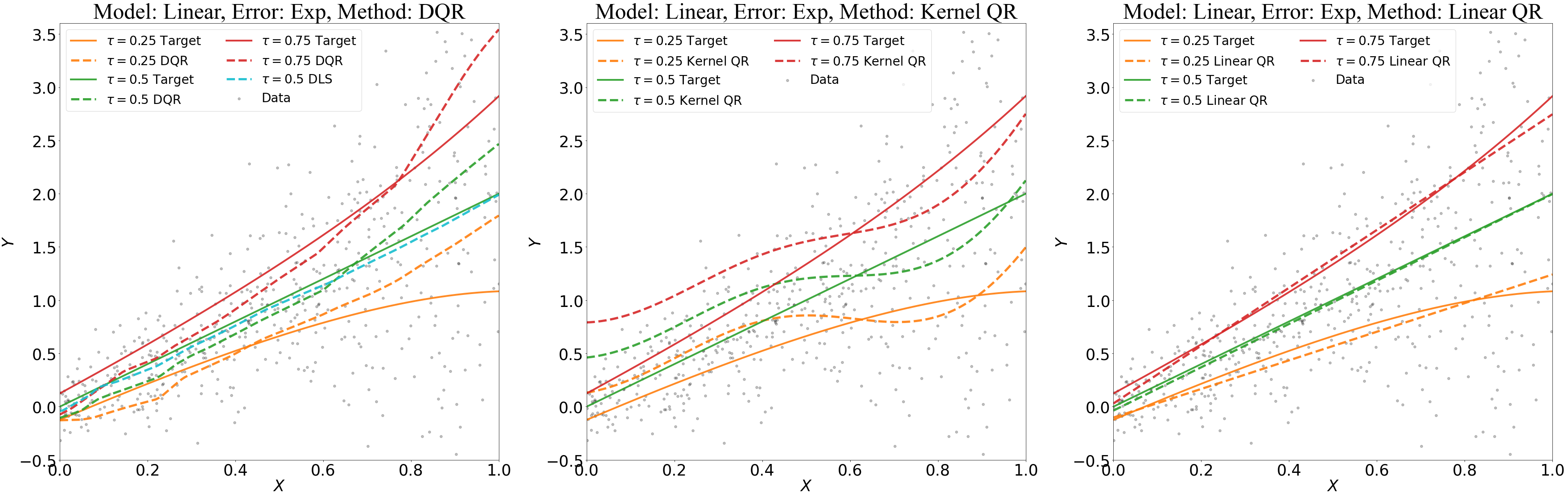}
	\end{subfigure}
\caption{The fitted quantile curves by different methods under the univariate model ``Linear'' with different errors. The training data is depicted as grey dots.The target quantile functions at the quantile levels $\tau=$0.25 (yellow), 0.5 (green), 0.75 (red)  are depicted as solid curves, and the estimated quantile functions are represented by dashed curves with the same color.
From the top to the bottom, the rows correspond to the errors $t(3)$, ``\textit{Sine}'' and ``\textit{Exp}''.
From the left to the right, the subfigures correspond to the methods \textit{DQR}, \textit{kernel QR} and \textit{linear QR}. The fitted \textit{DLS} curve (in blue)  is included in the left subfigure.}
	\label{fitted:linear1}
\end{figure}

\begin{figure}[H]
	\centering
	\begin{subfigure}{\textwidth}
		\includegraphics[width=1\textwidth]{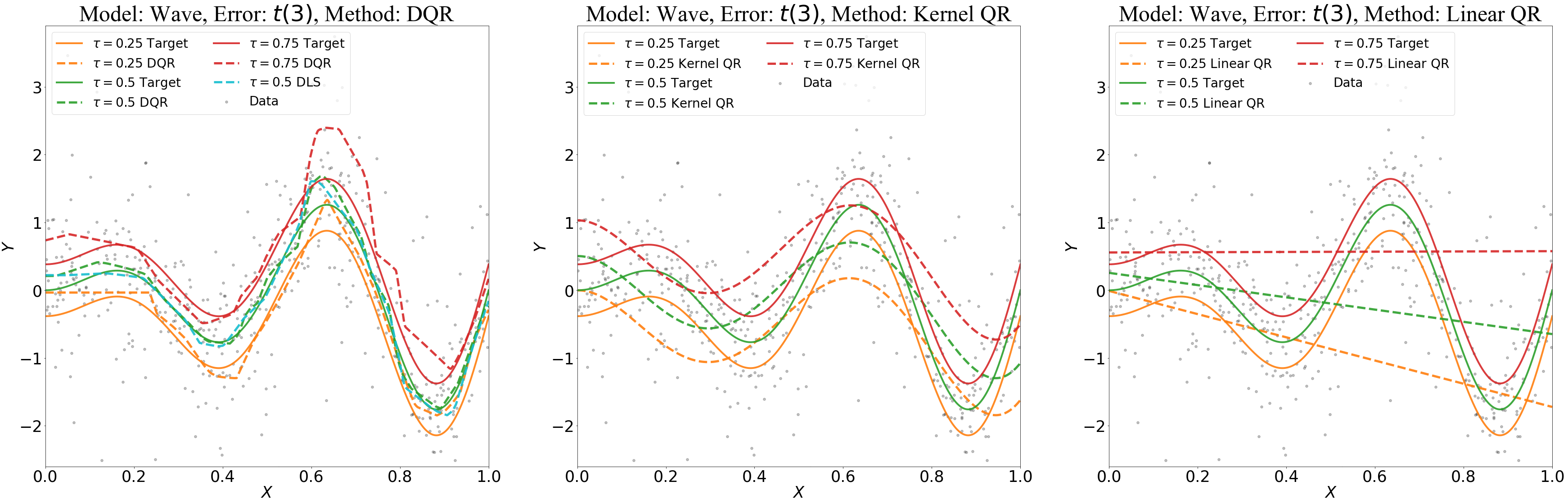}
	\end{subfigure}
	\begin{subfigure}{\textwidth}
		\includegraphics[width=1\textwidth]{fitted_wave_sine_1.png}
	\end{subfigure}
	\begin{subfigure}{\textwidth}
		\includegraphics[width=\textwidth]{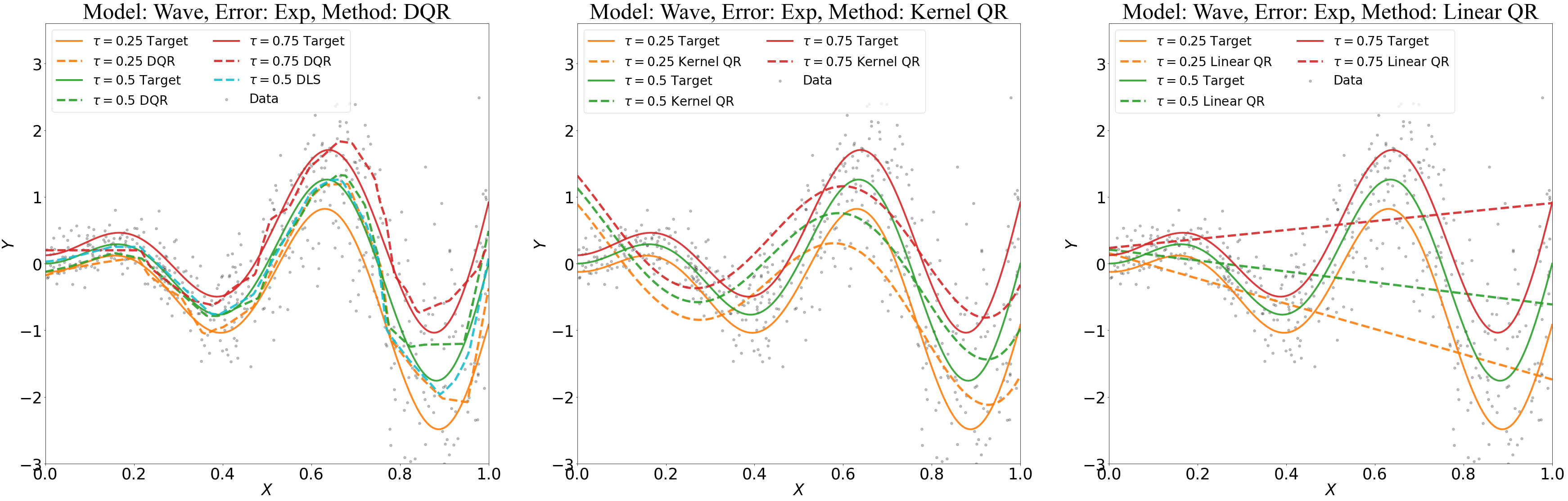}
	\end{subfigure}
\caption{The fitted quantile curves by different methods under the univariate model ``Wave'' with different errors. The training data is depicted as grey dots.The target quantile functions at the quantile levels $\tau=$0.25 (yellow), 0.5 (green), 0.75 (red)  are depicted as solid curves, and the estimated quantile functions are represented by dashed curves with the same color.
From the top to the bottom, the rows correspond to the errors $t(3)$, ``\textit{Sine}'' and ``\textit{Exp}''.
From the left to the right, the subfigures correspond to the methods \textit{DQR}, \textit{kernel QR} and \textit{linear QR}. The fitted \textit{DLS} curve (in blue)  is included in the left subfigure.}
	\label{fitted:wave1}
\end{figure}

\begin{figure}[H]
	\centering
	\begin{subfigure}{\textwidth}
		\includegraphics[width=1\textwidth]{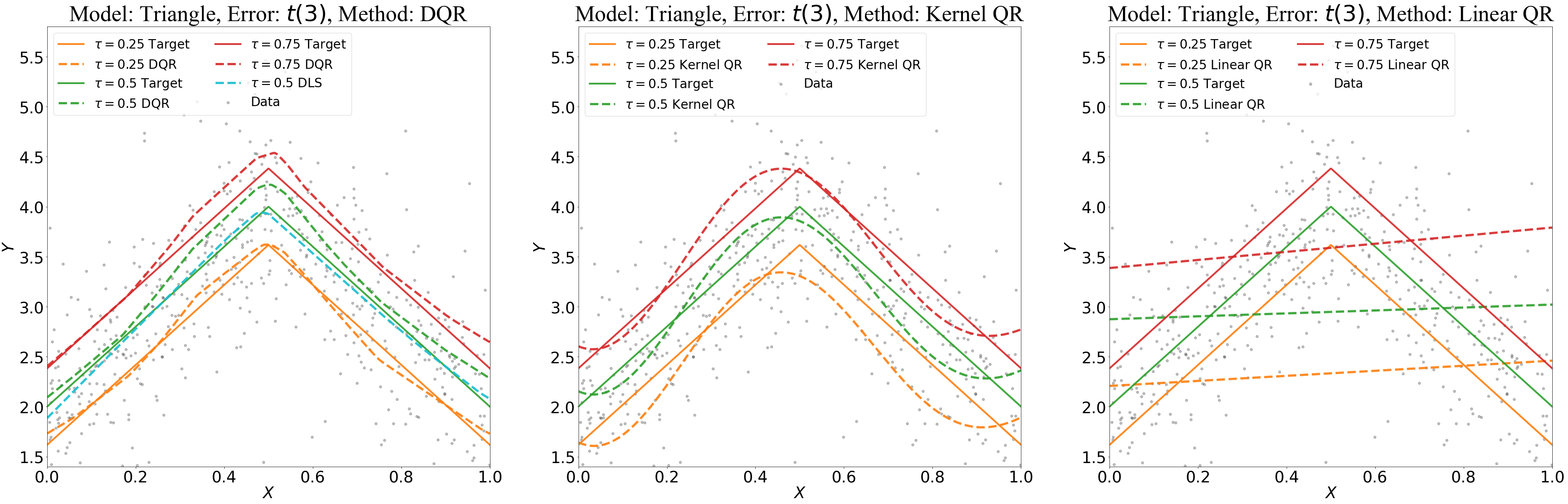}
	\end{subfigure}
	
	\begin{subfigure}{\textwidth}
		\includegraphics[width=1\textwidth]{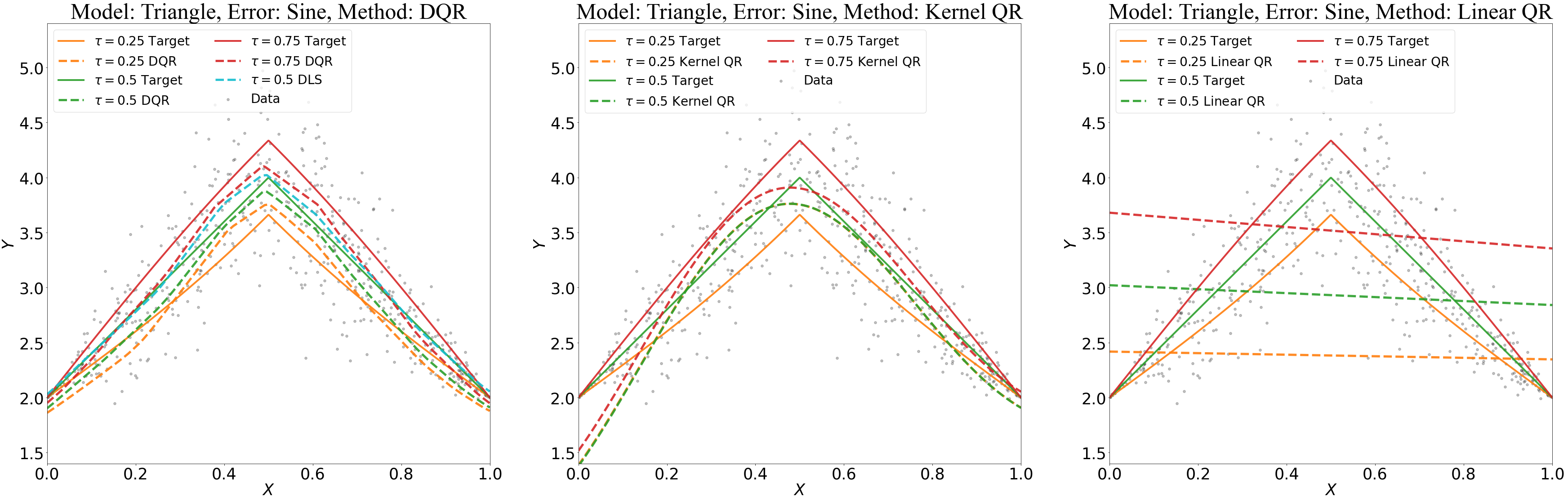}
	\end{subfigure}
	
	\begin{subfigure}{\textwidth}
		\includegraphics[width=\textwidth]{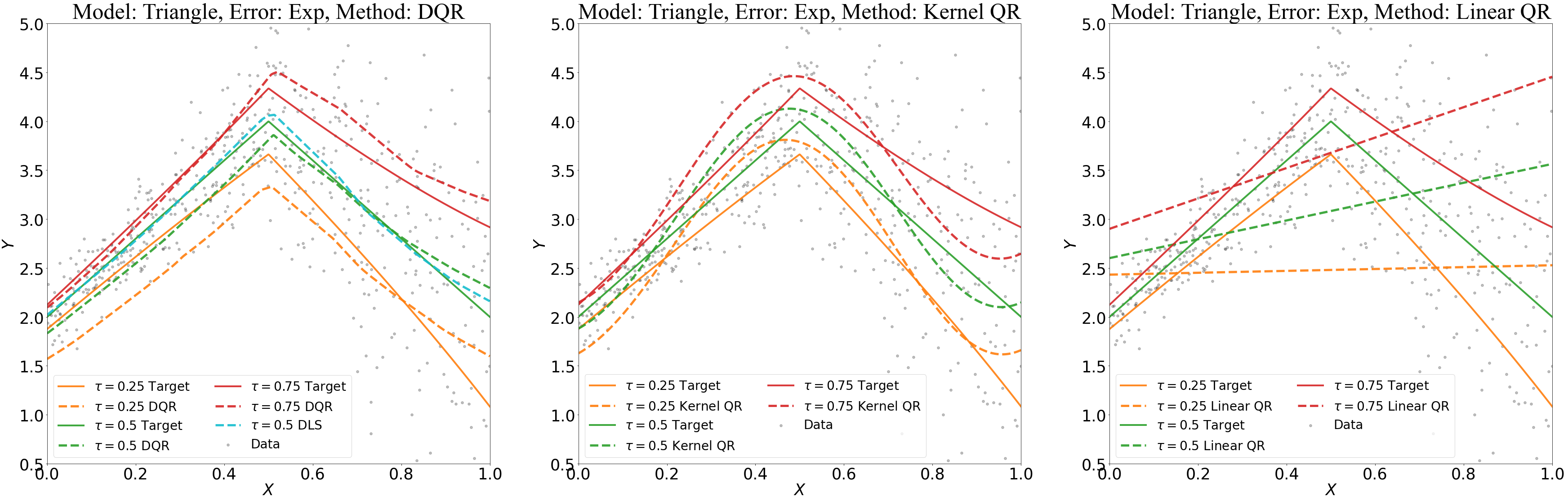}
	\end{subfigure}
\caption{The fitted quantile curves by different methods under the univariate model ``Triangle'' with different errors. The training data is depicted as grey dots.The target quantile functions at the quantile levels $\tau=$0.25 (yellow), 0.5 (green), 0.75 (red)  are depicted as solid curves, and the estimated quantile functions are represented by dashed curves with the same color.
From the top to the bottom, the rows correspond to the errors $t(3)$, ``\textit{Sine}'' and ``\textit{Exp}''.
From the left to the right, the subfigures correspond to the methods \textit{DQR}, \textit{kernel QR} and \textit{linear QR}. The fitted \textit{DLS} curve (in blue)  is included in the left subfigure.}
	\label{fitted:triangle1}
\end{figure}

We generate training data with sample sizes $n=128$ and set the batch size of Adam optimization to be $n/2$. In all settings, we implement the empirical risk minimization of \textit{DQR} and DLS by ReLU activated fixed width multilayer perceptrons, i.e., a class of 
ReLU activated multilayer perceptrons with 4 hidden layers, the width of the network are set to be $(1,256,256,256,256,1)$. All weights and biases in each layer are initialized by uniformly samples on bounded intervals according to the default initialization mechanism in \textit{PyTorch}. The fitted quantiles curves at $\tau=0.25,0.5,0.75$ are shown in Figures \ref{fitted:linear1}-\ref{fitted:triangle1}.
Summary measures including  the excess risks and the $L_1$
test errors are summarized in Tables \ref{tab:linear}-\ref{tab:triangle}.

Additional simulation results with $n=512$, including the estimated quantile curves at $\tau=$ 0.05, 0.25,0.5,0.75 and 0.95,
the corresponding excess risks, the $L_1$ and the $L_2$ test errors
are given in
Appendix \ref{AppFigures}.

It can be seen that for ``Linear'' model, the traditional \textit{linear QR}  works fine  but
 it does poorly in nonlinear models, e.g.,
in the  ``Wave'' and the ``Triangle'' models. This is not surprising since the linear model is misspecified here.  \textit{kernel QR} works reasonably well
in the three models considered, but has difficulty in fitting very winding or nonsmooth curves.
\textit{DQR} tends to perform better than \textit{kernel QR} across all the settings. In particular, \textit{DQR} successfully fits very winding and nonsmooth curves.
The performance of \textit{DLS} is similar to that of \textit{DQR} at the $0.5$-th quantile.

\begin{table}[H]
	\setlength{\tabcolsep}{0.5pt} 
	\renewcommand{\arraystretch}{1.0} 
	\centering
	\caption{Data is generated from ``Linear" model with training sample size $n = 128$ and the number of replications $R = 10$. The averaged excess risks and the $L_1$ test errors with the target with the corresponding standard deviations (in parentheses) are reported for the estimators trained by different methods.}
	\label{tab:linear}
	\resizebox{\textwidth}{!}{%
{\scriptsize
		\begin{tabular}{@{}ll|rr|rr|rr@{}}
			\toprule
			\multicolumn{2}{c|}{$n=128$}              & \multicolumn{2}{c|}{$t(3)$}                                            & \multicolumn{2}{c|}{\textit{Sine}}                                           & \multicolumn{2}{c}{\textit{Exp}}                                            \\
			Quantile                     & Method    & Excess risk           & $L_1$ error
           & Excess risk           & $L_1$ error
        & Excess risk           & $L_1$   error                         \\ \midrule
			\multirow{3}{*}{$\tau=0.25$}
			& DQR & 0.06(0.03)          & 0.31(0.09)
     & 0.02(0.01)          & 0.16(0.03)
             & 0.04(0.03)          & 0.30(0.07)          \\
			& Kernel QR & 0.04(0.03)          & 0.26(0.08)
       & 0.08(0.03)          & 0.32(0.06)
               & 0.01(0.01)          & 0.17(0.06)                \\
			& Linear QR & \textbf{0.01(0.01)} & \textbf{0.08(0.04)}
& \textbf{0.01(0.01)} & \textbf{0.11(0.02)}
 & \textbf{0.01(0.02)} & \textbf{0.09(0.04)} \\ \midrule
			\multirow{4}{*}{$\tau=0.5$}  & DLS       & 0.28(0.13)          & 0.37(0.08)
       & 0.06(0.04)          & 0.15(0.03)
               & 0.18(0.07)          & 0.27(0.05)                \\
			& DQR & 0.10(0.04)          & 0.38(0.11)
     & 0.02(0.01)          & 0.16(0.04)
           & 0.05(0.03)          & 0.28(0.05)                \\
			& Kernel QR & 0.03(0.01)          & 0.23(0.08)
      & 0.06(0.04)          & 0.22(0.10)
             & 0.02(0.01)          & 0.17(0.04)                \\
			& Linear QR & \textbf{0.01(0.01)} & \textbf{0.07(0.05)}
 & \textbf{0.01(0.01)} & \textbf{0.02(0.02)}
  & \textbf{0.01(0.01)} & \textbf{0.07(0.04)} \\ \midrule
			\multirow{3}{*}{$\tau=0.75$}
			& DQR & 0.08(0.05)          & 0.39(0.11)
       & 0.01(0.01)          & 0.20(0.05)
             & 0.05(0.02)          & 0.33(0.05)               \\
			& Kernel QR & 0.01(0.01)          & 0.20(0.08)
       & 0.05(0.03)          & 0.32(0.13)                 & 0.03(0.03)          & 0.20(0.08)         \\
			& Linear QR & \textbf{0.01(0.01)} & \textbf{0.09(0.04)}
& \textbf{0.01(0.01)} & \textbf{0.11(0.01)}
 & \textbf{0.01(0.01)} & \textbf{0.12(0.06)} \\
 \bottomrule
		\end{tabular}%
	}
}
\end{table}

\bigskip
\begin{table}[H]
	\setlength{\tabcolsep}{0.5pt} 
	\renewcommand{\arraystretch}{1.0} 
	\centering
	\caption{Data is generated from ``Wave" model with training sample size $n = 128$ and the number of replications $R = 10$. The averaged excess risks and the
$L_1$ test errors with the corresponding standard deviation (in parentheses)
are reported for the estimators trained by different methods.}
	\label{tab:wave}
	\resizebox{\textwidth}{!}{%
{\scriptsize
		\begin{tabular}{@{}ll|rr|rr|rr@{}}
			\toprule
			\multicolumn{2}{c|}{$n=128$}              & \multicolumn{2}{c|}{$t(3)$}                                            & \multicolumn{2}{c|}{\textit{Sine}}                         & \multicolumn{2}{c}{\textit{Exp}}       \\
			Quantile                     & Method    & Excess risk           & $L_1$  error
               & Excess risk           & $L_1$  error
               & Excess risk           & $L_1$  error               \\ \midrule
			\multirow{3}{*}{$\tau=0.25$}
			& DQR & \textbf{0.07(0.04)} & \textbf{0.34(0.07)}
& \textbf{0.02(0.01)} & \textbf{0.16(0.03)}
& \textbf{0.05(0.03)} & \textbf{0.35(0.06)} \\
	& Kernel QR & 0.13(0.01)          & 0.51(0.02)
& 0.17(0.01)          & 0.52(0.02)
& 0.13(0.02)          & 0.53(0.02)               \\
			& Linear QR & 0.25(0.02)          & 0.61(0.02)
   & 0.25(0.01)          & 0.61(0.01)
         & 0.23(0.04)          & 0.61(0.02)                \\ \midrule
\multirow{4}{*}{$\tau=0.5$}  & DLS       & 0.20(0.06)          & \textbf{0.33(0.05)}
& 0.05(0.02)          & \textbf{0.15(0.03)}
& 0.21(0.06)          & 0.30(0.05)             \\
			& DQR & \textbf{0.10(0.05)} & 0.35(0.10)
  & \textbf{0.02(0.01)} & 0.18(0.02)
       & \textbf{0.05(0.02)} & \textbf{0.29(0.06)}  \\
			& Kernel QR & 0.15(0.02)          & 0.50(0.04)
       & 0.17(0.01)          & 0.52(0.01)
             & 0.16(0.03)          & 0.53(0.02)           \\
			& Linear QR & 0.25(0.02)          & 0.60(0.02)
 & 0.28(0.02)          & 0.58(0.01)
     & 0.21(0.03)          & 0.60(0.02)             \\ \midrule
			\multirow{3}{*}{$\tau=0.75$}
			& DQR & \textbf{0.09(0.04)} & \textbf{0.44(0.13)}
    & \textbf{0.01(0.01)} & \textbf{0.20(0.03)}
     & \textbf{0.07(0.03)} & \textbf{0.35(0.06)}  \\
			& Kernel QR & 0.10(0.02)          & 0.52(0.02)
& 0.13(0.02)          & 0.52(0.02)
    & 0.13(0.03)          & 0.52(0.02)             \\
			& Linear QR & 0.14(0.01)          & 0.68(0.04)
  & 0.18(0.01)          & 0.76(0.05)
      & 0.12(0.01)          & 0.63(0.03)    \\
      \bottomrule
		\end{tabular}
	}
}
\end{table}

\begin{table}[H]
	\setlength{\tabcolsep}{0.5pt} 
	\renewcommand{\arraystretch}{1.0} 
	\centering
	\caption{ Data is generated from ``Triangle" model with training sample size $n = 128$  and the number of replications $R = 10$. The averaged excess risks and the $L_1$ test errors with the corresponding standard deviation (in parentheses) are reported for the estimators trained by different methods.}
	\label{tab:triangle}
	\resizebox{\textwidth}{!}{%
{\scriptsize
		\begin{tabular}{@{}ll|rr|rr|rr@{}}
			\toprule
			\multicolumn{2}{c|}{$n=128$}              & \multicolumn{2}{c|}{$t(3)$}                                            & \multicolumn{2}{c|}{\textit{Sine}}                                            & \multicolumn{2}{c}{\textit{Exp}}                                             \\
			Quantile                     & Method    & Excess risk           & $L_1$ error
           & Excess risk           & $L_1$ error
                   & Excess risk           & $L_1$  error                          \\ \midrule
			\multirow{3}{*}{$\tau=0.25$}
			& DQR & 0.05(0.03)          & 0.27(0.07)
       & \textbf{0.01(0.01)} & \textbf{0.12(0.03)}
        & 0.02(0.02)          & 0.25(0.06)                 \\
			& Kernel QR & \textbf{0.04(0.03)} & \textbf{0.23(0.09)}
& 0.10(0.05)          & 0.36(0.07)
       & \textbf{0.01(0.01)} & \textbf{0.20(0.05)} \\
			& Linear QR & 0.17(0.02)          & 0.55(0.03)
     & 0.17(0.02)          & 0.50(0.04)
           & 0.13(0.01)          & 0.59(0.02)                \\ \midrule
			\multirow{4}{*}{$\tau=0.5$}  & DLS       & 0.16(0.08)          & 0.29(0.07)
        & 0.02(0.02)          & \textbf{0.11(0.03)}
        & 0.11(0.05)          & 0.21(0.06)              \\
			& DQR & 0.06(0.03)          & 0.27(0.09)
        & \textbf{0.01(0.01)} & 0.15(0.03)
                & 0.07(0.05)          & 0.30(0.07)                \\
			& Kernel QR & \textbf{0.03(0.03)} & \textbf{0.20(0.10)}
& 0.05(0.03)          & 0.24(0.08)
         & \textbf{0.03(0.02)} & \textbf{0.19(0.07)} \\
			& Linear QR & 0.14(0.01)          & 0.51(0.01)
       & 0.19(0.01)          & 0.52(0.01)
& 0.11(0.01)          & 0.52(0.02)                 \\ \midrule
			\multirow{3}{*}{$\tau=0.75$}
			& DQR & 0.07(0.04)          & 0.38(0.10)
    & \textbf{0.01(0.01)} & \textbf{0.16(0.04)}
     & 0.04(0.02)          & 0.31(0.08)               \\
			& Kernel QR & \textbf{0.03(0.03)} & \textbf{0.23(0.11)}
 & 0.04(0.02)          & 0.26(0.09)
           & \textbf{0.03(0.01)} & \textbf{0.18(0.05)}  \\
			& Linear QR & 0.08(0.01)          & 0.53(0.02)
       & 0.14(0.01)          & 0.64(0.03)
& 0.07(0.01)          & 0.51(0.02)
         \\ \bottomrule
		\end{tabular}%
	}
}
\end{table}

\subsection{Data generation: multivariate models}
Throughout the multivariate model simulation, we set the input dimension $d=6$ and sample $X$ uniformly on $[0,1]^6$. We consider the models in Section \ref{compositionR} including single index model and additive model 
which correspond different specifications of $f_0$.  The formulae of are given below.
\begin{enumerate}[(a)]
	\item Single index model:
	$$f_0(x)=\exp(\theta^\top x),$$
	where $\theta=(2.2831,-1.4818,5.1966,0,0,0.0515)^\top\in\mathbb{R}^6$.
	\item Additive model:
	$$f_0(x)=\exp(4(x_1-0.5))+9(x_2-0.5)^2+10\sin(2\pi x_3)-7\vert x_4-0.5\vert,$$
	where $x=(x_1,\ldots,x_6)^\top\in[0,1]^6.$
\end{enumerate}
And we generate the error $\eta$ from following distributions,
\begin{enumerate}[(i)]
	\item $\eta$ follows a scaled Student's t distribution with degree of freedom 3, i.e., $\eta\sim 0.5\times t(3)$, denoted by $t(3)$;
	\item  Conditioning on $X=x$, the error $\eta$ follows a normal distribution of which variance depends on the covariate $X$, denoted by \textit{Sine}, i.e.,
	$$\eta\mid X=x\sim 0.5\times\mathcal{N}(0,\vert\sin(\pi \xi^\top x)\vert^2)$$
	where $\xi=(1.8100,-1.2999,0,0,-2.7874,0.3197)^\top\in\mathbb{R}^d$;
	\item Conditioning on $X=x$, the error $\eta$ follows a normal distribution of which variance depends on the covariate $X$, denoted by \textit{Exp}, i.e., $$\eta\mid X=x\sim0.5\times\mathcal{N}(0,\exp(4\xi^\top x-2))$$
	where $\xi=(1.8100,-1.2999,0,0,-2.7874,0.3197)^\top\in\mathbb{R}^d$.
\end{enumerate}
Similarly, the $\tau$-th conditional quantile $f_0^\tau(x)$ of response $Y$ given $X=x$ can be calculated by
$$f_0^\tau(x)=f_0(x)+F^{-1}_{\eta\mid X=x}(\tau),$$
where $F^{-1}_{\eta\mid X=x}(\cdot)$ is the inverse of the conditional cumulated distribution function of $\eta$ given $X=x$.

We generate training data with sample size $n=512$ and train the estimators in the same way as in the  univariate model simulations.  Summary measures including the excess risks and the $L_1$ test errors
based on $R=10$ replications are summarized in Tables \ref{tab:sim}-\ref{tab:add}.
Additional simulation results with $n=1024$, including the estimated quantile curves at $\tau=$ 0.05, 0.25,0.5,0.75 and 0.95,
the corresponding excess risks, the $L_1$ and the $L_2$ test errors
are given in Appendix \ref{AppFigures}.

We see that for the nonlinear multivariate models considered in the simulation studies, especially for single index model, \textit{DQR} performs significantly better than \textit{kernel QR} and
\textit{linear QR},
in the sense that \textit{DQR} estimates have smaller excess risks and $L_1$ test errors in all the scenarios.


\begin{table}[H]
	\setlength{\tabcolsep}{0.5pt} 
	\renewcommand{\arraystretch}{1.0} 
	\centering
	\caption{Data is generated from single index model with training sample size $n = 512$ and the number of replications $R = 10$. The averaged excess risks and $L_1$ test errors with the corresponding standard deviation (in parentheses) are reported for the estimators trained by different methods.}
	\label{tab:sim}
	\resizebox{1.0\textwidth}{!}{%
{\scriptsize
		\begin{tabular}{@{}ll|rr|rr|rr@{}}
			\toprule
			\multicolumn{2}{c|}{$n=512$}              & \multicolumn{2}{c|}{$t(3)$}                                                    & \multicolumn{2}{c|}{\textit{Sine}}                                                   & \multicolumn{2}{c}{\textit{Exp}}                                                        \\
			Quantile                     & Method    & Excess risk             & $L_1$  error         & Excess risk             & $L_1$   error                     & Excess risk             & $L_1$   error                          \\ \midrule
			\multirow{3}{*}{$\tau=0.25$}
			& DQR & \textbf{17.05(12.02)} & \textbf{2.50(0.78)} &
 \textbf{17.71(11.38)} & \textbf{2.08(0.45)}   &
 \textbf{17.66(10.83)} & \textbf{2.19(0.45)}       \\
			& Kernel QR & 1299.84(98.46)        & 26.37(0.39)
    & 1301.01(98.12)        & 26.40(0.37)
       & 1301.02(98.44)        & 26.444(0.38)                \\
			& Linear QR & 3406.75(88.87)        & 47.584(0.45)
      & 3408.65(80.14)        & 47.70(0.42)
        & 3402.89(84.86)        & 47.76(0.43)                \\ \midrule
			\multirow{4}{*}{$\tau=0.5$}  & DLS
& 97.952(46.86)          & \textbf{2.07(0.24)}
    & 98.78(38.52)          & \textbf{2.27(1.33)}
      & 87.27(26.93)          & \textbf{1.79(0.18)}       \\
			& DQR & \textbf{31.48(30.99)} & 6.08(3.39)
       & \textbf{24.20(22.29)} & 5.043(3.09)
              & \textbf{33.68(26.67)} & 4.26(2.69)                      \\
			& Kernel QR & 2358.61(213.79)       & 24.26(0.35)
          & 2362.17(213.94)       & 24.26(0.36)
             & 2363.70(213.50)       & 24.25(0.36)                   \\
			& Linear QR & 5664.87(282.83)       & 44.82(0.25)
       & 5669.25(287.26)       & 44.832(0.25)
            & 5667.21(289.23)       & 44.84(0.26)                  \\ \midrule
			\multirow{3}{*}{$\tau=0.75$}
			& DQR & \textbf{49.37(34.46)} & \textbf{5.59(4.87)}
  & \textbf{44.42(42.23)} & \textbf{3.007(1.90)}
    & \textbf{27.52(27.79)} & \textbf{9.20(4.74)}    \\
			& Kernel QR & 3293.03(311.42)       & 26.02(0.34)
         & 3298.69(308.74)       & 26.10(0.33)
             & 3299.32(308.81)       & 26.168(0.33)                 \\
			& Linear QR & 5410.39(496.13)       & 58.366(3.08)
         & 5419.965(499.63)       & 58.30(3.02)
           & 5422.32(499.93)       & 58.336(2.99)                \\
           \bottomrule
		\end{tabular}
	}
}
\end{table}

	\begin{table}[H]
		\setlength{\tabcolsep}{2pt} 
		\renewcommand{\arraystretch}{1.0} 
		\centering
		\caption{Data is generated from additive model with training sample size $n=512$ and the number of replications $R = 10$. The averaged excess risks and the $L_1$
test errors with the corresponding standard deviation (in parentheses) are reported for the estimators trained by different methods.}
		\label{tab:add}
		\resizebox{1.0\textwidth}{!}{%
{\scriptsize
			\begin{tabular}{@{}ll|rr|rr|rr@{}}
				\toprule
				\multicolumn{2}{c|}{$n=512$}              & \multicolumn{2}{c|}{$t(3)$}                                            & \multicolumn{2}{c|}{\textit{Sine}}               & \multicolumn{2}{c}{\textit{Exp}}           \\
		Quantile            & Method    & Excess risk           & $L_1$  error
& Excess risk           & $L_1$    error
& Excess risk           & $L_1$  error                   \\ \midrule
\multirow{4}{*}{$\tau=0.25$}
				& DQR & \textbf{0.28(0.04)} & \textbf{0.75(0.05)}
 & \textbf{0.14(0.03)} & \textbf{0.44(0.05)}
  & \textbf{0.12(0.03)} & \textbf{0.44(0.09)}\\
				& Kernel QR & 4.48(0.29)          & 3.291(0.06)
         & 4.218(0.23)          & 3.25(0.07)
             & 4.55(0.27)          & 3.40(0.07)             \\
				& Linear QR & 9.20(0.78)          & 4.79(0.17)
    & 8.97(0.40)          & 4.79(0.09)
           & 9.51(0.80)          & 4.96(0.19)               \\ \midrule
				\multirow{4}{*}{$\tau=0.5$}  & DLS       & 0.93(0.15)          & 0.72(0.05)
     & 0.28(0.04)          & \textbf{0.40(0.03)}
      & 0.261(0.07)          & \textbf{0.35(0.05)}  \\
				& DQR & \textbf{0.35(0.08)} & \textbf{0.72(0.06)}
 & \textbf{0.16(0.03)} & 0.45(0.03)
     & \textbf{0.16(0.05)} & 0.39(0.07)                  \\
				& Kernel QR & 3.63(0.49)          & 2.90(0.04)
       & 3.21(0.39)          & 2.85(0.04)
          & 3.50(0.47)          & 2.88(0.04)               \\
				& Linear QR & 7.04(0.78)          & 4.03(0.04)
        & 6.52(0.63)          & 4.04(0.03)
             & 7.17(0.74)          & 4.04(0.03)            \\ \midrule
				\multirow{4}{*}{$\tau=0.75$}
				& DQR & \textbf{0.45(0.08)} & \textbf{0.80(0.06)}
 & \textbf{0.18(0.04)} & \textbf{0.(0.047)}
  & \textbf{0.18(0.04)} & \textbf{0.41(0.07)}  \\
				& Kernel QR & 1.58(0.26)          & 3.21(0.07)
    & 1.41(0.16)          & 3.30(0.09)
        & 1.63(0.23)          & 3.31(0.09)            \\
				& Linear QR & 2.47(0.28)          & 4.69(0.11)
    & 2.44(0.32)          & 4.88(0.18)
            & 2.56(0.27)          & 4.84(0.11)             \\
     \bottomrule
			\end{tabular}%
		}}
	\end{table}

\section{Related work}
\label{related}
{\color{black}

There were several important early works on nonparametric quantile regression using neural networks. \citet{white1992} established the consistency of nonparametric conditional quantile estimators using shallow neural networks.
\citet{chen1999} obtained convergence rate in the Sobolev norm for a large class of single hidden layer feedforward neural networks with a smooth activation functions, assuming the target function satisfies
certain smoothness conditions.
\citet{chen2020efficient} considered quantile treatment effect estimation and established asymptotic distributional properties for the treatment effect estimator in the presence of a infinite-dimensional parameter that is estimated using deep neural networks.  In this semiparametric framework, to establish the asymptotic normality of a finite-dimensional parameter, it is necessary to derive the convergence rate
of the infinite-dimensional nuisance parameter.

Recently,
\citet{padilla2020quantile} studied the nonparametric quantile regression with  ReLU neural networks. They established an upper bound on the mean integrated squared error
of the empirical risk minimizer.
As a consequence, they derived a nearly optimal error bound
when the target quantile function is a composed of H\"older smooth functions. They also
derived a minimax nonparametric estimation rate with Gaussian errors when the target quantile regression function belongs to a Besov space without a compositional structure.
Their approach follows the method of  \cite{schmidt2020nonparametric}, which studied the
least squares nonparametric regression using ReLU neural networks to approximate the regression function. In particular, for approximating a composite function, \cite{padilla2020quantile}  used the approximation results from \cite{schmidt2020nonparametric}. Therefore, the error bounds obtained
by \citet{padilla2020quantile} are similar to the results of \citet{schmidt2020nonparametric}.
In particular,
the prefactor of their error bounds is of the order $O(2^d)$ unless the size $\mathcal{S}$ of the network grows exponentially with respect to the dimension $d$.
A prefactor of the order $O(2^d)$ is big even for a moderate $d$, which can dominate the error bound.

Another important difference between \citet{padilla2020quantile} and our work concerns the neural networks used in constructing the estimators. In  \citet{padilla2020quantile},
they assume that all the parameters (weights and biases) of the network are bounded by one and the
networks are sparse as in \citet{schmidt2020nonparametric}. We do not make such assumptions.
We note that such assumptions are usually not satisfied in training neural network models in practice.


A unique aspect of the quantile loss is that a bound on the excess risk does not automatically lead to a bound for the mean squared error of the estimated quantile regression function. This is different from the
squared loss whose excess risk bound directly leads to a bound on the mean squared error of the estimated regression function.
In \cite{steinwart2011estimating}, under the $\tau$-quantile of $p$-average type  condition on the joint distribution of $(X,Y)$, a general result is given: the $L^r(\nu)$ distance ($\nu$ denotes the distribution of the predictor) between any function $f$ and the target $f_0$ can be bound by the $q$-th root of the excess risk for some $r,q>0$.  This problem was also considered in  \cite{christmann2007svms, lv2018oracle, padilla2020quantile} and \cite{padilla2021adaptive}.
However, these existing results require that the conditional distribution of $Y$ given $X$ is bounded, which does not apply to our setting where we allow  the response to have heavy tails.
}

There are several recent important studies on least squares nonparametric regression using deep neural networks. Examples include
\citet{
bauer2019deep,
chen2019nonparametric,
nakada2019adaptive,
schmidt2019deep,
kohler2019estimation}
and \citet{farrell2021deep}.
In particular,
\citet{bauer2019deep} assumed that the activation function satisfies certain smoothness conditions, which excludes the use of ReLU activation; \citet{schmidt2020nonparametric} and \citet{farrell2021deep} considered the ReLU activation function. \citet{bauer2019deep} and \citet{schmidt2020nonparametric} assumed that the regression function has a compositional structure.
These studies
adopt a construction of function approximation using deep neural networks
similar to that of \cite{yarotsky2017error},
which will lead to a prefactor depending on the dimension $d$ exponentially.
For a large $d$, a prefactor that depends on $d$ exponentially will severely deteriorate the
quality of the error bound.
In comparison, the prefactor in the error bounds in our work has a polynomial dependence on $d$. Therefore, there is a significant improvement in our results in terms of mitigating the curse of dimensionality.

Finally, we should mention that there have been a great deal of efforts  to
deal with the curse of dimensionality by assuming that the distribution of the predictor   is supported on a lower dimensional manifold. Many methods have been developed under this condition, including
 local regression \citep{bickel2007local,cheng2013local,aswani2011regression},  kernel methods \citep{kpotufe2013adaptivity},
Gaussian process regression \citep{yang2016bayesian}, and deep neural networks \citep{nakada2019adaptive,schmidt2019deep,chen2019efficient,chen2019nonparametric,
kohler2019estimation, farrell2021deep,jiao2021deep}.
Several studies have focused on representing the data on the manifold itself, e.g.,  manifold learning or dimensionality reduction \citep{pelletier2005kernel,hendriks1990nonparametric,tenenbaum2000global,
donoho2003hessian,belkin2003laplacian,lee2007nonlinear}.
If a high-dimensional data vector can be well represented by a lower-dimensional feature,
 the problem of curse of dimensionality can be attenuated.

\section{Conclusion}
\label{conclusion}
In recent years, there have been intensive efforts devoted to understanding the properties of deep neural network modeling by researchers from various fields, including applied mathematics, machine learning,  and statistics. In particular, much work has been done to study the properties of the least squares nonparametric regression estimators using deep neural networks. This line of work showed that a key factor for the success of deep neural network modeling is its ability to accurately and adaptively approximate high-dimensional functions.
Indeed, although neural networks models had been developed many years ago and it had been shown that they can serve as universal approximators to multivariate functions, only recently the advantages of deep networks over shallow networks in approximating high-dimensional functions were clearly demonstrated.


In this work, we study the convergence properties of nonparametric quantile regression using deep neural networks. To mitigate the curse of dimensionality, we assume that the target quantile regression function has a compositional structure.  Based on the recent results on the approximation power of deep neural networks, we show that composite functions can be well approximated by neural networks with error rate determined by the intrinsic dimension of the function, not the ambient dimension.
We established non-asymptotic bounds for the excess risk of deep quantile regression and the mean squared error of the estimated quantile regression function. We explicitly describe how these bounds depend on the network parameters (e.g., depth and width), the intrinsic dimension and the ambient dimension. Our error bounds significantly improve over the existing ones in the sense that their prefactors depend linearly or quadratically on the ambient dimension $d$, instead of exponentially on $d$. We also provide explicit error bounds, including the prefactors, for several well-known semiparametric and nonparametric regression models that have been widely used to mitigate the curse of dimensionality.

Our results are obtained based on the key assumption that the conditional quantile function has a compositional structure. This assumption provides an effective way for mitigating the curse of dimensionality in nonparametric estimation problems. In the future work, it would be interesting to also consider other conditions that can help lessen the curse of dimensionality, such as the low-dimensional support assumption for the predictor that has been used in the context of least squares regression.
Another problem that deserves further study is to generalize the results in this work to the setting with a general convex losses, including robust loss functions, and other regression problems such as nonparametric Cox regression.  We hope to study these problems in the future.

%

\section*{Acknowledgements}
The work of Y. Jiao is supported in part by the National Science Foundation of China grant 11871474 and by the research fund of KLATASDSMOE of China.
The work of Y. Lin is supported by the Hong Kong Research Grants Council (Grant No.
14306219 and 14306620) and Direct Grants for Research,
The Chinese University of Hong Kong.
The work of J. Huang is partially supported by the U.S. National Science Foundation grant DMS-1916199.
	
\bibliographystyle{apalike}
\bibliography{dqr_bib.bib}    

\clearpage
\appendix
\numberwithin{equation}{section}
\makeatletter
\newcommand{\section@cntformat}{Appendix \thesection:\ }

\setcounter{table}{0}
\setcounter{figure}{0}
\setcounter{equation}{0}
\renewcommand{\thetable}{B.\arabic{table}}
\renewcommand{\thefigure}{B.\arabic{figure}}

\begin{center}
\begin{Large}
\textbf{Appendix}
\end{Large}
\end{center}

In the appendix, we give the proofs of the theoretical results in the paper and provide additional simulation results.

\section{Appendix: Proofs}

\subsection{Proof of Lemma \ref{lemma1}}
\begin{proof}
	By the definition of the empirical risk minimizer, for any $f\in\mathcal{F}_n$, we have $\mathcal{R}^\tau_n(\hat{f}_n)\leq \mathcal{R}^\tau_n(f)$. Therefore,
	\begin{align*} \mathcal{R}^\tau(\hat{f}_n)-\mathcal{R}^\tau(f_0)=&\mathcal{R}^\tau(\hat{f}_n)-\mathcal{R}^\tau_n(\hat{f}_n)+\mathcal{R}^\tau_n(\hat{f}_n)-\mathcal{R}^\tau_n(f)+\mathcal{R}^\tau_n(f)-\mathcal{R}^\tau(f)+\mathcal{R}^\tau(f)-\mathcal{R}^\tau(f_0)\\
		\leq&\mathcal{R}^\tau(\hat{f}_n)-\mathcal{R}^\tau_n(\hat{f}_n)+\mathcal{R}^\tau_n(f)-\mathcal{R}^\tau(f)+\mathcal{R}^\tau(f)-\mathcal{R}^\tau(f_0)\\
		=&\big\{\mathcal{R}^\tau(\hat{f}_n)-\mathcal{R}^\tau_n(\hat{f}_n)\big\}+\big\{\mathcal{R}^\tau_n(f)-\mathcal{R}^\tau(f)\big\}+\big\{\mathcal{R}^\tau(f)-\mathcal{R}^\tau(f_0)\big\}\\
		\leq& 2\sup_{f\in\mathcal{F}_n}\vert \mathcal{R}^\tau(f)-\mathcal{R}^\tau_n(f)\vert+\big\{\mathcal{R}^\tau(f)-\mathcal{R}^\tau(f_0)\big\}.
	\end{align*}
	Since the above inequality holds for any $f\in\mathcal{F}_n$, Lemma \ref{lemma1} is proved by choosing $f$  satisfying $f\in\arg\inf_{f\in\mathcal{F}_n}\mathcal{R}^\tau(f)$.
\end{proof}

\subsection{Proof of Lemma \ref{lemma2}}
\begin{proof}
	Let $S=\{Z_i=(X_i,Y_i)\}_{i=1}^n$ be a sample form the distribution of $Z=(X,Y)$ and $S^\prime=\{Z_i^\prime=(X^\prime_i,Y^\prime_i)\}_{i=1}^n$ be another sample independent with $S$. Define $g(f,Z_i)=\rho_\tau(f(X_i)-Y_i)-\rho_\tau(f_0(X_i)-Y_i)$ for any $f$ and sample $Z_i$. Note that the empirical risk minimizer $\hat{f}_\phi$ defined in Lemma \ref{lemma1} depends on the sample $S$, and its excess risk is $\mathbb{E}_{S^\prime} \{\sum_{i=1}^ng(\hat{f}_\phi,Z_i^\prime)/n\}$ and its prediction error (expected excess risk) is
	\begin{equation}\label{ine1}
		\mathbb{E}\big\{\mathcal{R}^\tau(\hat{f}_\phi)-\mathcal{R}^\tau(f_0)\big\}=\mathbb{E}_{S}[\mathbb{E}_{S^\prime} \{\frac{1}{n}\sum_{i=1}^ng(\hat{f}_\phi,Z_i^\prime)\}].
	\end{equation}
	Next we will take 3 steps to complete the proof of Lemma \ref{lemma2}.
	
	\subsubsection*{Step 1: Prediction error decomposition}
	
	Define the `best in class' estimator $f^*_\phi$ as the  estimator in the function class $\mathcal{F}_\phi=\mathcal{F}_{\mathcal{D},\mathcal{W},\mathcal{U},\mathcal{S},\mathcal{B}}$ with minimal $L$ risk:
	\begin{equation*}
		f^*_\phi=\arg\min_{f\in\mathcal{F}_{\phi}} \mathcal{R}^\tau(f).
	\end{equation*}
	The approximation error of $f^*_\phi$  is $\mathcal{R}^\tau(f^*_\phi)-\mathcal{R}^\tau(f_0)$. Note that the approximation error only depends on the function class $\mathcal{F}_{\mathcal{D},\mathcal{W},\mathcal{U},\mathcal{S},\mathcal{B}}$ and the distribution of data. By the definition of empirical risk minimizer, we have
	\begin{equation}  \label{ine2}
		\mathbb{E}_{S}\{\frac{1}{n}\sum_{i=1}^ng(\hat{f}_\phi,Z_i)\}\leq \mathbb{E}_S\{\frac{1}{n}\sum_{i=1}^ng(f^*_\phi,Z_i)\}.
	\end{equation}
	Multiply 2 by the both sides of (\ref{ine2}) and add it up with (\ref{ine1}), we have
	\begin{align} \nonumber
		\mathbb{E}\big\{\mathcal{R}^\tau(\hat{f}_\phi)-\mathcal{R}^\tau(f_0)\big\}&\leq\mathbb{E}_{S}\Big[ \frac{1}{n}\sum_{i=1}^n\big\{-2g(\hat{f}_\phi,Z_i)+\mathbb{E}_{S^\prime}g(\hat{f}_\phi,Z_i^\prime)\big\}\Big]+ 2\mathbb{E}_S\{\frac{1}{n}\sum_{i=1}^ng(f^*_\phi,Z_i)\}\\
		&\leq\mathbb{E}_{S}\Big[ \frac{1}{n}\sum_{i=1}^n\big\{-2g(\hat{f}_\phi,Z_i)+\mathbb{E}_{S^\prime}g(\hat{f}_\phi,Z_i^\prime)\big\}\Big]+ 2\big\{\mathcal{R}(f^*_\phi)-\mathcal{R}(f^*)\big\}. \label{bound0}
	\end{align}
	It is seen that the prediction error is upper bounded by the sum of a expectation of a stochastic term and approximation error.
	\subsubsection*{Step 2: Bounding the stochastic term}
	Next, we will focus on giving an upper bound of the first term on the right-hand side in (\ref{bound0}), and handle it with truncation and classical chaining technique of empirical process. In the following, for ease of presentation, we write $G(f,Z_i)=\mathbb{E}_{S^\prime}\{g(f,Z_i^\prime)\}-2g(f,Z_i)$ for $f\in\mathcal{F}_\phi$.

	Given a $\delta$-uniform covering of $\mathcal{F}_\phi$, we denote the centers of the balls by $f_j,j=1,2,...,\mathcal{N}_{2n},$ where $\mathcal{N}_{2n}=\mathcal{N}_{2n}(\delta,\Vert\cdot\Vert_\infty,\mathcal{F}_\phi)$ is the uniform covering number with radius $\delta$ ($\delta<\mathcal{B}$) under the norm $\Vert\cdot\Vert_\infty$, where $\mathcal{N}_{2n}(\delta,\Vert\cdot\Vert_\infty,\mathcal{F}_\phi)$
	is defined in (\ref{ucover}).
	By the definition of covering, there exists a (random) $j^*$ such that $\Vert\hat{f}_\phi(x) -f_{j^*}(x)\Vert_\infty\leq\delta$ on $x=(X_1,\ldots,X_n,X_1^\prime,\ldots,X_n^\prime)\in\mathcal{X}^{2n}$, i.e., $\vert \hat{f}_\phi(x) -f_{j^*}(x)\vert\leq\delta$ for all $x\in\{X_1,\ldots,X_n,X_1^\prime,\ldots,X_n^\prime\}$.
	Recall that $g(f,Z_i)=\rho_\tau(f(X_i)-Y_i)-\rho_\tau(f_0(X_i)-Y_i)$ and $\rho_\tau(a)=a(\tau-I(a<0))$. Denote $\lambda_\tau=\max\{\tau,1-\tau\}$, then by the Lipschitz property of $\rho_\tau$, for $a,b\in\mathbb{R}$
	\begin{align*}
		\vert\rho_\tau(a)-\rho_\tau(b)\vert\leq \max\{\tau,1-\tau\}\vert a-b\vert=\lambda_\tau\vert a-b\vert,
	\end{align*}
	and for $ i=1,\ldots,n$
	\begin{gather*}
		\vert g(\hat{f}_\phi,Z_i)-g(f_{j^*},Z_i) \vert \leq \lambda_\tau\delta,\\
		\vert \mathbb{E}_{S^\prime}\{g(\hat{f}_\phi,Z^\prime_i)\}-\mathbb{E}_{S^\prime}\{g(f_{j^*},Z^\prime_i)\}\vert \leq \lambda_\tau\delta.
	\end{gather*}
Then we have,
	\begin{align*}
		\mathbb{E}_{S}\big\{\frac{1}{n}\sum_{i=1}^n g(\hat{f}_\phi,Z_i)\}	\leq\frac{1}{n}\sum_{i=1}^n\mathbb{E}_{S}\big\{&g(f_{j^*},Z_i)\}+\lambda_\tau\delta
	\end{align*}
	and
	\begin{equation} \label{bound1}
		\mathbb{E}_{S}\Big[ \frac{1}{n}\sum_{i=1}^n G(\hat{f}_\phi,Z_i) \Big]
		\leq\mathbb{E}_{S}\Big[ \frac{1}{n}\sum_{i=1}^n G(f_{j^*},Z_i) \Big]+3\lambda_\tau\delta.
	\end{equation}
	
	Let $\beta_n\geq \mathcal{B}\geq1$ be a positive number who may depend on the sample size $n$. Denote $T_{\beta_n}$ as the truncation operator at level $\beta_n$, i.e., for any $Y\in\mathbb{R}$,
	$T_{\beta_n}Y=Y$ if $\vert Y\vert\leq\beta_n$ and $T_{\beta_n}Y= \beta_n\cdot {\rm sign}(Y)$ otherwise.
	Define the function $f^*_{\beta_n}:\mathcal{X}\to\mathbb{R}$ pointwisely by
	$$f^*_{\beta_n}(x)=\arg\min_{f(x):\Vert f\Vert_\infty\leq\beta_n}\mathbb{E}\big\{\rho_\tau(f(X)-T_{\beta_n}Y)|X=x\big\},$$
	for each $x\in\mathcal{X}$.
	Besides, recall that $\Vert f^*\Vert_\infty\leq\mathcal{B}\leq\beta_n$ and
	$$f_0(x)=\arg\min_{f(x):\Vert f\Vert_\infty\leq\beta_n} \mathbb{E}\big\{\rho_\tau(f(X)-Y)|X=x\big\}.$$
	Then for any $f$ satisfying $\Vert f \Vert_\infty\leq\beta_n$, the definition above implies that  $\mathbb{E}\{\rho_\tau(f^*_{\beta_n}(X_i)-T_{\beta_n}Y_i)\}\leq \mathbb{E}\{\rho_\tau(f(X_i)-T_{\beta_n}Y_i)\}$ and $\mathbb{E}\{\rho_\tau(f_0(X_i)-Y_i)\}\leq \mathbb{E}\{\rho_\tau(f(X_i)-Y_i)\}$.
	 For any $f\in\mathcal{F}_\phi$, we let $g_{\beta_n}(f,Z_i)=\rho_\tau(f(X_i)-T_{\beta_n}Y_i)-\rho_\tau(f^*_{\beta_n}(X_i)-T_{\beta_n}Y_i)$. Then we have
	\begin{align*}
		\mathbb{E}\{g(f,Z_i)\} &= \mathbb{E}\{g_{\beta_n}(f,Z_i)\}+ \mathbb{E}\{\rho_\tau(f(X_i),Y_i)-\rho_\tau(f(X_i),T_{\beta_n}Y_i)\}\\
		&\qquad\qquad\qquad\quad+\mathbb{E}\{\rho_\tau(f^*_{\beta_n}(X_i)-T_{\beta_n}Y_i)-\rho_\tau(f^*(X_i)-T_{\beta_n}Y_i)\}\\
		&\qquad\qquad\qquad\quad +\mathbb{E}\{\rho_\tau(f_0(X_i)-T_{\beta_n}Y_i)-\rho_\tau(f_0(X_i)-Y_i)\}\\
		&\leq \mathbb{E}\{g_{\beta_n}(f,Z_i)\}+ \mathbb{E}\vert \rho_\tau(f(X_i)-Y_i)-\rho_\tau(f(X_i)-T_{\beta_n}Y_i)\vert\\
		&\qquad\qquad\qquad\quad+\mathbb{E}\vert \rho_\tau(f_0(X_i)-T_{\beta_n}Y_i)-\rho_\tau(f_0(X_i)-Y_i)\vert\\
		&\leq \mathbb{E}\{g_{\beta_n}(f,Z_i)\} +2\lambda_\tau\mathbb{E}\{\vert T_{\beta_n}Y_i-Y_i\vert\}\\
		&\leq \mathbb{E}\{ g_{\beta_n}(f,Z_i)\}+2\lambda_\tau\mathbb{E}\big\{\vert \vert Y_i\vert I(\vert Y_i\vert>\beta_n)\big\}\\
		&\leq  \mathbb{E}\{ g_{\beta_n}(f,Z_i)\} + 2\lambda_\tau\mathbb{E}\{\vert Y_i\vert \vert Y_i\vert^{p-1}/\beta_n^{p-1}\}\\
		&\leq  \mathbb{E}\{ g_{\beta_n}(f,Z_i)\} + 2\lambda_\tau\mathbb{E}\vert Y_i\vert^{p}/\beta_n^{p-1}.
	\end{align*}
	
	By Assumption 2, the response $Y$ has finite $p$-moment and thus $\mathbb{E}\vert Y_i\vert^{p}<\infty$. Similarly,
	\begin{align*}
		\mathbb{E}\{g_{\beta_n}(f,Z_i)\}&= \mathbb{E}\{g(f,Z_i)\} + \mathbb{E}\{\rho_\tau(f_0(X_i)-Y_i)-\rho_\tau(f^*_{\beta_n}(X_i)-Y_i)\}\\
		&\qquad\qquad\qquad+\mathbb{E}\{\rho_\tau(f(X_i)-T_{\beta_n}Y_i)-\rho_\tau(f(X_i)-Y_i)\}\\
		&\qquad\qquad\qquad +\mathbb{E}\{\rho_\tau(f^*_{\beta_n}(X_i)-Y_i)-\rho_\tau(f^*_{\beta_n}(X_i)-T_{\beta_n}Y_i)\}\\
		&\leq \mathbb{E}\{g(f,Z_i)\} + \mathbb{E}\vert \rho_\tau(f(X_i)-T_{\beta_n}Y_i)-\rho_\tau(f(X_i)-Y_i)\vert\\
		&\qquad\qquad\qquad +\mathbb{E}\vert \rho_\tau(f^*_{\beta_n}(X_i)-Y_i)-\rho_\tau(f^*_{\beta_n}(X_i)-T_{\beta_n}Y_i)\vert\\
		&\leq \mathbb{E}\{ g(f,Z_i)\} + 2\lambda_\tau\mathbb{E}\vert Y_i\vert^{p}/\beta_n^{p-1}.
	\end{align*}
	Note that above inequalities also hold for $g(f,Z_i^\prime)$ and $g_{\beta_n}(f,Z_i^\prime)$.

	By Assumption \ref{moment}, the response $Y$ has finite $p$-moment and thus $\mathbb{E}\vert Y_i\vert^{p}<\infty$.
	Then for any $f\in\mathcal{F}_\phi$, define $G_{\beta_n}(f,Z_i)=\mathbb{E}_{S^\prime}\{g_{\beta_n}(f,Z_i^\prime)\}-2g_{\beta_n}(f,Z_i)$ and we have
	
	\begin{equation} \label{bound2}
		\mathbb{E}_{S}\Big[ \frac{1}{n}\sum_{i=1}^n G(f_{j^*},Z_i) \Big]
		\leq\mathbb{E}_{S}\Big[ \frac{1}{n}\sum_{i=1}^n G_{\beta_n}(f_{j^*},Z_i) \Big]+6\lambda_\tau\mathbb{E}\vert Y_i\vert^{p}/\beta_n^{p-1}.
	\end{equation}

	Besides, by Assumption \ref{moment}, for any $f\in\mathcal{F}_\phi$ we have $\vert g_{\beta_n}(f,Z_i)\vert\leq 4\lambda_\tau\beta_n$ and $\sigma^2_g(f):={\rm Var}(g_{\beta_n}(f,Z_i))\leq\mathbb{E} \{g_{\beta_n}(f,Z_i)^2\}\leq 4\lambda_\tau\beta_n\mathbb{E}\{g_{\beta_n}(f,Z_i)\}$. For each $f_j$ and any $t>0$, let $u=t/2+{\sigma_g^2(f_j)}/(8\lambda_\tau\beta_n)$,
	by applying the Bernstein inequality,
	\begin{align*}
		&P\Big\{\frac{1}{n}\sum_{i=1}^nG_{\beta_n}(f_j,Z_i)>t\Big\}\\
		=&P\Big\{\mathbb{E}_{S^\prime} \{g_{\beta_n}(f_j,Z_i^\prime)\}-\frac{2}{n}\sum_{i=1}^ng_{\beta_n}(f_j,Z_i)>t\Big\}\\
		=&P\Big\{\mathbb{E}_{S^\prime} \{g_{\beta_n}(f_j,Z_i^\prime)\}-\frac{1}{n}\sum_{i=1}^ng_{\beta_n}(f_j,Z_i)>\frac{t}{2}+\frac{1}{2}\mathbb{E}_{S^\prime} \{g_{\beta_n}(f_j,Z_i^\prime)\}\Big\}\\
		\leq& P\Big\{\mathbb{E}_{S^\prime} \{g_{\beta_n}(f_j,Z_i^\prime)\}-\frac{1}{n}\sum_{i=1}^ng_{\beta_n}(f_j,Z_i)>\frac{t}{2}+\frac{1}{2}\frac{\sigma_g^2(f_j)}{4\lambda_\tau\beta_n}\}\Big\}\\
		\leq& \exp\Big( -\frac{nu^2}{2\sigma_g^2(f_j)+16u\lambda_\tau\beta_n/3}\Big)\\
		\leq& \exp\Big( -\frac{nu^2}{16u\lambda_\tau\beta_n+16u\beta_n/3}\Big)\\
		\leq& \exp\Big( -\frac{1}{16+16/3}\cdot\frac{nu}{\lambda_\tau\beta_n}\Big)\\
		\leq& \exp\Big( -\frac{1}{32+32/3}\cdot\frac{nt}{\lambda_\tau\beta_n}\Big).
	\end{align*}
	This leads to a tail probability bound of $\sum_{i=1}^n G_{\beta_n}(f_{j^*},Z_i)/n$, which is
	$$P\Big\{\frac{1}{n}\sum_{i=1}^nG_{\beta_n}(f_{j^*},Z_i)>t\Big\}\leq 2\mathcal{N}_{2n}\exp\Big( -\frac{1}{43}\cdot\frac{nt}{\lambda_\tau\beta_n}\Big).$$
	Then for $a_n>0$,
	\begin{align*}
		\mathbb{E}_S\Big[ \frac{1}{n}\sum_{i=1}^nG_{\beta_n}(f_{j^*},Z_i)\Big]\leq& a_n +\int_{a_n}^\infty P\Big\{\frac{1}{n}\sum_{i=1}^nG_{\beta_n}(f_{j^*},Z_i)>t\Big\} dt\\
		\leq& a_n+ \int_{a_n}^\infty 2\mathcal{N}_{2n}\exp\Big( -\frac{1}{43}\cdot\frac{nt}{\lambda_\tau\beta_n}\Big) dt\\
		\leq& a_n+ 2\mathcal{N}_{2n}\exp\Big( -a_n\cdot\frac{n}{43\lambda_\tau\beta_n}\Big)\frac{43\lambda_\tau\beta_n}{n}.
	\end{align*}
	Choose $a_n=\log(2\mathcal{N}_{2n})\cdot{43\lambda_\tau\beta_n}/{n}$, we have
	\begin{equation} \label{bound3}
		\mathbb{E}_S\Big[ \frac{1}{n}\sum_{i=1}^nG_{\beta_n}(f_{j^*},Z_i)\Big]\leq \frac{43\lambda_\tau\beta_n(\log(2\mathcal{N}_{2n})+1)}{n}.
	\end{equation}
	Set $\delta=1/n$ and $\beta_n=c_1\max\{\mathcal{B},n^{1/p}\}$ and combine (\ref{bound0}), (\ref{bound1}), (\ref{bound2}) and  (\ref{bound3}), we get
	\begin{equation} \label{bound5}
		\mathbb{E}\big\{\mathcal{R}^\tau(\hat{f}_\phi)-\mathcal{R}^\tau(f_0)\big\}\leq\frac{c_2\lambda_\tau\mathcal{B}\log\mathcal{N}_{2n}(\frac{1}{n},\Vert\cdot\Vert_\infty,\mathcal{F}_{\phi})}{n^{1-1/p}}+ 2\big\{\mathcal{R}^\tau(f^*_\phi)-\mathcal{R}^\tau(f_0)\big\},
	\end{equation}
	where $c_2>0$ is a constant does not depend on $n,d,mathcal{B}$ and $\lambda_\tau$. This proves
	(\ref{entropy}).
	
	\subsubsection*{Step 3: Bounding the covering number}
	Lastly, we will give an upper bound on the  covering number by the VC dimension of $\mathcal{F}_\phi$ through its parameters. Denote ${\rm Pdim}(\mathcal{F}_\phi)$ by the pseudo dimension of $\mathcal{F}_\phi$, by Theorem 12.2 in \cite{anthony1999}, for $2n\geq {\rm Pdim}(\mathcal{F}_\phi)$
	$$\mathcal{N}_{2n}(\frac{1}{n},\Vert \cdot\Vert_\infty,\mathcal{F}_\phi)\leq\Big(\frac{2e\mathcal{B}n^2}{{\rm Pdim}(\mathcal{F}_\phi)}\Big)^{{\rm Pdim}(\mathcal{F}_\phi)}.$$
	Besides, based on Theorem 3 and 6 in \cite{bartlett2019nearly}, there exist universal constants $c$, $C$ such that
	$$c\cdot\mathcal{S}\mathcal{D}\log(\mathcal{S}/\mathcal{D})\leq{\rm Pdim}(\mathcal{F}_\phi)\leq C\cdot\mathcal{S}\mathcal{D}\log(\mathcal{S}).$$
	Combine the upper bound of the covering number and pseudo dimension with (\ref{bound5}), we have
	\begin{equation} \label{bound6}
		\mathbb{E}\big\{\mathcal{R}^\tau(\hat{f}_\phi)-\mathcal{R}^\tau(f_0)\big\}\leq c_3\lambda_\tau\mathcal{B}\frac{\log(n)\mathcal{S}\mathcal{D}\log(\mathcal{S})}{n^{1-1/p}}+ 2\big\{\mathcal{R}^\tau(f^*_\phi)-\mathcal{R}^\tau(f_0)\big\},
	\end{equation}
	for some constant $c_3>0$  not dependent on $n,d,\tau,\mathcal{B},\mathcal{S}$ and $\mathcal{D}$. Therefore,  (\ref{oracle}) follows.
	This completes the proof of Lemma \ref{lemma2}.
\end{proof}

\subsection{Proof of Lemma \ref{lemma3}}
Under Assumption \ref{moment}, the function $f_0$ is the risk minimizer. Then for any $f\in\mathcal{F}_\phi$, we have
\begin{align*}
	\mathcal{R}^\tau(f) -\mathcal{R}^\tau(f_0)&=\mathbb{E}\{\rho_\tau(f(X)-Y)-\rho_\tau(f_0(X)-Y)\}\leq \max\{\tau,1-\tau\}\mathbb{E}\{\vert f(X)-f_0(X)\vert\},
\end{align*}
thus $$\inf_{f\in\mathcal{F}_\phi}\{ \mathcal{R}^\tau(f)-\mathcal{R}^\tau(f_0)\}\leq \max\{\tau,1-\tau\}\inf_{f\in\mathcal{F}_{\phi}}\mathbb{E}\vert f(X)-f_0(X)\vert=\max\{\tau,1-\tau\}\inf_{f\in\mathcal{F}_{\phi}}\Vert f-f_0\Vert_{L^1(\nu)},$$
where $\nu$ denotes the marginal probability measure of $X$ and $\mathcal{F}_\phi=\mathcal{F}_{\mathcal{D},\mathcal{W},\mathcal{U},\mathcal{S},\mathcal{B}}$ denotes the class of feedforward neural networks with parameters $\mathcal{D},\mathcal{W},\mathcal{U},\mathcal{S}$ and $\mathcal{B}$.

\subsection{Proof of Lemma \ref{lemma4}}
As in the proof of Lemma \ref{lemma3}, for any $f\in\mathcal{F}_\phi$, we firstly have
\begin{align*}
	\mathcal{R}^\tau(f) -\mathcal{R}^\tau(f_0)&\leq \lambda_\tau\mathbb{E}\{\vert f(X)-f_0(X)\vert\},
\end{align*}
where $\lambda_\tau=\max\{\tau,1-\tau\}$.
Then for function $f\in\mathcal{F}_\phi$ satisfying $\Vert f-f_0\Vert_{L^\infty(\mathcal{X}^0)}>\delta^0_\tau$, we have
\begin{align*}
	\mathcal{R}^\tau(f) -\mathcal{R}^\tau(f_0)&\leq \lambda_\tau\mathbb{E}\{\vert f(X)-f_0(X)\vert\}\\
	&\leq \lambda_\tau\mathbb{E}\big\{\frac{\vert f(X)-f_0(X)\vert^2}{\delta^0_\tau}\big\}\\
	&\leq \frac{\lambda_\tau}{\delta^0_\tau}\Vert f(X)-f_0(X)\Vert^2_{L^2(\nu)}.
\end{align*}
Secondly, with Assumption \ref{quadratic}, we also have
$$\mathcal{R}^\tau(f) -\mathcal{R}^\tau(f_0)\leq c^0_\tau \Vert f-f_0 \Vert^2_{L^2(\nu)},$$
for any $f$ satisfying $\Vert f-f_0 \Vert_{L^\infty(\mathcal{X}^0)}\leq \delta^0_\tau$.
	
There exists a constant $c_\tau\geq \max\{c^0_\tau,\lambda_\tau/\delta^0_\tau\}$ such that
$$\mathcal{R}^\tau(f) -\mathcal{R}^\tau(f_0)\leq c_\tau \Vert f-f_0 \Vert^2_{L^2(\nu)},$$
for any $f\in\mathcal{F}_\phi$, where $\mathcal{X}^0$ is any subset of $\mathcal{X}$ such that $P(X\in\mathcal{X}^0)=P(X\in\mathcal{X})$.

\subsection{Proof of Lemma \ref{lemma6}}
\begin{proof}
	Consider the subnetworks approximating $h_{ij}$ in Lemma \ref{lemma5}, each of them with width $\max\{4t\lfloor N^{1/t}\rfloor+3t,12N+8\}$ and depth $12L+14$ has an approximation rate $18\sqrt{t}\omega(N^{-2/t}L^{-2/t})$ on its trifling region $\Omega_{j}:=\Omega([0,1]^t,K,\delta)$. Paralleling these $d$ equal-depth networks result in a wider network with width $d\times\max\{4t\lfloor N^{1/t}\rfloor+3t,12N+8\}$, depth $12L+14$ and trifling region $\Omega([0,1]^d,K,\delta)$ which covers the projection of all $\Omega_{j}$ onto $[0,1]^d$, i.e. $\cup_{j=1,\ldots,d}{\rm Proj}_{[0,1]^d}(\Omega_{j})\subset\Omega([0,1]^d,K,\delta)$.
\end{proof}

\subsection{Proof of Lemma \ref{lemma7}}

\begin{proof}
	Recall that $h_{ij}:\mathbb{R}^{t_{i}}\to\mathbb{R}$, $i=0,\ldots,q$ and $j=1,\ldots,d_{i+1}$ are H\"older continuous functions with order $\alpha_i\in[0,1]$ and constant $\lambda_i\ge0$ and $h_i=(h_{ij})_j^\top: \mathbb{R}^{d_{i}}\to\mathbb{R}^{d_{i+1}}$ are vectors of functions with domain $D_i$.
	Let $H_i=h_i\circ\ldots\circ h_0$ and $\tilde{H}_i=\tilde{h}_i\circ\ldots\circ \tilde{h}_0$ for $i=0,\ldots,q$. Let $S_{ij}\subset\{1,\ldots,d_{i+1}\}$ be the support of the $t_i$-variate function $h_{ij}$ and denote $x_{S_{ij}}$ by the $d_{i+1}$-dimensional vector $x$ restricted to the $t_i$-dimensional subspace according to the index $S_{ij}$. then
	\begin{align*}
		&\Vert h_q\circ\ldots h_0-\tilde{h}_q\circ\ldots \tilde{h}_0\Vert_{L^\infty(D_0)}\\
		=& \Vert h_q\circ H_{q-1}-h_q\circ \tilde{H}_{q-1}+h_q\circ\tilde{H}_{q-1}-\tilde{h}_q\circ\tilde{H}_{q-1}\Vert_{L^\infty(D_0)}\\
		\leq& \Vert h_q\circ H_{q-1}-h_q\circ \tilde{H}_{q-1}\Vert_{L^\infty(D_0)}+\Vert h_q\circ\tilde{H}_{q-1}-\tilde{h}_q\circ\tilde{H}_{q-1}\Vert_{L^\infty(D_0)}\\
		\leq& \max_{j=1,\ldots,d_{q+1}} \sup_{x\in D_0}\vert h_{qj}\circ H_{q-1}(x)-h_{qj}\circ \tilde{H}_{q-1}(x)\vert+\Vert h_q-\tilde{h}_q\Vert_{L^\infty(D_q)}\\
		\leq& \max_{j=1,\ldots,d_{q+1}} \omega_{h_{qj}}(\sup_{x\in D_0}\Vert H_{q-1}(x)_{S_{ij}}-{\tilde{H}_{q-1}(x)}_{{S_{ij}}}\Vert_2)+\Vert h_q-\tilde{h}_q\Vert_{L^\infty(D_q)}\\
		\leq& \max_{j=1,\ldots,d_{q+1}} \omega_{h_{qj}}(\sqrt{t_q}\Vert H_{q-1}-\tilde{H}_{q-1}\Vert_{L_\infty(D_0)})+\Vert h_q-\tilde{h}_q\Vert_{L^\infty(D_q)}\\
		\leq& \lambda_qt_q^{\alpha_q/2}\Vert H_{q-1}-\tilde{H}_{q-1}\Vert^{\alpha_{q}}_{L^\infty(D_0)}+\Vert h_q-\tilde{h}_q\Vert_{L^\infty(D_q)}\\
		\leq& \lambda_qt_q^{\alpha_q/2}\big(\lambda_{q-1}t_{q-1}^{\alpha_{q-1}/2}\Vert H_{q-2}-\tilde{H}_{q-2}\Vert_{L^\infty(D_0)}^{\alpha_{q-1}}+\Vert h_{q-1}-\tilde{h}_{q-1}\Vert_{L^\infty(D_{q-1})}\big)^{\alpha_{q}}\\
		&\qquad\qquad\qquad\qquad\qquad\qquad\qquad\qquad\qquad\qquad+\Vert h_q-\tilde{h}_q\Vert_{L^\infty(D_q)}\\
		\leq& \lambda_q\lambda_{q-1}^{\alpha_{q}}t_q^{\alpha_q/2}t_{q-1}^{\alpha_{q}\alpha_{q-1}/2}\Vert H_{q-2}-\tilde{H}_{q-2}\Vert_{L^\infty(D_0)}^{\alpha_{q}\alpha_{q-1}}\\
		&\qquad\qquad\qquad\qquad+\lambda_qt_{q}^{\alpha_{q}/2}\Vert h_{q-1}-\tilde{h}_{q-1}\Vert_{L^\infty(D_{q-1})}^{\alpha_{q}}+\Vert h_q-\tilde{h}_q\Vert_{L^\infty(D_q)}\\
		\leq &\sum_{i=0}^q \Pi_{j=i+1}^{q}\lambda_{j}^{\Pi_{k=j+1}^{q} \alpha_{k}} \Pi_{j=i+1}^{q}\sqrt{t_{j}}^{\Pi_{k=j}^{q} \alpha_{k}} \Vert h_i-\tilde{h}_i\Vert_{L^\infty(D_i)}^{\Pi_{j=i+1}^{q}\alpha_{j}}.
	\end{align*}
The third inequality follows from $\Vert x\Vert_2\leq\sqrt{d}\Vert x\Vert_\infty$ for a vector $x\in\mathbb{R}^d$. The fourth inequality follows from the definition of H\"older continuity. The second last inequality follows from $(a+b)^\alpha\leq a^\alpha+b^\alpha$ for all $a,b\ge0$ and $\alpha\in[0,1]$.
\end{proof}

\subsection{Proof of Lemma \ref{lemma8}}
\begin{proof}
	We start our proof from the most simple case where $h:\mathbb{R}^d\to\mathbb{R}$ be a linear combination operator, i.e., $h(x)=Tx+u$ with $T=(t_1,\ldots,t_d)\in\mathbb{R}^{1\times d}$ being a row vector and $u\in\mathbb{R}$ being a scalar. Then we can construct a three-layer ReLU neural network $\tilde{h}(x)=W_2\sigma(W_1x+b_1)+b_2$ with width $(d,2d,1)$ where $\sigma(\cdot)$ is the ReLU activation function, $b_1=\textbf{0}$, $b_2=u$,
	\begin{equation*}
		W_1=\left[
		\begin{array}{ccccccc}
			1& 0 & 0 &\cdots & \cdots &0 & 0\\
			-1& 0& 0 &\cdots & \cdots &0 & 0\\
			0 &  1& 0& 0 & 0 & \cdots  &0\\
			0 &  -1 & 0& 0 & 0 & \cdots  &0\\
			\vdots & \ddots & \ddots &   \ddots& \ddots & \ddots  &\vdots\\
			0 & \cdots & \cdots  &  \cdots & \cdots  & 0 & 1\\
			0 & \cdots & \cdots  &  \cdots & \cdots  & 0 & -1\\
		\end{array}
		\right],
	\end{equation*}
	and $W_2=(t_1,-t_1,t_2,-t_2,\ldots,t_{d-1},-t_{d-1},t_d,-t_d)_{1\times 2d}$ is a $2d$-dimensional row vector. And it is easy to verify that $\tilde{h}(x)=h(x)$, for any $x\in\mathbb{R}^d$. More generally, when $T=(t_{ij})\in\mathbb{R}^{m\times d}$ and $u\in\mathbb{R}^{m}$, we can construct the three-layer network with width $(d,2d,m)$ in a similar manner where $W_1$, $b_1$ and $b_2$ are kept the same as above but $W_2\in\mathbb{R}^{m\times 2d}$ is constructed analogically by stacking $m$ many $2d$-dimensional vectors together, i.e.,
	\begin{equation*}
		W_2=\left[
		\begin{array}{ccccccc}
			t_{11}& -t_{11} & t_{12} &-t_{12} & \cdots &t_{1d} & -t_{1d}\\
			\vdots & \ddots & \ddots &   \ddots& \ddots & \ddots  &\vdots\\
			t_{m1}& -t_{m1} & t_{m2} &-t_{m2} & \cdots &t_{md} & -t_{md}\\
		\end{array}
		\right].
	\end{equation*}
	In such a way, the constructed $\tilde{h}$ satisfies $\tilde{h}(x)=h(x)$ for any $x\in\mathbb{R}^d$.
	
\end{proof}

\subsection{Proof of Theorem \ref{thm1}}
\begin{proof}

In Lemma \ref{lemma5} and Lemma \ref{lemma6}, the domain of the approximated functions are required to be $[0,1]^{d}$. In light of this, the Lemmas can not be directly applied to each $h_i$ of the composition since in general neither the  domain of $h_i$  is $[0,1]^{d_{i}}$ nor the range of $h_i$ is $[0,1]^{d_{i+1}}$.
 Thus the domain of the constructed ReLU networks have to be aligned with the approximated functions $h_{i}$. Considering this, we can add an additional invertible linear layer ${A}_i(\cdot):D_{i}\to[0,1]^{d_i}$ at the beginning of each of the subnetworks $\tilde{h}_{i}$ in Lemma \ref{lemma6} for $0=1,\ldots,q$ to accommodate to general $h_i$. In the following, we introduce the accommodation in details.

Note that all $h_i$, $i=0,\ldots,q$ are continuous functions on bounded domain $D_i$, where $D_0=[a,b]^d$ and $h_{i-1}\circ\ldots\circ h_0([a,b]^d)\subseteq D_{i}$ for $i=1,\ldots,q$. Without loss of generality, we can let $a_{i}:=\min_{j=1,\ldots,d_{i-1}}  \inf_{x\in[a,b]^d}h_{(i-1)j}\circ\ldots\circ h_0(x)$ and $b_{i}:=\max_{j=1,\ldots,d_{i-1}}  \sup_{x\in[a,b]^d}h_{(i-1)j}\circ\ldots\circ h_0(x)$ for $i=1,\ldots,q$. Then we can view $h_i$ as functions with domain $[a_i,b_i]^{d_i}$. Further, for each $i\in\{0,\ldots,q\}$, these exists an invertible linear transformation $A_i(x)=\sigma(W_ix+b_i)$ where $W_i\in\mathbb{R}^{d_i\times d_i}$ is a diagonal matrix with equivalent entries $1/(b_i-a_i)$, $b_i\in\mathbb{R}^{d_i}$ is a vector with equivalent components $-a_i/(b_i-a_i)$ and $\sigma(\cdot)$ is the ReLU activation function such that $A_i$ is an invertible transformation from $[a_i,b_i]^{d_i}$ to $[0,1]^{d_i}$. Now we can apply Lemma \ref{lemma6} to build up networks approximate $h_i$ on domains $[a_i,b_i]^{d_i}$.

For any $L_i\in\mathbb{N}^+$ and $N_i\in\mathbb{N}^+$, there exists functions $\tilde{h}_i$ for $i\in J^c$ implemented by ReLU FNNs with width $d_{i}\max\{4t_i\lfloor N_i^{1/t_i}\rfloor+3t_i,12N_i+8\}$ and depth $12L_i+15$ such that $\Vert\tilde{h}_i\Vert_{L_i^\infty(\mathbb{R}^{d_i})}\leq \max_{j=1,\ldots,d_i} \vert h_{ij}(\textbf{0})\vert+\omega(\sqrt{t_i})$ and
$$\vert \tilde{h}_i(x)-h_i(x)\vert\leq18\sqrt{t_i}\lambda_i(N_iL_i)^{-2\alpha_i/t_i}, \quad{\rm for\ any\ } x\in D_i\backslash\ A^{-1}_i(\Omega([0,1]^{d_i},K,\delta)),$$
where $A^{-1}_i:[a_i,b_i]^{d_i}\to[0,1]^{d_i}$ is the inverse of above defined linear transformation $A_i$ (the first layer of $\tilde{h}_i$), $K_i=\lfloor N_i^{1/d_i}\rfloor^2\lfloor L_i^{1/d_i}\rfloor^2$ and $\delta_i$ is an arbitrary number in $(0,1/(3K_i)]$. And the trifling region $\Omega([0,1]^d,K,\delta)$ of $[0,1]^d$ is defined as $$\Omega([0,1]^d,K,\delta)=\cup_{i=1}^d\{x=[x_1,x_2,...,x_d]^T:x_i\in\cup_{k=1}^{K-1}(k/K-\delta,k/K)\},$$
and
$$A^{-1}_i(\Omega([0,1]^{d_i},K,\delta))=\{x\in\mathbb{R}^{d_i}:A(x)\in\Omega([0,1]^{d_i},K,\delta\}.$$

By Lemma \ref{lemma8}, for $j\in J$, there exists functions $\tilde{h}_j$ implemented by 3-layer ReLU FNNs with width vector $(d_j,2d_j,d_{j+1})$ such that
$$\vert \tilde{h}_j(x)-h_j(x)\vert=0 \quad{\rm for\ any\ } x\in\mathbb{R}^{d_j}.$$
To approximate the composited function $H_q=h_q\circ\ldots\circ h_0: [a,b]^d\to\mathbb{R}$, we let  $\tilde{H}_q=\tilde{h}_q\circ\ldots\circ \tilde{h}_0$ be the composition of above defined $\tilde{h}_i$, which is a  function implemented by ReLU FNN with width $\max\{\max_{i\in J^c}d_{i}\max\{4t_i\lfloor N_i^{1/t_i}\rfloor+3t_i,12N_i+8\},\max_{j\in J}2d_j\}$ and depth $\sum_{i\in J^c}(12L_i+15)+2\vert J\vert$.  Then by applying Lemma \ref{lemma7}, we have
	\begin{align*}
		&\vert \tilde{H}_q(x)-H_q(x)\vert\\
		\leq&\sum_{i\in J^c} \Pi_{j=i+1}^{q}\lambda_{j}^{\Pi_{k=j+1}^{q} \alpha_{k}} \Pi_{j=i+1}^{q}\sqrt{t_{j}}^{\Pi_{k=j}^{q} \alpha_{k}} \big( 18\sqrt{t_i}\lambda_i\big)^{\Pi_{j=i+1}^{q}\alpha_{j}}(N_iL_i)^{-2(\Pi_{j=i}^{q}\alpha_{j})/t_i}\\
		\leq&\sum_{i\in J^c} 18^{\Pi_{j=i+1}^{q}\alpha_{j}}\Pi_{j=i}^{q}\lambda_{j}^{\Pi_{k=j+1}^{q} \alpha_{k}} \frac{\Pi_{j=i}^{q}\sqrt{t_{j}}^{\Pi_{k=j}^{q} \alpha_{k}}}{\sqrt{t_i}^{\alpha_i}} (N_iL_i)^{-2(\Pi_{j=i}^{q}\alpha_{j})/t_i}\\
		=&\sum_{i\in J^c} C_i^*\lambda_i^* t_i^*(N_iL_i)^{-2\alpha_i^*/t_i}, \qquad {\rm for\ any\ }x\in [a,b]^d\backslash \Omega_0,
	\end{align*}
	where $\lambda_j=\alpha_j=1$ for $j\in J$, $C_i^*=18^{\Pi_{j=i+1}^{q}\alpha_{j}}$, $\lambda_i^*=\Pi_{j=i}^{q}\lambda_{j}^{\Pi_{k=j+1}^{q} \alpha_{k}}$, $\alpha_i^*=\Pi_{j=i}^q \alpha_j$, $t_i^*={(\Pi_{j=i}^{q}\sqrt{t_{j}}^{\Pi_{k=j}^{q} \alpha_{k}})}/{\sqrt{t_i}^{\alpha_i}}$ and $\Omega_0$ is a subset of $[a,b]^d$ which satisfies
	$$ \Omega([0,1]^{d_i},K_i,\delta_i)\subseteq A_i\circ\tilde{h}_{i-1}\circ\ldots\circ \tilde{h}_0(\Omega_0), \qquad {\rm for\ } i=0,\ldots,q,$$
	where $A_j$ is defined as identity map for $j\in J$.
	Note that since $\alpha_i\in[0,1]$, further we have $C_i^*\leq 18$ and $t_i^*\leq \Pi_{j=i}^q\sqrt{t_j}\leq \Pi_{j=0}^q\sqrt{t_j}$.
\end{proof}

\subsection{Proof of Theorem \ref{thm2}}
\begin{proof}
	By Theorem \ref{thm1}, given any $N_i,L_i\in\mathbb{N}^+, i\in J^c$, for the function class of ReLU multi-layer perceptrons $\mathcal{F}_{\phi}=\mathcal{F}_{\mathcal{D},\mathcal{W},\mathcal{U},\mathcal{S},\mathcal{B}}$
	with width $\mathcal{W}=\max\{\max_{i\in J^c}d_{i}\max\{4t_i\lfloor N_i^{1/t_i}\rfloor+3t_i,12N_i+8\},\max_{j\in J}2d_j\}$ and depth $\mathcal{D}=\sum_{i\in J^c}(12L_i+15)+2\vert J\vert$, there exists a $f^*_\phi$ such that
	\begin{align*}
		\vert f^*_\phi(x)-f_0(x)\vert\leq\sum_{i\in J^c} C_i^*\lambda_i^* t_i^*(N_iL_i)^{-2\alpha_i^*/t_i}, \qquad {\rm for\ any\ }x\in [a,b]^d\backslash \Omega_0,
	\end{align*}
	where $C_i^*=18^{\Pi_{j=i+1}^{q}\alpha_{j}}$, $\lambda_i^*=\Pi_{j=i}^{q}\lambda_{j}^{\Pi_{k=j+1}^{q} \alpha_{k}}$, $\alpha_i^*=\Pi_{j=i}^q \alpha_j$, $t_i^*={(\Pi_{j=i}^{q}\sqrt{t_{j}}^{\Pi_{k=j}^{q} \alpha_{k}})}/{\sqrt{t_i}^{\alpha_i}}$ and $\Omega_0$ is a subset of $[a,b]^d$ which satisfies
	$$ \Omega([0,1]^{d_i},K_i,\delta_i)\subseteq A_i\circ\tilde{h}_{i-1}\circ\ldots\circ \tilde{h}_0(\Omega_0), \qquad {\rm for\ } i=0,\ldots,q,$$
	where $A_i$ are defined as in Theorem \ref{thm1}.
	Note that the Lebesgue measure of each $\Omega([0,1]^{d_i},K_i,\delta_i)$ is no more than $\delta_i(K_i-1)d$ which can be arbitrarily small since $\delta_i\in(0,1/(3K_i))$ can be arbitrarily small. Thus the preimage or inverse image of $\Omega([0,1]^{d_i},K_i,\delta_i)$ under $A_i\circ\tilde{h}_{i-1}\circ\ldots\circ \tilde{h}_0$ can has arbitrarily small Lebesgue measure since all $A_i, \tilde{h}_i$ are continuous mappings. As a consequence, the Lebesgue measure of $\Omega_0$ can be arbitrarily small by choosing arbitrarily small $\delta_i$. Besides, $\nu$ (the probability measure of $X$) is absolutely continuous with respect to Lebesgue measure, then we have
	\begin{align*}
		\mathbb{E}_{X}\vert f^*_\phi(X)-f_0(X)\vert=\Vert f^*_\phi-f_0\Vert_{L^2(\nu)}\leq \sum_{i\in J^c} C_i^*\lambda_i^* t_i^*(N_iL_i)^{-2\alpha_i^*/t_i}.
	\end{align*}
	Combining Lemma \ref{lemma2}-\ref{lemma3}, we have for $2n \ge \text{Pdim}(\mathcal{F}_\phi)$,
	the prediction error of the \textit{DQR} estimator $\hat{f}_\phi$ satisfies
	\begin{equation*}
		\mathbb{E}\big\{\mathcal{R}^\tau(\hat{f}_\phi)-\mathcal{R}^\tau(f_0)\big\}\leq C\frac{\lambda_\tau\mathcal{B}\mathcal{S}\mathcal{D}\log(\mathcal{S})\log(n)}{n^{1-1/p}}+ 2\lambda_\tau\sum_{i\in J^c} C_i^*\lambda_i^* t_i^*(N_iL_i)^{-2\alpha_i^*/t_i},
	\end{equation*}
	where $\lambda_\tau=\max\{\tau,1-\tau\}$ and $C>0$ is a constant does not depend on $n,d,\tau,\mathcal{B},\mathcal{S},\mathcal{D},C_i^*,\lambda_i^*,\alpha_i^*,N_i$ or $L_i$, and $C_i^*=18^{\Pi_{j=i+1}^{q}\alpha_{j}}$, $\lambda_i^*=\Pi_{j=i}^{q}\lambda_{j}^{\Pi_{k=j+1}^{q} \alpha_{k}}$, $\alpha_i^*=\Pi_{j=i}^q \alpha_j$ and $t_i^*={(\Pi_{j=i}^{q}\sqrt{t_{j}}^{\Pi_{k=j}^{q} \alpha_{k}})}/{\sqrt{t_i}^{\alpha_i}}$.
	If Assumption \ref{quadratic} additionally holds, then combining Lemma \ref{lemma2},\ref{lemma4}, the approximation result can be directly applied,
	\begin{equation*}
		\mathbb{E}\big\{\mathcal{R}^\tau(\hat{f}_\phi)-\mathcal{R}^\tau(f_0)\big\}\leq C\frac{\lambda_\tau\mathcal{B}\mathcal{S}\mathcal{D}\log(\mathcal{S})\log(n)}{n^{1-1/p}}+ 2c_\tau\big[\sum_{i\in J^c} C_i^*\lambda_i^* t_i^*(N_iL_i)^{-2\alpha_i^*/t_i}\big]^2,
	\end{equation*}
	where $c_\tau>0$ is a constant defined in Lemma \ref{lemma4}.
\end{proof}

\subsection{Proof of Lemma \ref{lemma9}}
\begin{proof}
By equation (B.3) in \cite{belloni2011l1}, for any scalar $w,v\in\mathbb{R}$ we have
\begin{align*}
	\rho_\tau(w-v)-\rho_\tau(w)=-v\{\tau-I(w\leq0)\}+\int_0^v\{I(w\leq z)-I(w\leq0)\}dz.
\end{align*}
Given any $f$ and $X=x$, let $w=Y-f_0(X)$, $v=f(X)-f_0(X)$ with $\vert f(x)-f_0(x)\vert\leq \gamma$. Then given $X=x$, taking conditional expectation on above equation with respect to $Y\mid X=x$, we have
\begin{align*}
	&\mathbb{E}\{\rho_\tau(Y-f(X))-\rho_\tau(Y-f_0(X))\mid X=x\}\\
	=&\mathbb{E}\big[-\{f(X)-f_0(X)\}\{\tau-I(Y-f(X)\leq0)\}\mid X=x\big]\\
	&+\mathbb{E}\big[\int_0^{f(X)-f_0(X)}\{I(Y-f_0(X)\leq z)-I(Y-f_0(X)\leq0)\}dz\mid X=x\big]\\
	=&0+\mathbb{E}\big[\int_0^{f(X)-f_0(X)}\{I(Y-f_0(X)\leq z)-I(Y-f_0(X)\leq0)\}dz\mid X=x\big]\\
	=&\int_0^{f(x)-f_0(x)} \{P_{Y|X}(f_0(x)+z)-P_{Y|X}(f_0(x))\}dz\\
	\geq&\int_0^{f(x)-f_0(x)} \kappa \vert z\vert dz\\
	=&\frac{\kappa}{2}\vert f(x)-f_0(x)\vert^2.
\end{align*}
Suppose $f(x)-f_0(x)> \gamma$, then similarly we have
\begin{align*}
	&\mathbb{E}\{\rho_\tau(Y-f(X))-\rho_\tau(Y-f_0(X))\mid X=x\}\\
	=&\int_0^{f(x)-f_0(x)} \{P_{Y|X}(f_0(x)+z)-P_{Y|X}(f_0(x))\}dz\\
	\geq&\int_{\gamma/2}^{ f(x)-f_0(x)} \{P_{Y|X}(f_0(x)+\gamma/2)-P_{Y|X}(f_0(x))\}dz\\
	\geq&(f(x)-f_0(x)-\gamma/2)(\kappa \gamma/2)\\
	\geq &\frac{\kappa\gamma}{4}\vert f(x)-f_0(x)\vert.
\end{align*}

The case $f(x)-f_0(x)\leq -\gamma$ can be handled similarly as in \cite{padilla2021adaptive}. The conclusion follows combining the three different cases and taking expectation with respect to $X$ of above obtained inequality.
\end{proof}

\subsection{Proof of Theorem \ref{thm3}}
\begin{proof}
Theorem \ref{thm3} follows directly from Theorem \ref{thm2} and Lemma \ref{lemma9}.
\end{proof}


\section{Additional simulation results}
\label{AppFigures}
In this section, we provide additional simulation results, including the estimated quantile curves at $\tau=$ 0.05, 0.25,0.5,0.75, and 0.95,
the corresponding excess risks, the $L_1$ and the $L_2$ test errors.
To make this section self-contained as much as possible, we also include the detailed description of the simulation studies in the main text of the paper.

We consider the following quantile regression methods:
\begin{itemize}
	\item The traditional linear quantile regression as described in \cite{koenker1978}, denoted by \textit{linear QR}. Without regularization, the empirical risk is minimized over the parameter space (intercept included) $\mathbb{R}^{d+1}$ to give an linear estimator. These estimation are implemented on Python via package \textit{statsmodels}.
	\item   Kernel-based nonparametric quantile regression  as described in \cite{sangnier2016joint}, denoted by \textit{kernel QR}.
	This is a joint quantile regression method based on vector-valued reproducing kernel Hilbert space (RKHS), which enjoys few quantile crossing and enhanced performances compared to independent estimations and hard non-crossing constraints. In our implementation, the radial basis function (RBF) kernel is chosen and a coordinate descent primal-dual algorithm \citep{fercoq2019coordinate} is used via Python package \textit{qreg}.
	\item Deep quantile regression as described in Section \ref{sec2}, denoted by \textit{DQR}. We implement it in Python via \textit{Pytorch} and use \textit{Adam} \citep{kingma2014adam} as the optimization algorithm with default learning rate 0.01 and default $\beta=(0.9,0.99)$ (coefficients used for computing running averages of gradients and their squares).
	\item Deep least squares regression, denoted by \textit{DLS}. We minimize the mean square error on the training data to get the nonparametric least square estimator using deep neural networks. Similarly we implement it on Python via \textit{Pytorch} and use \textit{Adam} as the optimization algorithm with default settings. The comparison with DLS mainly focuses on the $0.5$-th quantile curve since the conditional mean and the conditional median coincident with each other when error is symmetric.
\end{itemize}

\subsection{Estimations and Evaluations}
We consider estimating the quantile curves at 5 different levels for each simulated model, i.e., we estimate quantile curves for $\tau\in\{0.05,0.25,0.5,0.75,0.95\}$. For each model $f_0$ and each error $\eta$, according to model (\ref{model}) we generate the training data $(X_i^{train},Y_i^{train})_{i=1}^n$ with sample size $n$ to train the empirical risk minimizer at $\tau\in\{0.05,0.25,0.5,0.75,0.95\}$ by different methods, i.e.
\begin{align*}
	\hat{f}^\tau_n\in\arg\min_{f\in\mathcal{F}} \frac{1}{n}\sum_{i=1}^n\rho_\tau(Y_i^{train}-f(X_i^{train})),
\end{align*}
where $\mathcal{F}$ is the class of linear functions, RKHS
or the class of
ReLU neural network functions. For each $f_0$ and each error $\eta$, we also generate the testing data $(X_t^{test},Y_t^{test})_{t=1}^T$ with sample size $T$ from the same distribution of the training data. Then for each obtained $\hat{f}^\tau_n$, we calculate its testing risk on $(X_t^{test},Y_t^{test})_{t=1}^T$, i.e.,
 \begin{align*}
 	\mathcal{R}^\tau(\hat{f}^\tau_n)=\frac{1}{T}\sum_{t=1}^T \rho_\tau(Y_t^{test}-\hat{f}^\tau_n(X_t^{test})).
 \end{align*}
Moreover, for each obtained $\hat{f}^\tau_n$, we calculate the $L_1$ distance between $\hat{f}^\tau_n$ and the corresponding risk minimizer $f_0^\tau$, i.e.
 \begin{align*}
\Vert\hat{f}^\tau_n-f_0^\tau\Vert_{L^1(\nu)}=\frac{1}{T}\sum_{t=1}^T \vert\hat{f}^\tau_n(X_t^{test})-f_0^\tau(X_t^{test})\vert,
\end{align*}
and we also calculate the $L_2$ distance between $\hat{f}^\tau_n$ and the corresponding risk minimizer $f_0^\tau$, i.e.
 \begin{align*}
	\Vert\hat{f}^\tau_n-f_0^\tau\Vert^2_{L^2(\nu)}=\frac{1}{T}\sum_{t=1}^T \vert\hat{f}^\tau_n(X_t^{test})-f_0^\tau(X_t^{test})\vert^2.
\end{align*}
All the $L_2$ test error results are provided in the appendix.
The specific forms of $f_0^\tau$ are given in the part on the data generation models below.

In the simulation studies, we take $T=100,000$ as the sample size of testing data for each data generation model. We report the mean and standard deviation of statistics including excess risk $\mathcal{R}^\tau(\hat{f}^\tau_n)-\mathcal{R}^\tau({f}^\tau_0)$, $L_1$ distance and $L^2_2$ distance  over $R = 10$ replications under different scenarios.  For \textit{DLS}, the testing risk and the excess risk are calculated in terms of mean squares loss function other than the check loss $\rho_\tau$.

\subsection{Data generation: univariate models}
We generate data according to model (\ref{model}), i.e., $Y=f_0(X)+\eta$. We consider three basic univariate models, including ``Linear'', ``Wave'' and ``Triangle'', which corresponds to different specifications of $f_0$. The formulae are given below.
\begin{enumerate}
\setlength\itemsep{-0.05 cm}
	\item Linear:
	\begin{align*} f_0(x)=2x.
\end{align*}
	\item Wave:
	$$f_0(x)=2x\sin(4\pi x).$$
	\item Triangle:
	$$f_0(x)=4(1-\vert x-0.5\vert).$$
\end{enumerate}
We use the linear model as a baseline model in our simulations and expect all the methods
perform well under the linear model. The ``Wave'' is a nonlinear smooth model and the ``Triangle'' is a
nonlinear continuous but non-differentiable model. These models
are chosen so that we can evaluate the
performance of DQR,  $kernel QR$ and $linear QR$ under different types of  models.

\begin{figure}[H]
\centering
\includegraphics[width=\textwidth]{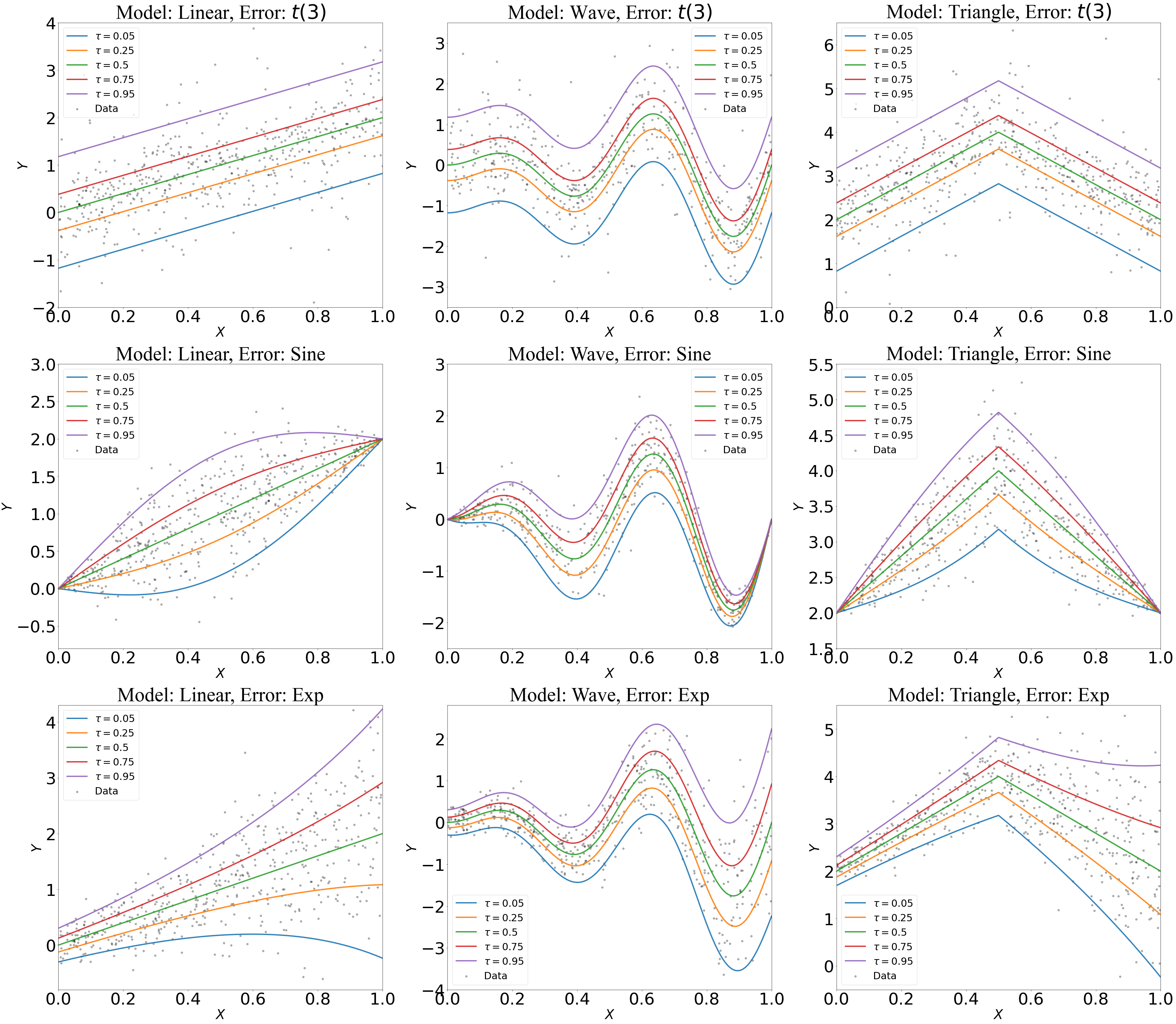}
\caption{The target quantiles curves at $\tau=0.05,0.25,0.5,0.75,0.95$
under different models and error distributions. From the left to the right, each column corresponds a data generation model, i.e., ``Linear'', ``Wave'' and ``Triangle''. From the top to the bottom, each row corresponds a error distribution, i.e. $t(3)$, ``\textit{Sine}'' and ``\textit{Exp}''.}
\label{fig:targetA}
\end{figure}

For these models, we generate $X$ uniformly from the unit interval $[0,1]$. We generate the error $\eta$ from the following distributions.
\begin{enumerate}
\setlength\itemsep{-0.05 cm}
	\item $\eta$ follows a scaled Student's t distribution with degrees of freedom 3, i.e., $\eta\sim 0.5\times t(3)$, denoted by $t(3)$;
	\item  Conditioning on $X=x$, the error $\eta$ follows a normal distribution whose variance depends on the covariate $X$, i.e., $\eta\mid X=x\sim 0.5\times\mathcal{N}(0,[\sin(\pi x)]^2)$, denoted by \textit{Sine};
	\item Conditioning on $X=x$, the error $\eta$ follows a normal distribution whose variance depends on the covariate $X$, i.e., $\eta\mid X=x\sim0.5\times\mathcal{N}(0,\exp(4x-2))$, denoted by \textit{Exp}.
\end{enumerate}
Note that except for $t(3)$, other two types of errors depend on the predictor $X$.
The $\tau$-th conditional quantile $f_0^\tau(x)$ of the response $Y$ given $X=x$ can be calculated by
$$f_0^\tau(x)=f_0(x)+F^{-1}_{\eta\mid X=x}(\tau),$$
where $F^{-1}_{\eta\mid X=x}(\cdot)$ is the inverse of the conditional cumulated distribution function of $\eta$ given $X=x$. For $t(3)$ error, $\eta$ is independent with $X$, then $F^{-1}_{\eta\mid X=x}(\cdot)$ is simply the inverse of distributional function of the $2 t(3)$. For the \textit{Sine} error, $F^{-1}_{\eta\mid X=x}(\tau)=0.5\times\sin(\pi x)\times\Phi^{-1}(\tau)$ where $\Phi^{-1}(\cdot)$ is the inverse of the CDF of a standard normal random variable. Similarly, for the \textit{Exp} error, $F^{-1}_{\eta\mid X=x}(\tau)=0.5\times\exp(2x-1)\times\Phi^{-1}(\tau)$.
Figure \ref{fig:targetA} shows all these univariate data generation models and their corresponding conditional quantiles at $\tau=$ 0.05,0.25,0.5,0.75,0.95.

\begin{figure}[H]
	\centering
	\begin{subfigure}{\textwidth}
		\includegraphics[width=1\textwidth]{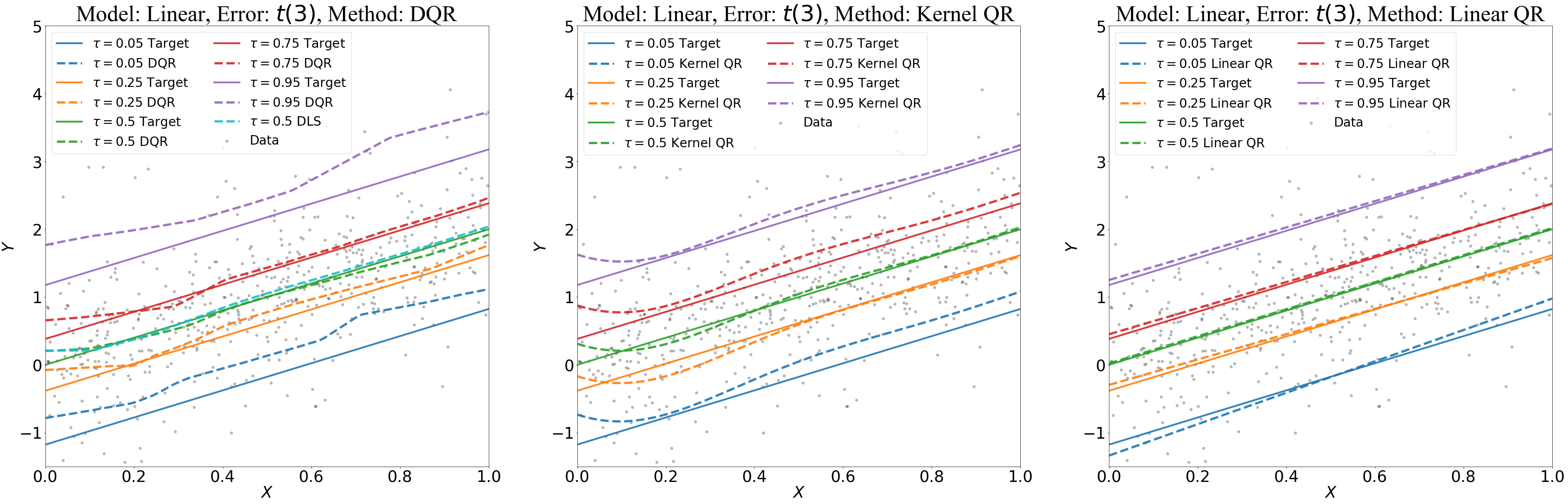}
	\end{subfigure}
	\begin{subfigure}{\textwidth}
		\includegraphics[width=1\textwidth]{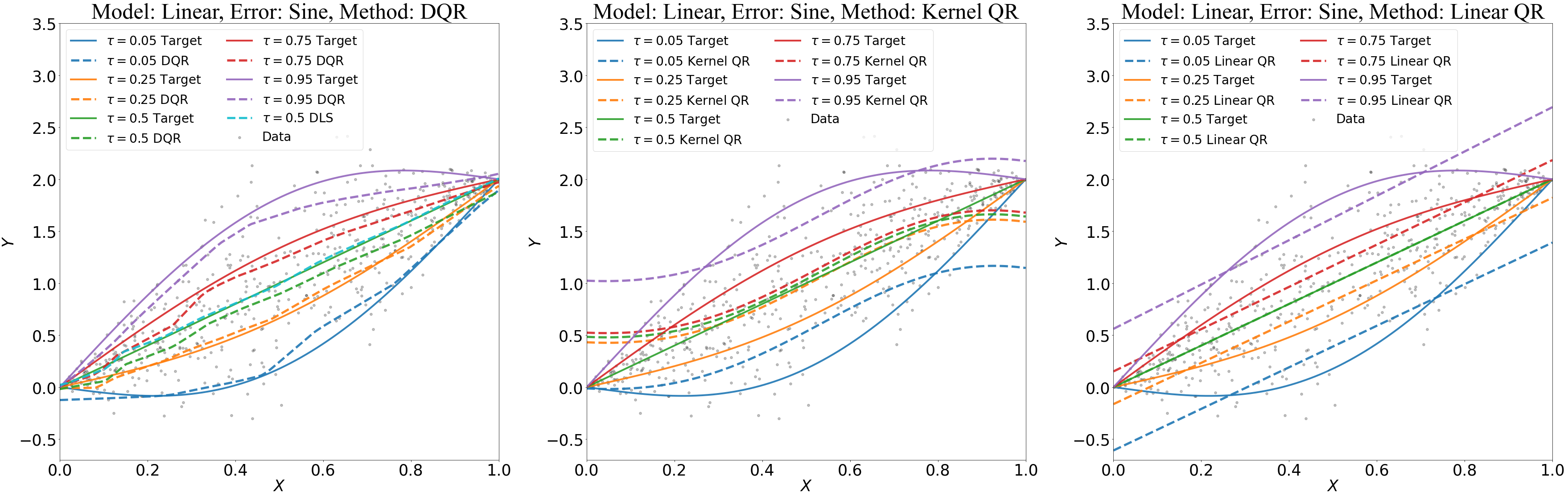}
	\end{subfigure}
	\begin{subfigure}{\textwidth}
		\includegraphics[width=\textwidth]{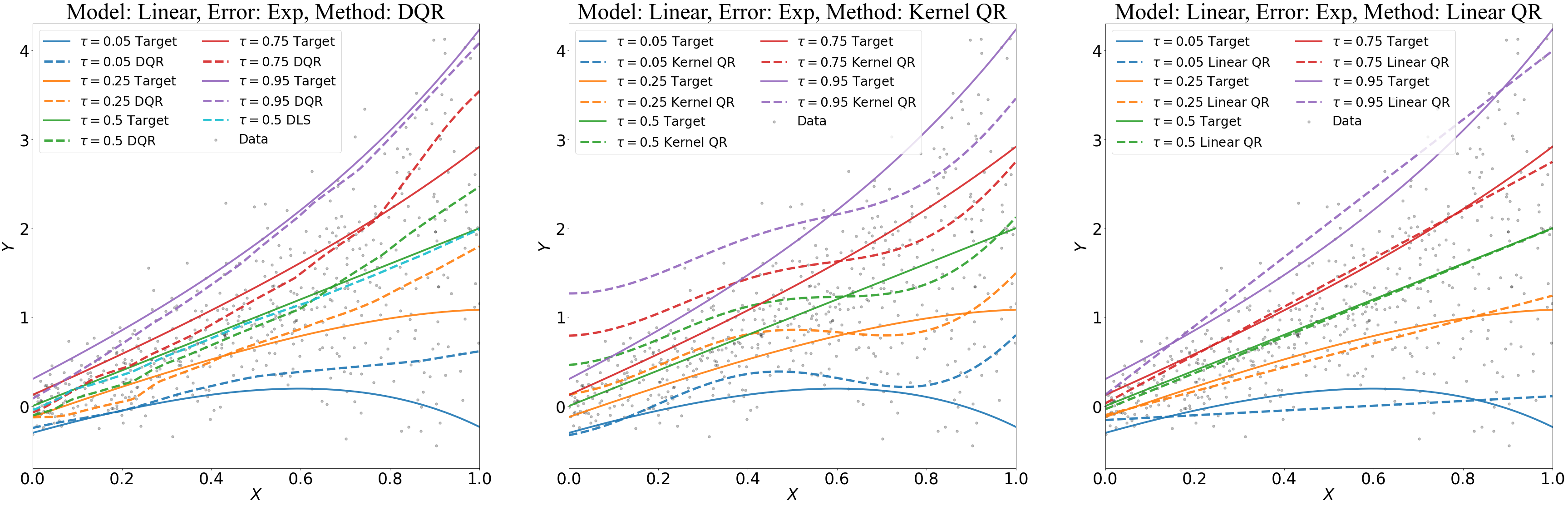}
	\end{subfigure}
	\caption{The fitted quantile curves by different methods under univariate model ``Linear'' with different errors. The training data is depicted as grey dots.The target quantile functions at $\tau=0.05,0.25,0.5,0.75,0.95$ are depicted as solid curves in different colors, and colored dashed curves represent the corresponding estimates. From the top to the bottom, each row corresponds a certain type of error: $t(3)$, ``\textit{Sine}'' and ``\textit{Exp}''. From the left to right, each column corresponds a certain estimation method: \textit{DQR}, \textit{kernel QR} and \textit{linear QR}. Fitted \textit{DLS} curves are contained in the \textit{DQR} plots.}
	\label{fitted:linear1A}
\end{figure}

\begin{figure}[H]
	\centering
	\begin{subfigure}{\textwidth}
		\includegraphics[width=1\textwidth]{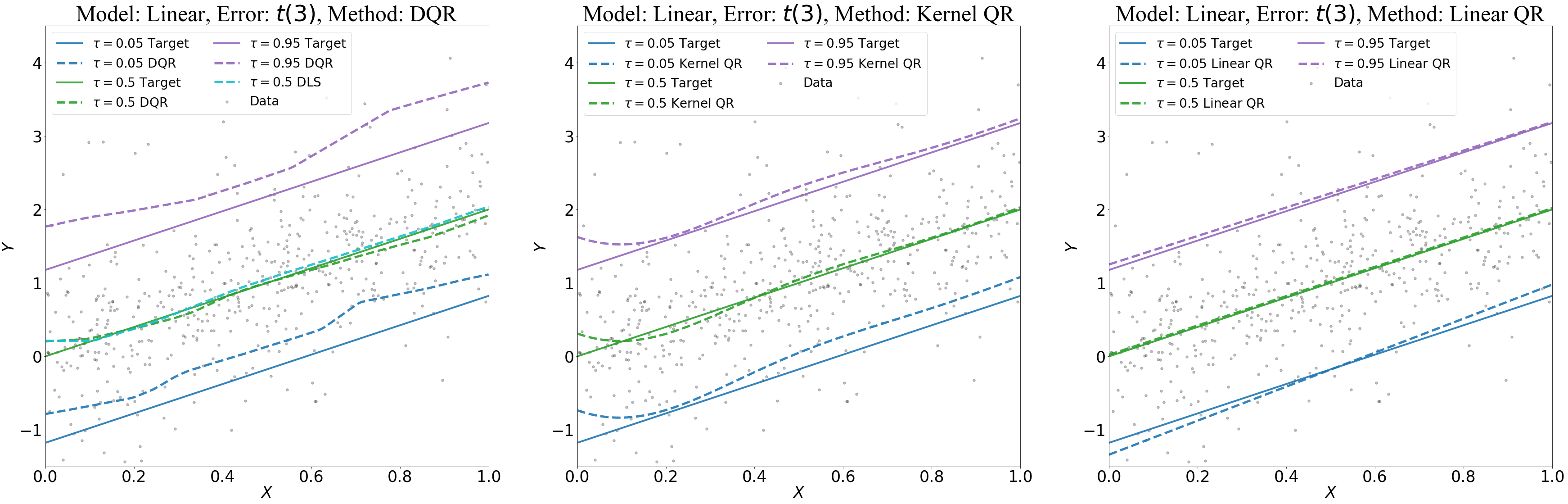}
	\end{subfigure}
	\begin{subfigure}{\textwidth}
		\includegraphics[width=1\textwidth]{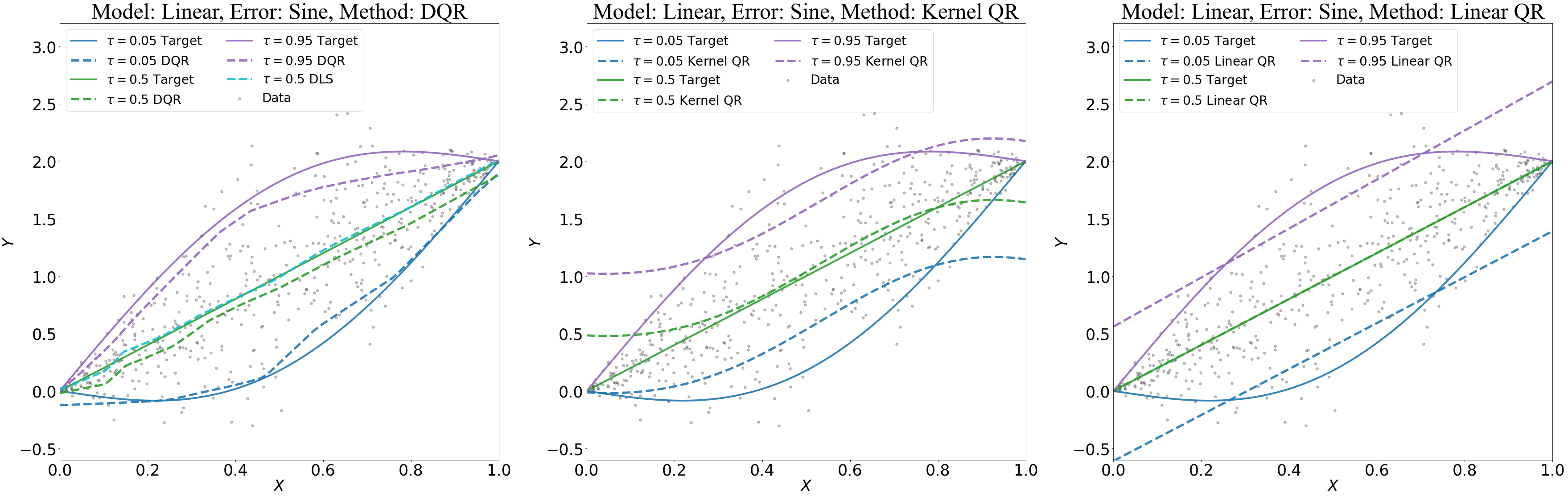}
	\end{subfigure}
	\begin{subfigure}{\textwidth}
		\includegraphics[width=\textwidth]{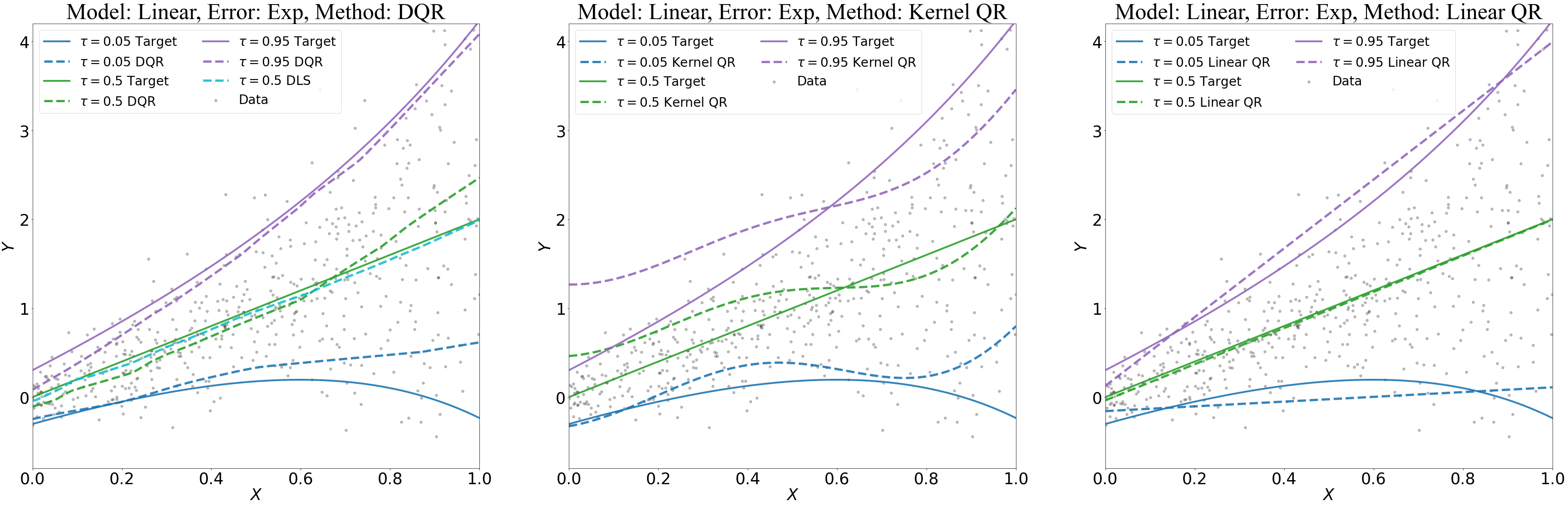}
	\end{subfigure}
	\caption{The fitted quantile curves by different methods under univariate model ``Linear'' with different errors. The training data is depicted as grey dots.The target quantile functions at $\tau=0.05,0.5,0.95$ are depicted as solid curves in different colors, and colored dashed curves represent the corresponding estimates. From the top to the bottom, each row corresponds a certain type of error: $t(3)$, ``\textit{Sine}'' and ``\textit{Exp}''. From the left to right, each column corresponds a certain estimation method: \textit{DQR}, \textit{kernel QR} and \textit{linear QR}. Fitted \textit{DLS} curves are contained in the \textit{DQR} plots.}
	\label{fitted:linear2A}
\end{figure}

\begin{figure}[H]
	\centering
	\begin{subfigure}{\textwidth}
		\includegraphics[width=1\textwidth]{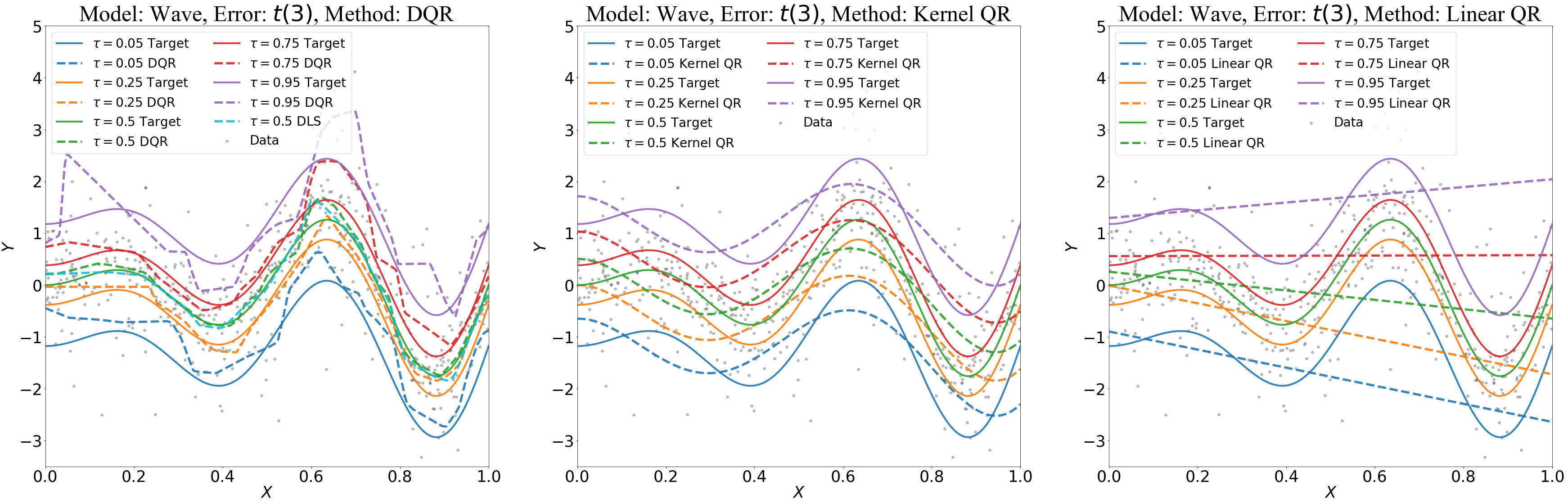}
	\end{subfigure}
	\begin{subfigure}{\textwidth}
		\includegraphics[width=1\textwidth]{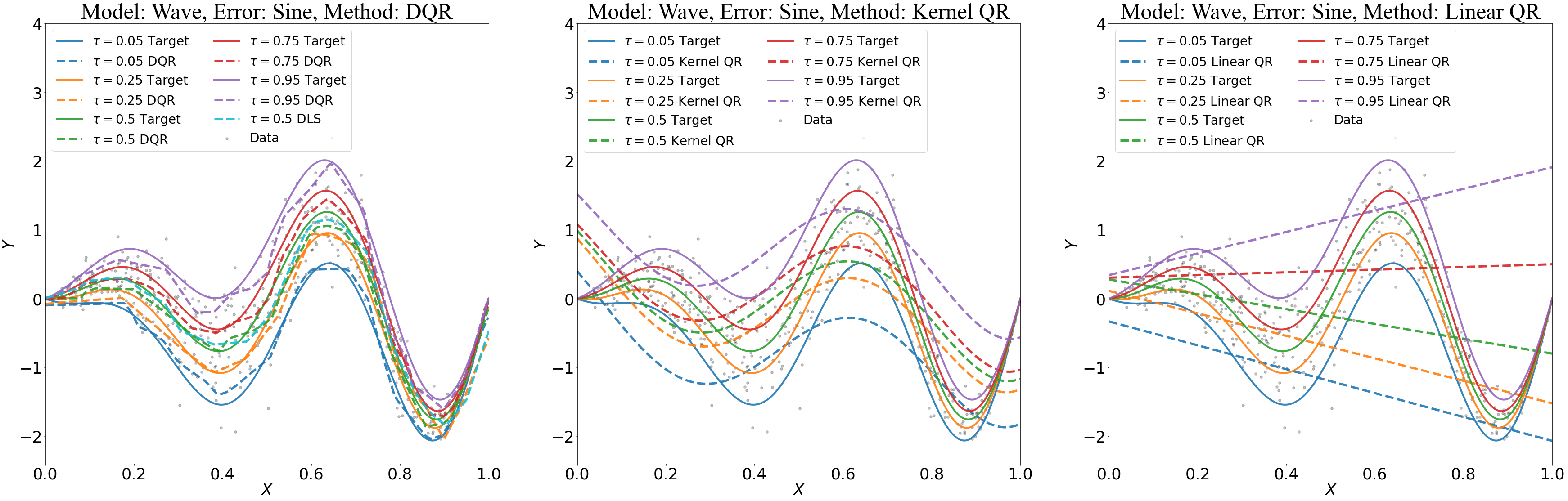}
	\end{subfigure}
	\begin{subfigure}{\textwidth}
		\includegraphics[width=\textwidth]{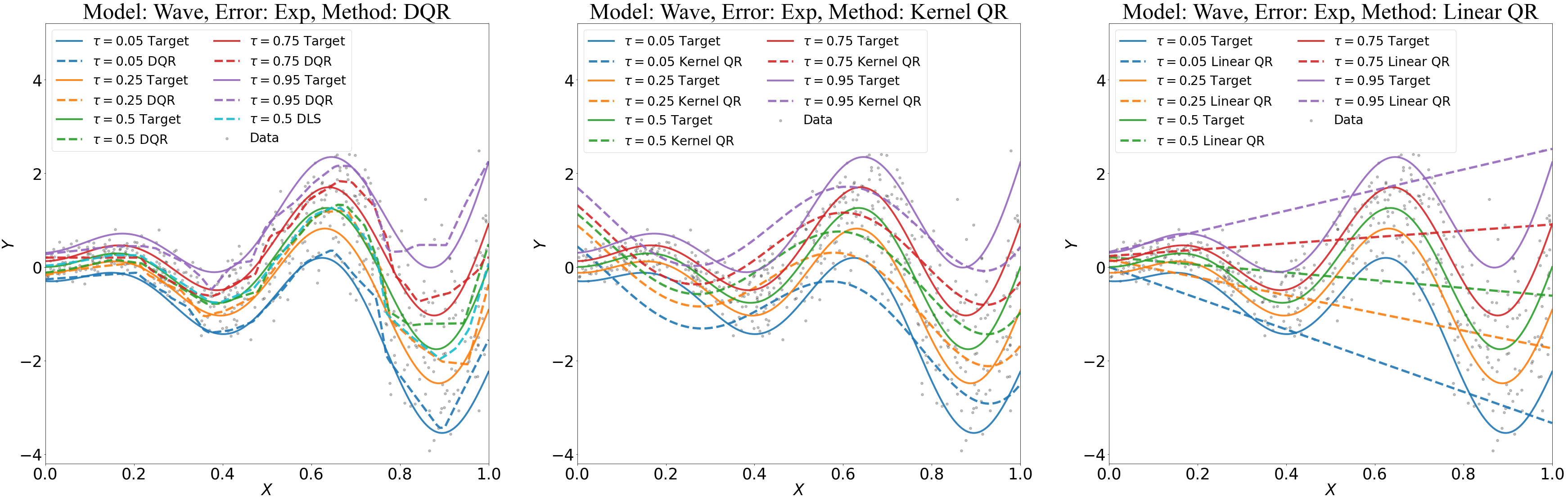}
	\end{subfigure}
	\caption{The fitted quantile curves by different methods under univariate model ``Wave'' with different errors. The training data is depicted as grey dots.The target quantile functions at $\tau=0.05,0.25,0.5,0.75,0.95$ are depicted as solid curves in different colors, and colored dashed curves represent the corresponding estimates. From the top to the bottom, each row corresponds a certain type of error: $t(3)$, ``Sine'' and ``Exp''. From the left to right, each column corresponds a certain estimation method: \textit{DQR, kernel QR} and \textit{linear QR}.}
	\label{fitted:wave1A}
\end{figure}

\begin{figure}[H]
	\centering
	\begin{subfigure}{\textwidth}
		\includegraphics[width=1\textwidth]{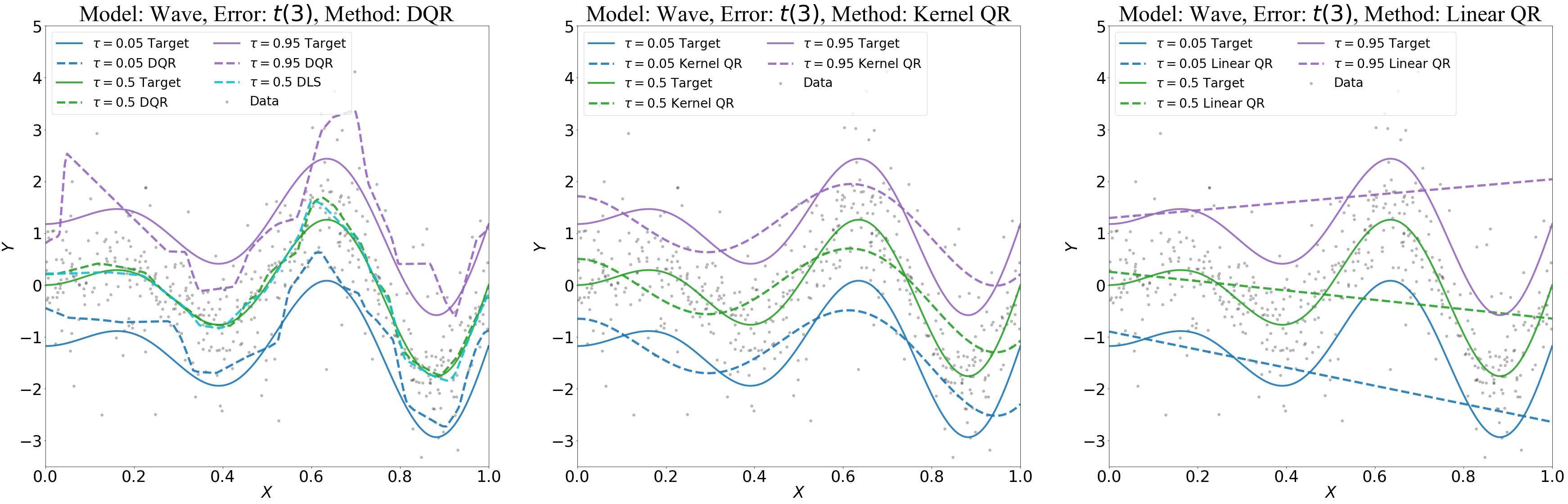}
	\end{subfigure}
	\begin{subfigure}{\textwidth}
		\includegraphics[width=1\textwidth]{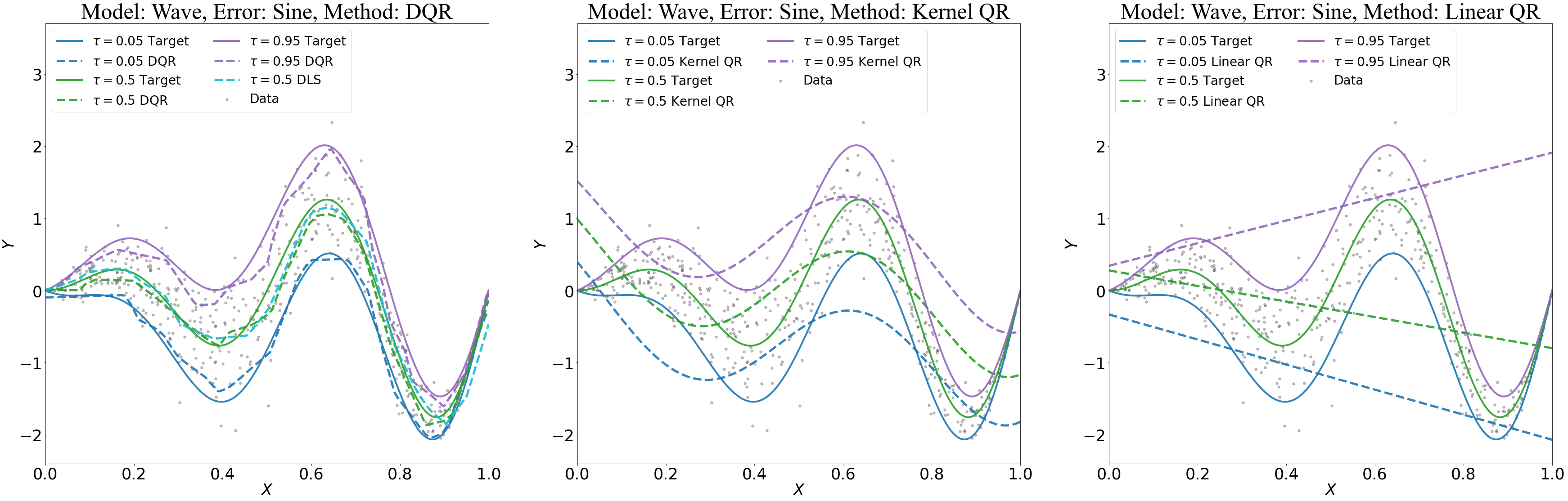}
	\end{subfigure}
	\begin{subfigure}{\textwidth}
		\includegraphics[width=\textwidth]{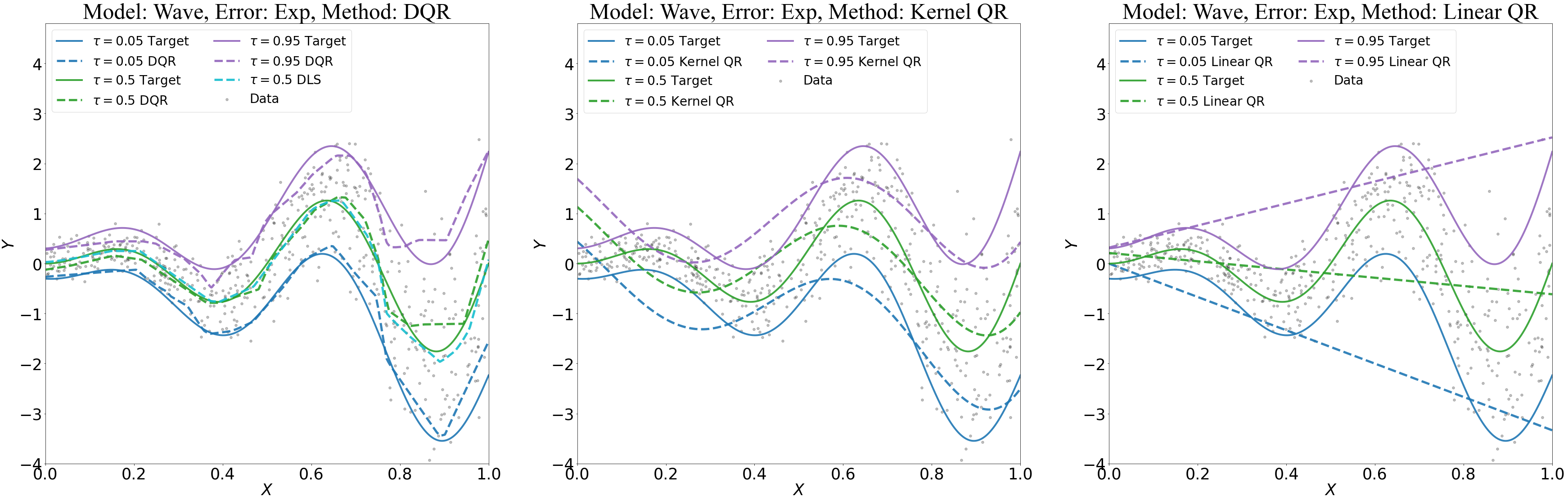}
	\end{subfigure}
	\caption{The fitted quantile curves by different methods under univariate model ``Wave'' with different errors. The training data is depicted as grey dots.The target quantile functions at $\tau=0.05,0.5,0.95$ are depicted as solid curves in different colors, and colored dashed curves represent the corresponding estimates. From the top to the bottom, each row corresponds a certain type of error: $t(3)$, ``\textit{Sine}'' and ``\textit{Exp}''. From the left to right, each column corresponds a certain estimation method: \textit{DQR}, \textit{kernel QR} and \textit{linear QR}. Fitted \textit{DLS} curves are contained in the \textit{DQR} plots.}
	\label{fitted:wave2A}
\end{figure}

\begin{figure}[H]
	\centering
	\begin{subfigure}{\textwidth}
		\includegraphics[width=1\textwidth]{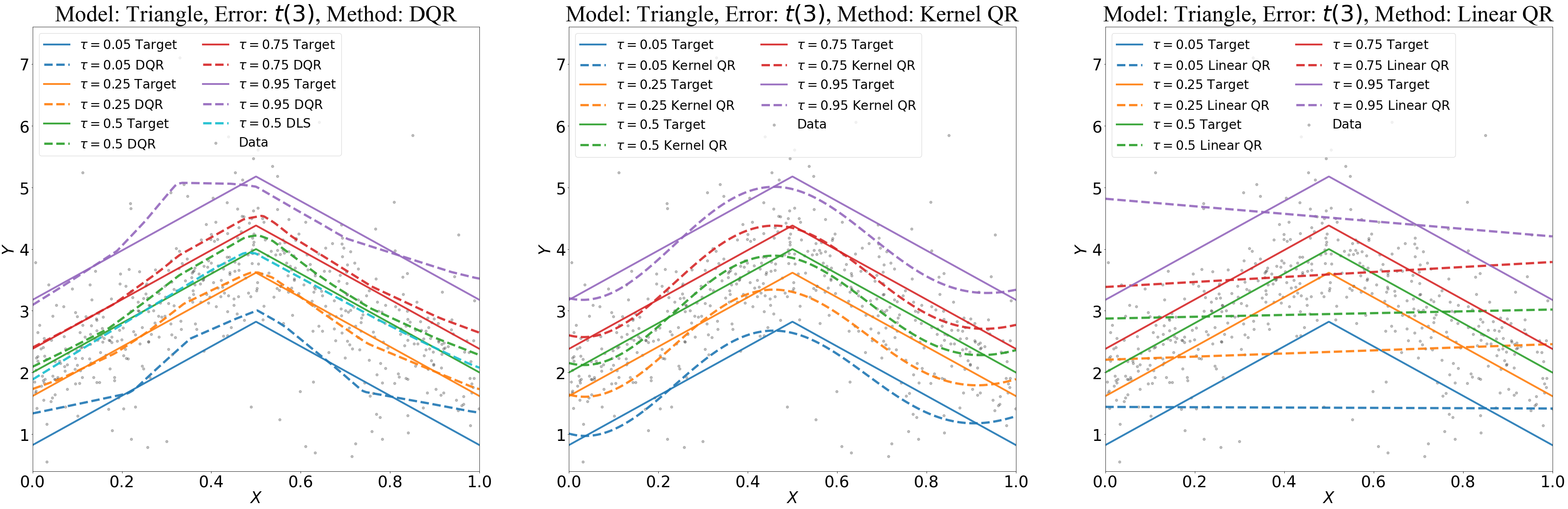}
	\end{subfigure}
	
	\begin{subfigure}{\textwidth}
		\includegraphics[width=1\textwidth]{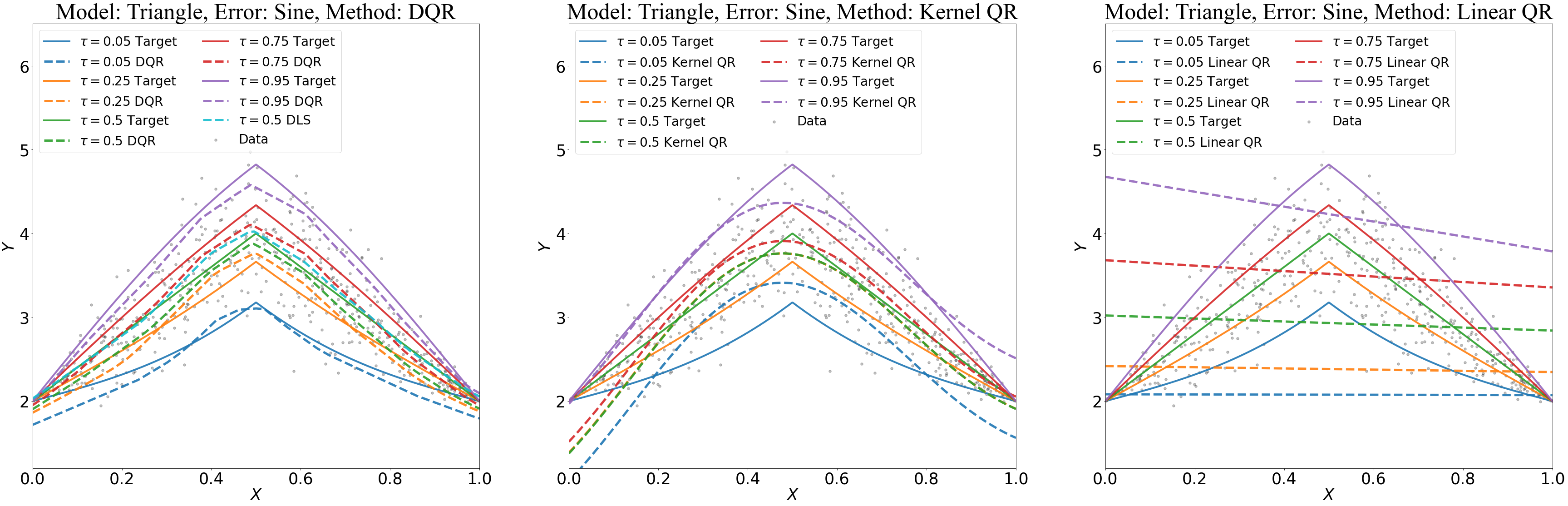}
	\end{subfigure}
	
	\begin{subfigure}{\textwidth}
		\includegraphics[width=\textwidth]{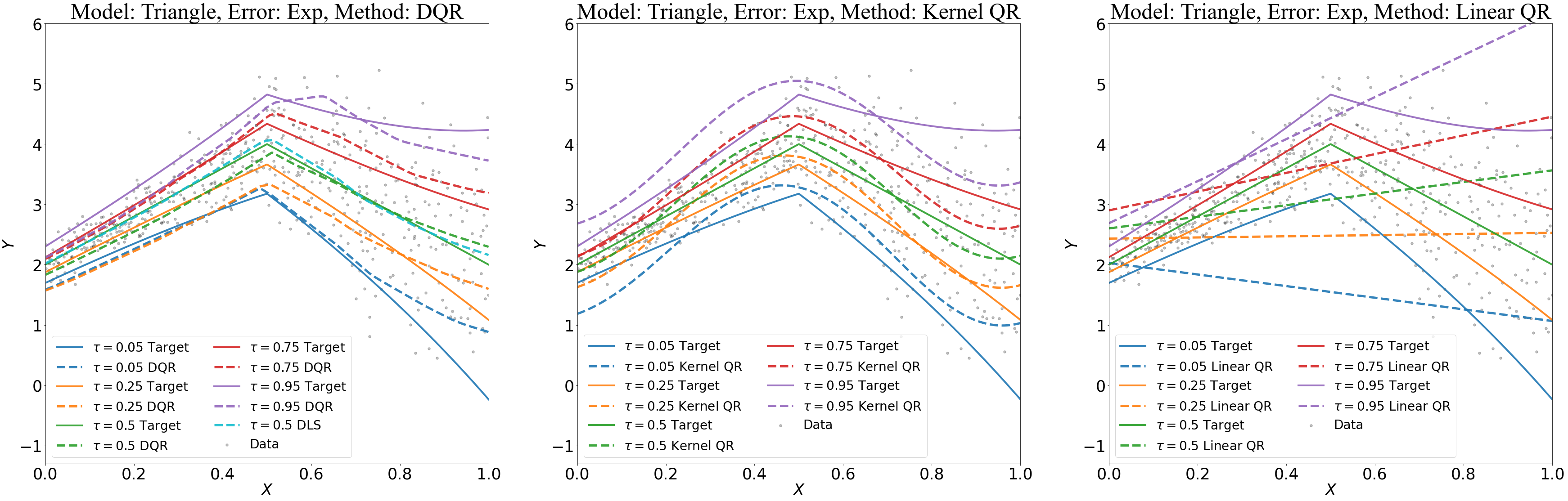}
	\end{subfigure}
	\caption{The fitted quantile curves by different methods under univariate model ``Triangle'' with different errors. The training data is depicted as grey dots.The target quantile functions at $\tau=0.05,0.25,0.5,0.75,0.95$ are depicted as solid curves in different colors, and colored dashed curves represent the corresponding estimates. From the top to the bottom, each row corresponds a certain type of error: $t(3)$, ``\textit{Sine}'' and ``\textit{Exp}''. From the left to right, each column corresponds a certain estimation method: \textit{DQR}, \textit{kernel QR} and \textit{linear QR}. Fitted \textit{DLS} curves are contained in the \textit{DQR} plots.}
	\label{fitted:triangle1A}
\end{figure}

\begin{figure}[H]
	\centering
	\begin{subfigure}{\textwidth}
		\includegraphics[width=1\textwidth]{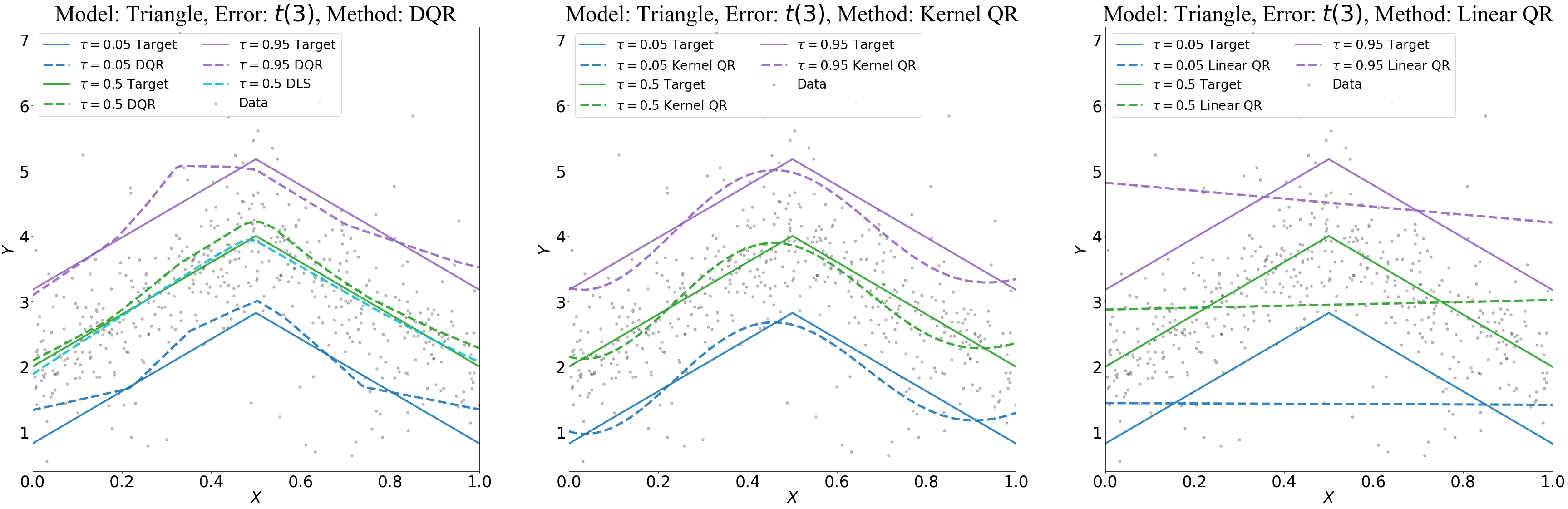}
	\end{subfigure}
	
	\begin{subfigure}{\textwidth}
		\includegraphics[width=1\textwidth]{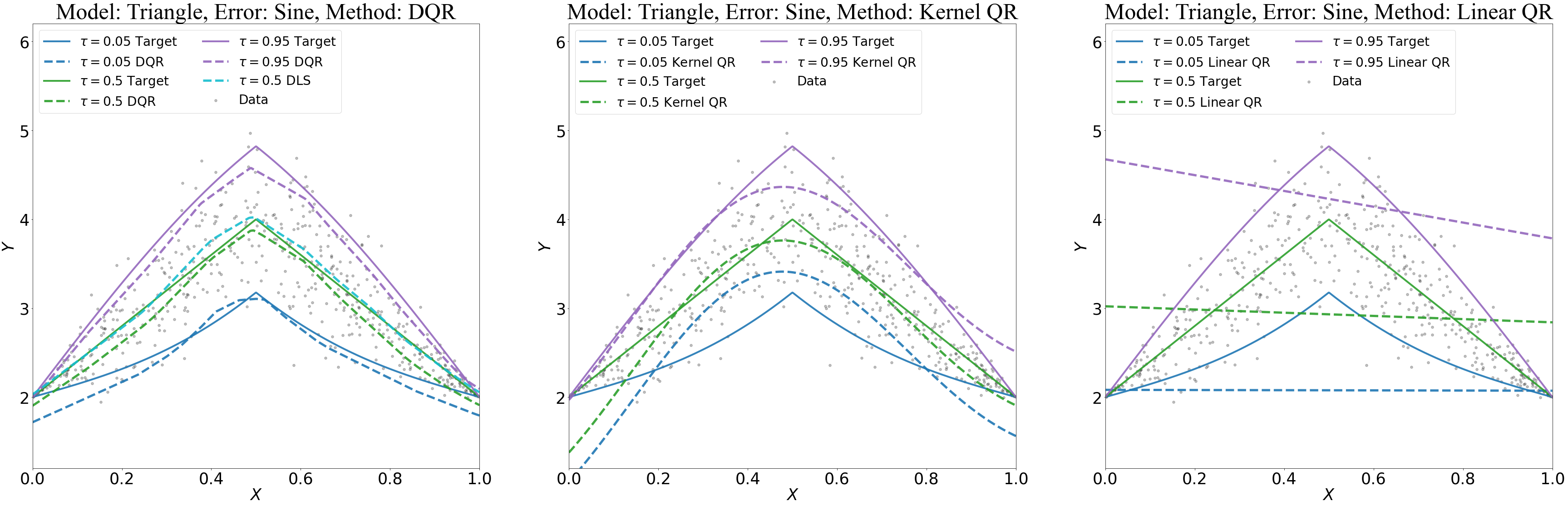}
	\end{subfigure}
	
	\begin{subfigure}{\textwidth}
		\includegraphics[width=\textwidth]{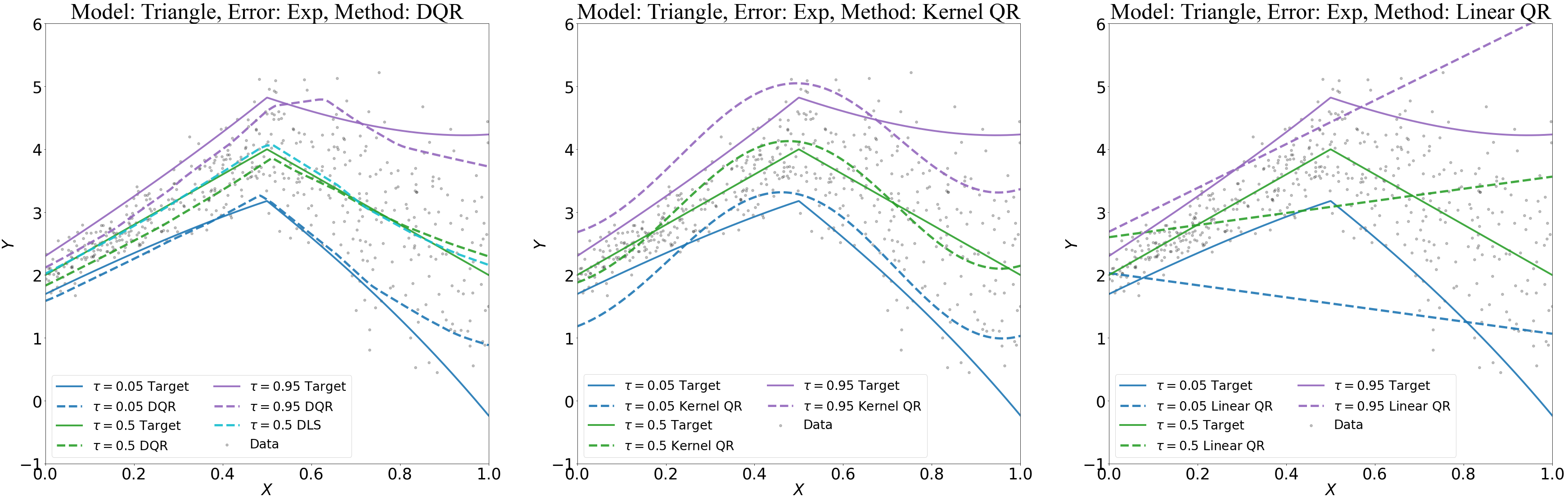}
	\end{subfigure}
	\caption{The fitted quantile curves by different methods under univariate model ``Triangle'' with different errors. The training data is depicted as grey dots.The target quantile functions at $\tau=0.05,0.5,0.95$ are depicted as solid curves in different colors, and colored dashed curves represent the corresponding estimates. From the top to the bottom, each row corresponds a certain type of error: $t(3)$, ``\textit{Sine}'' and ``\textit{Exp}''. From the left to right, each column corresponds a certain estimation method: \textit{DQR}, \textit{kernel QR} and \textit{linear QR}. Fitted \textit{DLS} curves are contained in the \textit{DQR} plots.}
	\label{fitted:triangle2A}
\end{figure}

We generate training data with sample sizes $n=128, 512$ and set the batch size of Adam optimization to be $n/2$. In all settings, we implement the empirical risk minimization of \textit{DQR} and DLS by ReLU activated fixed width multilayer perceptrons, i.e., a class of 
ReLU activated multilayer perceptrons with 4 hidden layers, the width of the network are set to be $(1,256,256,256,256,1)$. All weights and biases in each layer are initialized by uniformly samples on bounded intervals according to the default initialization mechanism in \textit{PyTorch}. The fitted quantiles curves at $\tau=0.05, 0.25,0.5,0.75, 0.95$ are shown in Figure \ref{fitted:linear1A}-\ref{fitted:triangle2A}.
Summary measures including the excess risks and the $L_1$ test and the $L_2^2$ errors
based on $R=10$ replications are summarized
are summarized in Tables \ref{tab:linearA}-\ref{tab:triangleA}.

\begin{table}[H]
	\setlength{\tabcolsep}{2pt} 
	\renewcommand{\arraystretch}{1.4} 
	\centering
	\caption{\footnotesize Data is generated from ``Linear" model with training sample size $n = 128$ or $512$ and the number of replications $R = 10$. The averaged excess risks, the $L_1$ and the $L_2^2$ test errors with the corresponding standard deviations (in parentheses) are reported for the estimators trained by different methods.}
	\label{tab:linearA}
	\resizebox{\textwidth}{!}{%
		\begin{tabular}{@{}ll|ccc|ccc|ccc@{}}
			\toprule
			\multicolumn{2}{c|}{$n=128$}              & \multicolumn{3}{c|}{$t(3)$}                                            & \multicolumn{3}{c|}{\textit{Sine}}                                           & \multicolumn{3}{c}{\textit{Exp}}                                            \\
			Quantile                     & Method    & Excess risk           & $L_1$                 & $L_2^2$               & Excess risk           & $L_1$                 & $L_2^2$               & Excess risk           & $L_1$                 & $L_2^2$               \\ \midrule
			\multirow{3}{*}{$\tau=0.05$}
			& DQR & 0.028(0.016)          & 0.563(0.074)          & 0.450(0.114)          & 0.022(0.011)          & \textbf{0.225(0.043)} & 0.103(0.038)          & 0.033(0.019)          & 0.463(0.085)          & 0.451(0.181)          \\
			& Kernel QR & 0.003(0.009)          & 0.223(0.118)          & 0.083(0.075)          & 0.094(0.084)          & 0.478(0.167)          & 0.333(0.181)          & \textbf{0.001(0.008)} & 0.339(0.102)          & 0.194(0.098)          \\
			& Linear QR & \textbf{0.002(0.010)} & \textbf{0.182(0.073)} & \textbf{0.049(0.041)} & \textbf{0.008(0.003)} & 0.244(0.021)          & \textbf{0.100(0.035)} & 0.002(0.014)          & \textbf{0.183(0.065)} & \textbf{0.063(0.050)} \\ \midrule
			\multirow{3}{*}{$\tau=0.25$}
			& DQR & 0.058(0.029)          & 0.305(0.090)          & 0.187(0.093)          & 0.021(0.009)          & 0.160(0.027)          & 0.049(0.020)          & 0.037(0.026)          & 0.302(0.068)          & 0.197(0.083)          \\
			& Kernel QR & 0.039(0.025)          & 0.255(0.082)          & 0.098(0.059)          & 0.083(0.034)          & 0.323(0.058)          & 0.138(0.041)          & 0.005(0.011)          & 0.172(0.061)          & 0.045(0.031)          \\
			& Linear QR & \textbf{0.007(0.007)} & \textbf{0.084(0.039)} & \textbf{0.010(0.008)} & \textbf{0.006(0.003)} & \textbf{0.109(0.015)} & \textbf{0.016(0.004)} & \textbf{0.003(0.016)} & \textbf{0.085(0.043)} & \textbf{0.013(0.013)} \\ \midrule
			\multirow{4}{*}{$\tau=0.5$}  & DLS       & 0.284(0.134)          & 0.366(0.078)          & 0.284(0.133)          & 0.057(0.044)          & 0.146(0.031)          & 0.057(0.043)          & 0.176(0.071)          & 0.269(0.052)          & 0.174(0.070)          \\
			& DQR & 0.097(0.042)          & 0.379(0.112)          & 0.313(0.177)          & 0.015(0.007)          & 0.155(0.035)          & 0.044(0.017)          & 0.054(0.025)          & 0.282(0.048)          & 0.168(0.068)          \\
			& Kernel QR & 0.028(0.013)          & 0.230(0.083)          & 0.086(0.046)          & 0.056(0.039)          & 0.220(0.097)          & 0.090(0.057)          & 0.024(0.010)          & 0.168(0.035)          & 0.040(0.012)          \\
			& Linear QR & \textbf{0.004(0.004)} & \textbf{0.071(0.048)} & \textbf{0.008(0.010)} & \textbf{0.001(0.003)} & \textbf{0.022(0.021)} & \textbf{0.001(0.002)} & \textbf{0.006(0.009)} & \textbf{0.070(0.035)} & \textbf{0.008(0.007)} \\ \midrule
			\multirow{3}{*}{$\tau=0.75$}
			& DQR & 0.077(0.049)          & 0.392(0.110)          & 0.327(0.206)          & 0.012(0.013)          & 0.196(0.046)          & 0.070(0.034)          & 0.046(0.017)          & 0.329(0.052)          & 0.200(0.054)          \\
			& Kernel QR & 0.011(0.007)          & 0.195(0.077)          & 0.063(0.047)          & 0.049(0.034)          & 0.326(0.129)          & 0.205(0.182)          & 0.028(0.026)          & 0.198(0.076)          & 0.072(0.067)          \\
			& Linear QR & \textbf{0.002(0.002)} & \textbf{0.087(0.039)} & \textbf{0.012(0.008)} & \textbf{0.003(0.003)} & \textbf{0.108(0.012)} & \textbf{0.016(0.004)} & \textbf{0.007(0.005)} & \textbf{0.121(0.062)} & \textbf{0.026(0.023)} \\ \midrule
			\multirow{3}{*}{$\tau=0.95$}
			& DQR & 0.073(0.034)          & 0.684(0.145)          & 0.669(0.256)          & 0.014(0.008)          & 0.310(0.050)          & 0.155(0.048)          & 0.038(0.021)          & 0.448(0.094)          & 0.361(0.140)          \\
			& Kernel QR & 0.020(0.022)          & 0.331(0.164)          & 0.183(0.158)          & 0.029(0.015)          & 0.429(0.103)          & 0.274(0.123)          & 0.023(0.009)          & 0.444(0.069)          & 0.283(0.100)          \\
			& Linear QR & \textbf{0.003(0.005)} & \textbf{0.224(0.088)} & \textbf{0.077(0.064)} & \textbf{0.006(0.002)} & \textbf{0.237(0.017)} & \textbf{0.090(0.024)} & \textbf{0.005(0.004)} & \textbf{0.187(0.094)} & \textbf{0.058(0.062)} \\ \midrule
			\multicolumn{2}{c|}{$n=512$}              & \multicolumn{3}{c|}{$t(3)$}                                            & \multicolumn{3}{c|}{\textit{Sine}}                                           & \multicolumn{3}{c}{\textit{Exp}}                                            \\
			Quantile                     & Method    & Excess risk           & $L_1$                 & $L_2^2$               & Excess risk           & $L_1$                 & $L_2^2$               & Excess risk           & $L_1$                 & $L_2^2$               \\ \midrule
			\multirow{3}{*}{$\tau=0.05$}
			& DQR & \textbf{0.001(0.005)} & 0.401(0.038)          & 0.223(0.043)          & \textbf{0.003(0.002)} & \textbf{0.118(0.022)} & \textbf{0.026(0.011)} & \textbf{0.001(0.004)} & 0.303(0.039)          & 0.209(0.050)          \\
			& Kernel QR & 0.006(0.007)          & 0.203(0.059)          & 0.065(0.035)          & 0.012(0.007)          & 0.273(0.089)          & 0.145(0.101)          & 0.005(0.005)          & 0.266(0.096)          & 0.140(0.067)          \\
			& Linear QR & 0.002(0.009)          & \textbf{0.137(0.070)} & \textbf{0.032(0.035)} & 0.005(0.001)          & 0.224(0.006)          & 0.079(0.010)          & 0.001(0.005)          & \textbf{0.132(0.016)} & \textbf{0.026(0.005)} \\ \midrule
			\multirow{3}{*}{$\tau=0.25$}
			& DQR & 0.019(0.012)          & 0.192(0.044)          & 0.070(0.036)          & 0.004(0.003)          & \textbf{0.083(0.020)} & 0.013(0.005)          & \textbf{0.002(0.008)} & 0.207(0.030)          & 0.094(0.041)          \\
			& Kernel QR & 0.027(0.011)          & 0.180(0.057)          & 0.053(0.027)          & 0.034(0.018)          & 0.177(0.051)          & 0.048(0.021)          & 0.002(0.009)          & 0.159(0.067)          & 0.043(0.031)          \\
			& Linear QR & \textbf{0.001(0.003)} & \textbf{0.042(0.015)} & \textbf{0.003(0.002)} & \textbf{0.004(0.001)} & 0.099(0.006)          & \textbf{0.013(0.001)} & 0.005(0.004)          & \textbf{0.056(0.012)} & \textbf{0.007(0.004)} \\ \midrule
			\multirow{4}{*}{$\tau=0.5$}  & DLS       & 0.074(0.040)          & 0.186(0.062)          & 0.074(0.040)          & 0.010(0.006)          & 0.066(0.019)          & 0.010(0.006)          & 0.045(0.022)          & 0.136(0.041)          & 0.046(0.023)          \\
			& DQR & 0.030(0.013)          & 0.195(0.036)          & 0.084(0.046)          & 0.002(0.004)          & 0.110(0.023)          & 0.020(0.008)          & 0.001(0.006)          & 0.157(0.022)          & 0.050(0.018)          \\
			& Kernel QR & 0.016(0.014)          & 0.137(0.062)          & 0.035(0.027)          & 0.037(0.029)          & 0.171(0.076)          & 0.058(0.047)          & 0.020(0.018)          & 0.133(0.049)          & 0.032(0.020)          \\
			& Linear QR & \textbf{0.001(0.001)} & \textbf{0.036(0.016)} & \textbf{0.002(0.002)} & \textbf{0.001(0.002)} & \textbf{0.009(0.008)} & \textbf{0.001(0.001)} & \textbf{0.000(0.003)} & \textbf{0.029(0.017)} & \textbf{0.002(0.001)} \\ \midrule
			\multirow{3}{*}{$\tau=0.75$}
			& DQR & 0.032(0.006)          & 0.277(0.024)          & 0.165(0.040)          & 0.003(0.003)          & 0.140(0.026)          & 0.028(0.011)          & 0.006(0.007)          & 0.190(0.032)          & 0.065(0.029)          \\
			& Kernel QR & 0.008(0.007)          & 0.146(0.054)          & 0.037(0.025)          & 0.015(0.013)          & 0.200(0.076)          & 0.059(0.035)          & 0.017(0.011)          & 0.141(0.059)          & 0.037(0.021)          \\
			& Linear QR & \textbf{0.001(0.000)} & \textbf{0.033(0.012)} & \textbf{0.002(0.001)} & \textbf{0.003(0.002)} & \textbf{0.100(0.006)} & \textbf{0.013(0.002)} & \textbf{0.004(0.001)} & \textbf{0.060(0.012)} & \textbf{0.006(0.003)} \\ \midrule
			\multirow{3}{*}{$\tau=0.95$}
			& DQR & 0.028(0.011)          & 0.598(0.123)          & 1.281(1.633)          & \textbf{0.001(0.003)} & \textbf{0.205(0.034)} & \textbf{0.063(0.025)} & 0.005(0.004)          & 0.266(0.056)          & 0.119(0.049)          \\
			& Kernel QR & 0.004(0.003)          & 0.168(0.061)          & 0.047(0.026)          & 0.006(0.005)          & 0.247(0.090)          & 0.101(0.069)          & 0.010(0.006)          & 0.307(0.089)          & 0.148(0.081)          \\
			& Linear QR & \textbf{0.001(0.002)} & \textbf{0.137(0.082)} & \textbf{0.036(0.041)} & 0.006(0.001)          & 0.224(0.005)          & 0.077(0.007)          & \textbf{0.003(0.001)} & \textbf{0.124(0.010)} & \textbf{0.026(0.010)} \\ \bottomrule
		\end{tabular}%
	}
\end{table}

\begin{table}[H]
	\setlength{\tabcolsep}{2pt} 
	\renewcommand{\arraystretch}{1.4} 
	\centering
	\caption{\footnotesize Data is generated from ``Wave" model with training sample size $n = 128$ or $512$ and the number of replications $R = 10$. The averaged excess risks, the $L_1$ and the $L_2^2$ test errors with the corresponding standard deviations (in parentheses) are reported for the estimators trained by different methods.}
	\label{tab:waveA}
	\resizebox{\textwidth}{!}{%
		\begin{tabular}{@{}ll|ccc|ccc|ccc@{}}
			\toprule
			\multicolumn{2}{c|}{$n=128$}              & \multicolumn{3}{c|}{$t(3)$}                                            & \multicolumn{3}{c|}{\textit{Sine}}                                            & \multicolumn{3}{c}{\textit{Exp}}                                             \\
			Quantile                     & Method    & Excess risk           & $L_1$                 & $L_2^2$               & Excess risk           & $L_1$                 & $L_2^2$               & Excess risk           & $L_1$                 & $L_2^2$               \\ \midrule
			\multirow{3}{*}{$\tau=0.05$}
			& DQR & 0.040(0.026)          & 0.634(0.100)          & 0.593(0.249)          & \textbf{0.027(0.014)} & \textbf{0.253(0.061)} & \textbf{0.122(0.056)} & \textbf{0.023(0.014)} & \textbf{0.452(0.082)} & \textbf{0.409(0.130)} \\
			& Kernel QR & \textbf{0.033(0.004)} & \textbf{0.511(0.037)} & \textbf{0.366(0.051)} & 0.068(0.010)          & 0.542(0.031)          & 0.415(0.049)          & 0.048(0.007)          & 0.636(0.074)          & 0.570(0.151)          \\
			& Linear QR & 0.088(0.026)          & 0.672(0.063)          & 0.847(0.208)          & 0.079(0.007)          & 0.713(0.056)          & 0.888(0.103)          & 0.078(0.012)          & 0.660(0.016)          & 0.816(0.072)          \\ \midrule
			\multirow{3}{*}{$\tau=0.25$}
			& DQR & \textbf{0.072(0.038)} & \textbf{0.339(0.071)} & \textbf{0.251(0.143)} & \textbf{0.021(0.009)} & \textbf{0.161(0.032)} & \textbf{0.050(0.021)} & \textbf{0.048(0.027)} & \textbf{0.350(0.059)} & \textbf{0.287(0.109)} \\
			& Kernel QR & 0.130(0.014)          & 0.506(0.023)          & 0.358(0.037)          & 0.169(0.007)          & 0.517(0.015)          & 0.358(0.024)          & 0.126(0.015)          & 0.529(0.023)          & 0.365(0.034)          \\
			& Linear QR & 0.245(0.021)          & 0.608(0.017)          & 0.707(0.067)          & 0.252(0.011)          & 0.612(0.012)          & 0.745(0.048)          & 0.229(0.044)          & 0.610(0.017)          & 0.690(0.086)          \\ \midrule
			\multirow{4}{*}{$\tau=0.5$}  & DLS       & 0.203(0.063)          & \textbf{0.327(0.047)} & \textbf{0.203(0.064)} & 0.050(0.023)          & \textbf{0.148(0.031)} & \textbf{0.050(0.023)} & 0.205(0.060)          & 0.296(0.049)          & 0.206(0.061)          \\
			& DQR & \textbf{0.099(0.054)} & 0.348(0.101)          & 0.412(0.573)          & \textbf{0.020(0.006)} & 0.179(0.020)          & 0.058(0.013)          & \textbf{0.050(0.024)} & \textbf{0.290(0.059)} & \textbf{0.187(0.084)} \\
			& Kernel QR & 0.145(0.020)          & 0.504(0.043)          & 0.373(0.092)          & 0.172(0.011)          & 0.519(0.011)          & 0.360(0.022)          & 0.157(0.028)          & 0.526(0.019)          & 0.368(0.038)          \\
			& Linear QR & 0.247(0.023)          & 0.595(0.016)          & 0.586(0.019)          & 0.280(0.016)          & 0.583(0.002)          & 0.581(0.019)          & 0.210(0.027)          & 0.597(0.017)          & 0.595(0.029)          \\ \midrule
			\multirow{3}{*}{$\tau=0.75$}
			& DQR & \textbf{0.090(0.038)} & \textbf{0.443(0.134)} & 0.376(0.274)          & \textbf{0.013(0.008)} & \textbf{0.202(0.030)} & \textbf{0.070(0.023)} & \textbf{0.065(0.029)} & \textbf{0.347(0.055)} & \textbf{0.239(0.076)} \\
			& Kernel QR & 0.095(0.015)          & 0.516(0.019)          & \textbf{0.365(0.038)} & 0.127(0.016)          & 0.516(0.015)          & 0.358(0.024)          & 0.128(0.032)          & 0.523(0.022)          & 0.403(0.062)          \\
			& Linear QR & 0.135(0.010)          & 0.679(0.040)          & 0.712(0.093)          & 0.176(0.005)          & 0.756(0.050)          & 1.064(0.198)          & 0.124(0.013)          & 0.629(0.030)          & 0.653(0.075)          \\ \midrule
			\multirow{3}{*}{$\tau=0.95$}
			& DQR & 0.077(0.030)          & 0.766(0.249)          & 2.252(4.679)          & \textbf{0.016(0.011)} & \textbf{0.304(0.058)} & \textbf{0.155(0.057)} & 0.055(0.027)          & \textbf{0.515(0.091)} & \textbf{0.451(0.166)} \\
			& Kernel QR & 0.026(0.011)          & \textbf{0.528(0.027)} & \textbf{0.397(0.047)} & 0.052(0.042)          & 0.578(0.046)          & 0.494(0.080)          & \textbf{0.041(0.007)} & 0.565(0.044)          & 0.498(0.087)          \\
			& Linear QR & \textbf{0.021(0.005)} & 0.866(0.133)          & 1.122(0.285)          & 0.047(0.004)          & 0.987(0.087)          & 2.085(0.458)          & 0.029(0.001)          & 0.660(0.044)          & 0.723(0.102)          \\ \midrule
			\multicolumn{2}{c|}{$n=512$}              & \multicolumn{3}{c|}{$t(3)$}                                            & \multicolumn{3}{c|}{\textit{Sine}}                                            & \multicolumn{3}{c}{\textit{Exp}}                                             \\
			Quantile                     & Method    & Excess risk           & $L_1$                 & $L_2^2$               & Excess risk           & $L_1$                 & $L_2^2$               & Excess risk           & $L_1$                 & $L_2^2$               \\ \midrule
			\multirow{3}{*}{$\tau=0.05$}
			& DQR & \textbf{0.004(0.004)} & \textbf{0.367(0.049)} & \textbf{0.179(0.048)} & \textbf{0.004(0.002)} & \textbf{0.139(0.026)} & \textbf{0.033(0.012)} & \textbf{0.006(0.005)} & \textbf{0.277(0.040)} & \textbf{0.169(0.060)} \\
			& Kernel QR & 0.025(0.007)          & 0.427(0.027)          & 0.259(0.042)          & 0.061(0.011)          & 0.473(0.019)          & 0.320(0.021)          & 0.036(0.012)          & 0.494(0.036)          & 0.323(0.045)          \\
			& Linear QR & 0.082(0.009)          & 0.622(0.019)          & 0.770(0.073)          & 0.075(0.002)          & 0.684(0.026)          & 0.831(0.049)          & 0.091(0.011)          & 0.665(0.021)          & 0.886(0.080)          \\ \midrule
			\multirow{3}{*}{$\tau=0.25$}
			& DQR & \textbf{0.015(0.007)} & \textbf{0.183(0.037)} & \textbf{0.055(0.020)} & \textbf{0.005(0.002)} & \textbf{0.098(0.015)} & \textbf{0.016(0.004)} & \textbf{0.005(0.006)} & \textbf{0.203(0.036)} & \textbf{0.096(0.037)} \\
			& Kernel QR & 0.091(0.008)          & 0.424(0.021)          & 0.244(0.020)          & 0.160(0.017)          & 0.478(0.035)          & 0.311(0.052)          & 0.102(0.022)          & 0.447(0.019)          & 0.262(0.031)          \\
			& Linear QR & 0.240(0.017)          & 0.597(0.011)          & 0.684(0.051)          & 0.243(0.007)          & 0.602(0.009)          & 0.711(0.035)          & 0.227(0.019)          & 0.602(0.007)          & 0.675(0.037)          \\ \midrule
			\multirow{4}{*}{$\tau=0.5$}  & DLS       & 0.055(0.019)          & \textbf{0.174(0.029)} & \textbf{0.055(0.019)} & 0.013(0.004)          & \textbf{0.082(0.011)} & \textbf{0.013(0.004)} & 0.044(0.023)          & \textbf{0.135(0.026)} & \textbf{0.044(0.024)} \\
			& DQR & \textbf{0.029(0.012)} & 0.200(0.034)          & 0.099(0.092)          & \textbf{0.002(0.003)} & 0.119(0.015)          & 0.022(0.006)          & \textbf{0.014(0.018)} & 0.216(0.041)          & 0.102(0.045)          \\
			& Kernel QR & 0.096(0.018)          & 0.415(0.042)          & 0.237(0.055)          & 0.152(0.014)          & 0.476(0.019)          & 0.310(0.032)          & 0.113(0.012)          & 0.449(0.014)          & 0.260(0.019)          \\
			& Linear QR & 0.242(0.009)          & 0.581(0.004)          & 0.565(0.004)          & 0.280(0.010)          & 0.579(0.002)          & 0.577(0.013)          & 0.216(0.015)          & 0.584(0.007)          & 0.569(0.008)          \\ \midrule
			\multirow{3}{*}{$\tau=0.75$}
			& DQR & \textbf{0.040(0.014)} & \textbf{0.326(0.059)} & 0.306(0.214)          & \textbf{0.001(0.003)} & \textbf{0.149(0.024)} & \textbf{0.035(0.012)} & \textbf{0.012(0.007)} & \textbf{0.223(0.030)} & \textbf{0.084(0.020)} \\
			& Kernel QR & 0.065(0.013)          & 0.418(0.023)          & \textbf{0.244(0.025)} & 0.116(0.015)          & 0.494(0.026)          & 0.333(0.044)          & 0.088(0.010)          & 0.464(0.021)          & 0.295(0.036)          \\
			& Linear QR & 0.132(0.008)          & 0.661(0.028)          & 0.718(0.083)          & 0.175(0.006)          & 0.726(0.034)          & 0.949(0.128)          & 0.123(0.004)          & 0.632(0.022)          & 0.663(0.058)          \\ \midrule
			\multirow{3}{*}{$\tau=0.95$}
			& DQR & 0.033(0.015)          & 0.664(0.309)          & 2.215(4.597)          & \textbf{0.002(0.003)} & \textbf{0.219(0.037)} & \textbf{0.073(0.025)} & \textbf{0.011(0.006)} & \textbf{0.311(0.055)} & \textbf{0.172(0.058)} \\
			& Kernel QR & 0.019(0.005)          & \textbf{0.432(0.025)} & \textbf{0.254(0.034)} & 0.036(0.010)          & 0.544(0.036)          & 0.448(0.061)          & 0.029(0.005)          & 0.497(0.047)          & 0.374(0.076)          \\
			& Linear QR & \textbf{0.018(0.002)} & 0.796(0.073)          & 1.108(0.183)          & 0.045(0.001)          & 0.977(0.034)          & 2.108(0.174)          & 0.030(0.001)          & 0.729(0.039)          & 0.935(0.115)          \\ \bottomrule
		\end{tabular}%
	}
\end{table}

\begin{table}[H]
	\setlength{\tabcolsep}{2pt} 
	\renewcommand{\arraystretch}{1.4} 
	\centering
	\caption{\footnotesize Data is generated from ``Triangle" model with training sample size $n = 128$ or $512$ and the number of replications $R = 10$. The averaged excess risks, the $L_1$ and the $L_2^2$ test errors with the corresponding standard deviations (in parentheses) are reported for the estimators trained by different methods.}
	\label{tab:triangleA}
	\resizebox{\textwidth}{!}{%
		\begin{tabular}{@{}ll|ccc|ccc|ccc@{}}
			\toprule
			\multicolumn{2}{c|}{$n=128$}              & \multicolumn{3}{c|}{$t(3)$}                                            & \multicolumn{3}{c|}{\textit{Sine}}                                            & \multicolumn{3}{c}{\textit{Exp}}                                             \\
			Quantile                     & Method    & Excess risk           & $L_1$                 & $L_2^2$               & Excess risk           & $L_1$                 & $L_2^2$               & Excess risk           & $L_1$                 & $L_2^2$               \\ \midrule
			\multirow{3}{*}{$\tau=0.05$}
			& DQR & 0.025(0.013)          & 0.519(0.096)          & 0.386(0.130)          & \textbf{0.011(0.006)} & \textbf{0.176(0.035)} & \textbf{0.060(0.024)} & 0.013(0.014)          & \textbf{0.350(0.058)} & 0.316(0.122)          \\
			& Kernel QR & \textbf{0.012(0.014)} & \textbf{0.270(0.077)} & \textbf{0.108(0.049)} & 0.073(0.078)          & 0.409(0.153)          & 0.244(0.166)          & \textbf{0.002(0.004)} & 0.378(0.061)          & \textbf{0.237(0.080)} \\
			& Linear QR & 0.081(0.038)          & 0.650(0.160)          & 0.621(0.325)          & 0.039(0.001)          & 0.401(0.023)          & 0.263(0.024)          & 0.086(0.023)          & 0.774(0.069)          & 0.857(0.163)          \\ \midrule
			\multirow{3}{*}{$\tau=0.25$}
			& DQR & 0.045(0.025)          & 0.268(0.071)          & 0.144(0.079)          & \textbf{0.012(0.007)} & \textbf{0.122(0.028)} & \textbf{0.028(0.014)} & 0.019(0.022)          & 0.250(0.064)          & 0.176(0.094)          \\
			& Kernel QR & \textbf{0.037(0.026)} & \textbf{0.227(0.089)} & \textbf{0.083(0.062)} & 0.097(0.050)          & 0.355(0.074)          & 0.175(0.068)          & \textbf{0.009(0.012)} & \textbf{0.199(0.052)} & \textbf{0.075(0.059)} \\
			& Linear QR & 0.167(0.021)          & 0.553(0.028)          & 0.436(0.053)          & 0.166(0.015)          & 0.502(0.042)          & 0.387(0.070)          & 0.125(0.013)          & 0.586(0.017)          & 0.485(0.040)          \\ \midrule
			\multirow{4}{*}{$\tau=0.5$}  & DLS       & 0.155(0.082)          & 0.289(0.071)          & 0.155(0.082)          & 0.024(0.016)          & \textbf{0.111(0.033)} & \textbf{0.024(0.016)} & 0.110(0.049)          & 0.214(0.056)          & 0.110(0.050)          \\
			& DQR & 0.058(0.033)          & 0.272(0.092)          & 0.157(0.102)          & \textbf{0.013(0.006)} & 0.146(0.030)          & 0.039(0.014)          & 0.071(0.049)          & 0.301(0.065)          & 0.217(0.098)          \\
			& Kernel QR & \textbf{0.032(0.030)} & \textbf{0.202(0.097)} & \textbf{0.073(0.065)} & 0.053(0.031)          & 0.236(0.076)          & 0.089(0.060)          & \textbf{0.031(0.023)} & \textbf{0.192(0.074)} & \textbf{0.057(0.040)} \\
			& Linear QR & 0.140(0.010)          & 0.508(0.012)          & 0.355(0.031)          & 0.189(0.006)          & 0.519(0.009)          & 0.375(0.022)          & 0.110(0.007)          & 0.518(0.015)          & 0.378(0.039)          \\ \midrule
			\multirow{3}{*}{$\tau=0.75$}
			& DQR & 0.070(0.038)          & 0.380(0.101)          & 0.256(0.144)          & \textbf{0.003(0.005)} & \textbf{0.164(0.041)} & \textbf{0.044(0.021)} & 0.041(0.024)          & 0.306(0.081)          & 0.193(0.088)          \\
			& Kernel QR & \textbf{0.025(0.027)} & \textbf{0.228(0.106)} & \textbf{0.091(0.077)} & 0.038(0.021)          & 0.260(0.088)          & 0.103(0.060)          & \textbf{0.025(0.013)} & \textbf{0.180(0.052)} & \textbf{0.052(0.030)} \\
			& Linear QR & 0.075(0.005)          & 0.528(0.016)          & 0.405(0.041)          & 0.135(0.004)          & 0.636(0.033)          & 0.612(0.118)          & 0.066(0.005)          & 0.507(0.024)          & 0.393(0.067)          \\ \midrule
			\multirow{3}{*}{$\tau=0.95$}
			& DQR & 0.054(0.024)          & 0.583(0.123)          & 0.498(0.158)          & \textbf{0.004(0.006)} & \textbf{0.224(0.057)} & \textbf{0.084(0.046)} & 0.036(0.035)          & \textbf{0.394(0.087)} & 0.318(0.188)          \\
			& Kernel QR & \textbf{0.017(0.015)} & \textbf{0.328(0.109)} & \textbf{0.162(0.097)} & 0.053(0.054)          & 0.423(0.145)          & 0.264(0.126)          & 0.028(0.019)          & 0.447(0.077)          & \textbf{0.305(0.137)} \\
			& Linear QR & 0.015(0.007)          & 0.633(0.103)          & 0.718(0.298)          & 0.042(0.006)          & 0.907(0.118)          & 1.526(0.603)          & \textbf{0.017(0.003)} & 0.501(0.068)          & 0.422(0.151)          \\ \midrule
			\multicolumn{2}{c|}{$n=512$}              & \multicolumn{3}{c|}{$t(3)$}                                            & \multicolumn{3}{c|}{\textit{Sine}}                                            & \multicolumn{3}{c}{\textit{Exp}}                                             \\
			Quantile                     & Method    & Excess risk           & $L_1$                 & $L_2^2$               & Excess risk           & $L_1$                 & $L_2^2$               & Excess risk           & $L_1$                 & $L_2^2$               \\ \midrule
			\multirow{3}{*}{$\tau=0.05$}
			& DQR & 0.004(0.005)          & 0.409(0.049)          & 0.208(0.046)          & \textbf{0.002(0.002)} & \textbf{0.098(0.021)} & \textbf{0.021(0.009)} & \textbf{0.003(0.004)} & \textbf{0.270(0.036)} & \textbf{0.171(0.050)} \\
			& Kernel QR & \textbf{0.001(0.011)} & \textbf{0.236(0.093)} & \textbf{0.088(0.071)} & 0.020(0.009)          & 0.365(0.067)          & 0.221(0.070)          & 0.005(0.004)          & 0.330(0.056)          & 0.186(0.058)          \\
			& Linear QR & 0.058(0.009)          & 0.552(0.028)          & 0.438(0.056)          & 0.039(0.001)          & 0.407(0.014)          & 0.270(0.014)          & 0.073(0.010)          & 0.754(0.029)          & 0.798(0.068)          \\ \midrule
			\multirow{3}{*}{$\tau=0.25$}
			& DQR & \textbf{0.010(0.004)} & \textbf{0.171(0.037)} & 0.044(0.014)          & \textbf{0.003(0.002)} & \textbf{0.093(0.028)} & \textbf{0.013(0.006)} & \textbf{0.003(0.005)} & 0.202(0.020)          & 0.095(0.021)          \\
			& Kernel QR & 0.033(0.029)          & 0.200(0.102)          & \textbf{0.077(0.075)} & 0.042(0.014)          & 0.292(0.048)          & 0.123(0.036)          & 0.009(0.008)          & \textbf{0.157(0.033)} & \textbf{0.038(0.013)} \\
			& Linear QR & 0.146(0.010)          & 0.525(0.011)          & 0.384(0.023)          & 0.160(0.005)          & 0.484(0.018)          & 0.359(0.031)          & 0.132(0.015)          & 0.586(0.009)          & 0.480(0.017)          \\ \midrule
			\multirow{4}{*}{$\tau=0.5$}  & DLS       & 0.047(0.030)          & 0.157(0.044)          & 0.046(0.029)          & 0.008(0.004)          & \textbf{0.060(0.014)} & \textbf{0.008(0.004)} & 0.031(0.020)          & \textbf{0.114(0.037)} & \textbf{0.030(0.020)} \\
			& DQR & \textbf{0.017(0.007)} & \textbf{0.155(0.037)} & \textbf{0.046(0.021)} & \textbf{0.002(0.001)} & 0.095(0.022)          & 0.012(0.005)          & \textbf{0.004(0.006)} & 0.172(0.033)          & 0.061(0.030)          \\
			& Kernel QR & 0.030(0.035)          & 0.201(0.106)          & 0.072(0.076)          & 0.029(0.015)          & 0.192(0.060)          & 0.067(0.036)          & 0.031(0.030)          & 0.186(0.079)          & 0.066(0.052)          \\
			& Linear QR & 0.135(0.007)          & 0.505(0.003)          & 0.342(0.005)          & 0.186(0.005)          & 0.508(0.006)          & 0.348(0.011)          & 0.103(0.002)          & 0.523(0.005)          & 0.394(0.015)          \\ \midrule
			\multirow{3}{*}{$\tau=0.75$}
			& DQR & 0.035(0.015)          & 0.317(0.076)          & 0.305(0.277)          & \textbf{0.005(0.001)} & \textbf{0.146(0.026)} & \textbf{0.027(0.008)} & \textbf{0.004(0.007)} & 0.196(0.031)          & 0.057(0.022)          \\
			& Kernel QR & \textbf{0.018(0.022)} & \textbf{0.223(0.106)} & \textbf{0.082(0.076)} & 0.029(0.022)          & 0.264(0.083)          & 0.098(0.052)          & 0.018(0.016)          & \textbf{0.171(0.067)} & \textbf{0.049(0.035)} \\
			& Linear QR & 0.071(0.002)          & 0.535(0.013)          & 0.407(0.033)          & 0.131(0.002)          & 0.625(0.010)          & 0.567(0.030)          & 0.061(0.002)          & 0.493(0.013)          & 0.351(0.028)          \\ \midrule
			\multirow{3}{*}{$\tau=0.95$}
			& DQR & 0.034(0.016)          & 0.814(0.341)          & 1.721(1.850)          & \textbf{0.001(0.001)} & \textbf{0.184(0.030)} & \textbf{0.046(0.012)} & \textbf{0.002(0.003)} & \textbf{0.236(0.067)} & \textbf{0.085(0.042)} \\
			& Kernel QR & \textbf{0.006(0.010)} & \textbf{0.224(0.108)} & \textbf{0.086(0.086)} & 0.018(0.028)          & 0.290(0.157)          & 0.129(0.118)          & 0.014(0.006)          & 0.372(0.075)          & 0.187(0.067)          \\
			& Linear QR & 0.010(0.005)          & 0.588(0.051)          & 0.554(0.146)          & 0.039(0.001)          & 0.858(0.023)          & 1.178(0.081)          & 0.015(0.001)          & 0.466(0.024)          & 0.310(0.046)          \\ \bottomrule
		\end{tabular}%
	}
\end{table}

\subsection{Data generation: multivariate models}
Throughout the multivariate model simulation, we set the input dimension $d=6$ and sample $X$ uniformly on $[0,1]^6$. We consider the models in Section \ref{compositionR} including single index model and additive model 
which correspond different specifications of $f_0$.  The formulae of are given below.
\begin{enumerate}
	\item Single index model:
	$$f_0(x)=\exp(\theta^\top x),$$
	where $\theta=(2.2831,-1.4818,5.1966,0,0,0.0515)^\top\in\mathbb{R}^6$.
	\item Additive model:
	$$f_0(x)=\exp(4(x_1-0.5))+9(x_2-0.5)^2+10\sin(2\pi x_3)-7\vert x_4-0.5\vert,$$
	where $x=(x_1,\ldots,x_6)^\top\in[0,1]^6.$
\end{enumerate}
And we generate the error $\eta$ from following distributions,
\begin{enumerate}
	\item $\eta$ follows a scaled Student's t distribution with degree of freedom 3, i.e., $\eta\sim 0.5\times t(3)$, denoted by $t(3)$;
	\item  Conditioning on $X=x$, the error $\eta$ follows a normal distribution whose variance depends on the covariate $X$, denoted by \textit{Sine}, i.e.,
	$$\eta\mid X=x\sim 0.5\times\mathcal{N}(0,\vert\sin(\pi \xi^\top x)\vert^2)$$
	where $\xi=(1.8100,-1.2999,0,0,-2.7874,0.3197)^\top\in\mathbb{R}^d$;
	\item Conditioning on $X=x$, the error $\eta$ follows a normal distribution whose variance depends on the covariate $X$, denoted by \textit{Exp}, i.e., $$\eta\mid X=x\sim0.5\times\mathcal{N}(0,\exp(4\xi^\top x-2))$$
	where $\xi=(1.8100,-1.2999,0,0,-2.7874,0.3197)^\top\in\mathbb{R}^d$.
\end{enumerate}
Similarly, the $\tau$-th conditional quantile $f_0^\tau(x)$ of  $Y$ given $X=x$ can be calculated by
$$f_0^\tau(x)=f_0(x)+F^{-1}_{\eta\mid X=x}(\tau),$$
where $F^{-1}_{\eta\mid X=x}(\cdot)$ is the inverse of the conditional cumulated distribution function of $\eta$ given $X=x$.

We generate training data with sample sizes $n=512, 1024$  and train the estimators in the same way as in the univariate model simulations.
Summary measures including the excess risks and the $L_1$ test and the $L_2^2$ errors
based on $R=10$ replications are summarized
 in  Tables \ref{tab:sim1A}-\ref{tab:add2A}.

We see that for the nonlinear multivariate models considered, especially for single index model, \textit{DQR} performs significantly better than \textit{kernel QR} and \textit{linear QR} across all settings of error distributions.

\begin{landscape}
\begin{table}[H]
	\setlength{\tabcolsep}{0.5pt} 
	\renewcommand{\arraystretch}{1.6} 
	\centering
	\caption{\footnotesize Data is generated from single index model with training sample size $n = 512$ and the number of replications $R = 10$. The averaged excess risks, the $L_1$ and the $L_2^2$ test errors with the corresponding standard deviations (in parentheses) are reported for the estimators trained by different methods.}
	\label{tab:sim1A}
	\resizebox{1.4\textwidth}{!}{%
		\begin{tabular}{@{}cc|ccc|ccc|ccc@{}}
			\toprule
			\multicolumn{2}{c|}{$n=512$}              & \multicolumn{3}{c|}{$t(3)$}                                                    & \multicolumn{3}{c|}{\textit{Sine}}                                                   & \multicolumn{3}{c}{\textit{Exp}}                                                        \\
			Quantile                     & Method    & Excess risk             & $L_1$                 & $L_2^2$                     & Excess risk             & $L_1$                 & $L_2^2$                    & Excess risk             & $L_1$                   & $L_2^2$                      \\ \midrule
			\multirow{3}{*}{$\tau=0.05$}
			& DQR       & \textbf{4.217(1.309)}   & \textbf{2.151(0.217)} & \textbf{88.164(21.657)}     & \textbf{4.307(1.347)}   & \textbf{2.075(0.280)} & \textbf{89.774(23.802)}    & \textbf{4.238(1.478)}   & \textbf{2.325(0.548)}   & \textbf{93.056(26.654)}      \\
			& Kernel QR & 294.188(17.050)         & 32.817(1.023)         & 5806.695(341.938)           & 294.540(17.043)         & 33.336(1.001)         & 5853.665(341.673)          & 294.514(16.857)         & 33.551(0.991)           & 5850.878(336.681)            \\
			& Linear QR & 713.114(11.642)         & 49.666(0.071)         & 14146.989(232.853)          & 714.349(12.146)         & 50.166(0.085)         & 14234.740(242.944)         & 714.942(13.323)         & 50.388(0.118)           & 14239.648(265.850)           \\ \midrule
			\multirow{3}{*}{$\tau=0.25$}
			& DQR       & \textbf{17.049(12.019)} & \textbf{2.498(0.779)} & \textbf{112.648(42.008)}    & \textbf{17.706(11.381)} & \textbf{2.082(0.452)} & \textbf{95.818(37.391)}    & \textbf{17.663(10.833)} & \textbf{2.185(0.449)}   & \textbf{99.205(29.985)}      \\
			& Kernel QR & 1299.841(98.461)        & 26.366(0.389)         & 5209.186(385.579)           & 1301.009(98.121)        & 26.400(0.371)         & 5220.212(384.629)          & 1301.021(98.439)        & 26.444(0.377)           & 5216.858(385.396)            \\
			& Linear QR & 3406.749(88.870)        & 47.584(0.446)         & 13594.697(354.328)          & 3408.650(80.138)        & 47.701(0.418)         & 13616.544(319.787)         & 3402.889(84.862)        & 47.760(0.433)           & 13589.463(338.430)           \\ \midrule
			\multirow{4}{*}{$\tau=0.5$}  & DLS       & 97.952(46.861)          & \textbf{2.069(0.237)} & \textbf{97.997(46.882)}     & 98.776(38.518)          & \textbf{2.270(1.326)} & \textbf{98.758(38.524)}    & 87.272(26.925)          & \textbf{1.790(0.183)}   & \textbf{87.180(26.877)}      \\
			& DQR       & \textbf{31.482(30.991)} & 6.080(3.396)          & 286.932(241.121)            & \textbf{24.204(22.286)} & 5.043(3.085)          & 308.133(260.270)           & \textbf{33.683(26.674)} & 4.260(2.693)            & 214.805(213.108)             \\
			& Kernel QR & 2358.607(213.788)       & 24.260(0.345)         & 4832.728(399.260)           & 2362.170(213.940)       & 24.258(0.356)         & 4839.095(398.853)          & 2363.701(213.499)       & 24.254(0.356)           & 4840.900(398.880)            \\
			& Linear QR & 5664.872(282.831)       & 44.817(0.248)         & 11472.181(525.495)          & 5669.251(287.263)       & 44.832(0.253)         & 11479.596(534.587)         & 5667.212(289.233)       & 44.838(0.260)           & 11475.803(538.182)           \\ \midrule
			\multirow{3}{*}{$\tau=0.75$}
			& DQR       & \textbf{49.370(34.458)} & \textbf{5.586(4.870)} & \textbf{265.765(265.944)}   & \textbf{44.422(42.234)} & \textbf{3.007(1.903)} & \textbf{144.606(97.917)}   & \textbf{27.519(27.786)} & \textbf{9.203(4.743)}   & \textbf{577.012(424.183)}    \\
			& Kernel QR & 3293.028(311.423)       & 26.022(0.338)         & 4665.851(367.340)           & 3298.689(308.738)       & 26.101(0.334)         & 4670.822(363.761)          & 3299.316(308.809)       & 26.168(0.332)           & 4678.233(365.136)            \\
			& Linear QR & 5410.394(496.133)       & 58.366(3.076)         & 9347.631(283.711)           & 5419.965(499.627)       & 58.299(3.023)         & 9340.783(292.640)          & 5422.321(499.931)       & 58.336(2.989)           & 9351.314(295.154)            \\ \midrule
			\multirow{3}{*}{$\tau=0.95$}
			& DQR       & \textbf{36.993(37.340)} & \textbf{9.674(9.633)} & \textbf{859.259(1392.082)}  & \textbf{52.139(52.257)} & \textbf{8.890(6.794)} & \textbf{891.320(1202.299)} & \textbf{17.656(22.966)} & \textbf{15.673(20.930)} & \textbf{4280.954(10557.121)} \\
			& Kernel QR & 2642.302(523.477)       & 88.606(22.027)        & 10448.269(4384.084)         & 2641.926(519.902)       & 89.516(22.610)        & 10625.592(4586.608)        & 2644.186(521.020)       & 89.755(22.492)          & 10668.437(4585.526)          \\
			& Linear QR & 1807.445(532.122)       & 168.879(20.744)       & 35177.445(7645.388)         & 1808.274(523.959)       & 169.236(20.580)       & 35218.465(7600.995)        & 1816.628(526.056)       & 169.271(20.542)         & 35265.383(7589.261)          \\ \midrule
		\end{tabular}%
	}
\end{table}
\end{landscape}

\begin{landscape}
	\begin{table}[H]
		\setlength{\tabcolsep}{0.5pt} 
		\renewcommand{\arraystretch}{1.6} 
		\centering
		\caption{\footnotesize Data is generated from single index model with training sample size $n=1024$ and the number of replications $R = 10$. The averaged excess risks, the $L_1$ and the $L_2^2$ test errors with the corresponding standard deviations (in parentheses) are reported for the estimators trained by different methods.}
		\label{tab:sim2A}
		\resizebox{1.4\textwidth}{!}{%
\begin{tabular}{@{}cc|ccc|ccc|ccc@{}}
\toprule
\multicolumn{2}{c|}{$n=1024$}             & \multicolumn{3}{c|}{$t(3)$}                                                    & \multicolumn{3}{c|}{\textit{Sine}}                   &\multicolumn{3}{c}{\textit{Exp}}                                                        \\
Quantile                     & Method    & Excess risk             & $L_1$                 & $L_2^2$                     & Excess risk             & $L_1$                 & $L_2^2$                    & Excess risk             & $L_1$                   & $L_2^2$                      \\ \midrule
\multirow{3}{*}{$\tau=0.05$}
& DQR & \textbf{1.754(0.679)}   & \textbf{1.794(0.782)} & \textbf{33.969(14.685)}     & \textbf{2.037(0.684)}   & \textbf{1.437(0.298)} & \textbf{39.895(12.900)}    & \textbf{1.377(0.733)}   & \textbf{1.493(0.368)}   & \textbf{28.311(13.635)}      \\
& Kernel QR & 238.917(8.306)          & 31.335(1.258)         & 4703.625(168.559)           & 247.563(8.299)          & 30.801(0.929)         & 4917.161(166.789)          & 238.968(8.434)          & 32.075(1.240)           & 4743.118(168.673)            \\
& Linear QR & 708.536(9.963)          & 49.697(0.070)         & 14053.810(199.273)          & 716.733(11.356)         & 50.193(0.039)         & 14282.439(227.090)         & 711.296(9.158)          & 50.409(0.031)           & 14166.998(182.722)           \\ \midrule
\multirow{3}{*}{$\tau=0.25$}
& DQR & \textbf{4.828(3.762)}   & \textbf{1.731(0.538)} & \textbf{36.095(14.611)}     & \textbf{7.171(2.536)}   & \textbf{1.631(0.674)} & \textbf{36.184(8.533)}     & \textbf{4.617(2.491)}   & \textbf{1.326(0.389)}   & \textbf{28.375(14.217)}      \\
& Kernel QR & 1004.738(65.562)        & 23.378(0.402)         & 4031.585(255.732)           & 1071.526(55.534)        & 23.422(0.291)         & 4297.444(217.730)          & 1005.130(66.363)        & 23.474(0.389)           & 4036.355(258.944)            \\
& Linear QR & 3385.564(63.690)        & 47.578(0.186)         & 13508.011(254.175)          & 3415.326(65.336)        & 47.705(0.326)         & 13642.956(260.718)         & 3386.482(63.084)        & 47.693(0.192)           & 13523.764(251.627)           \\ \midrule
\multirow{4}{*}{$\tau=0.5$}  & DLS       & 21.687(8.463)           & \textbf{1.149(0.142)} & \textbf{21.697(8.463)}      & 33.510(9.535)           & \textbf{1.211(0.257)} & \textbf{33.506(9.538)}     & 21.354(9.077)           & \textbf{1.057(0.153)}   & \textbf{21.324(9.057)}       \\
& DQR & \textbf{9.027(7.164)}   & 3.180(1.869)          & 126.745(124.557)            & \textbf{12.998(8.342)}  & 2.523(2.185)          & 89.498(119.938)            & \textbf{7.250(5.798)}   & 2.736(2.025)            & 85.600(98.329)               \\
& Kernel QR & 1783.861(143.639)       & 21.357(0.201)         & 3691.269(260.875)           & 1923.828(120.525)       & 21.365(0.198)         & 3945.973(222.487)          & 1785.665(143.829)       & 21.343(0.194)           & 3694.289(261.772)            \\
& Linear QR & 5665.659(133.961)       & 44.579(0.086)         & 11446.959(251.709)          & 5732.970(190.640)       & 44.607(0.066)         & 11585.740(356.092)         & 5677.712(125.335)       & 44.573(0.081)           & 11469.584(235.753)           \\ \midrule
\multirow{3}{*}{$\tau=0.75$}
& DQR & \textbf{15.079(15.889)} & \textbf{5.938(7.075)} & \textbf{358.676(685.721)}   & \textbf{19.547(19.361)} & \textbf{6.709(6.778)} & \textbf{670.137(957.505)}  & \textbf{7.565(9.392)}   & \textbf{6.098(5.699)}   & \textbf{334.212(525.110)}    \\
& Kernel QR & 2451.441(213.667)       & 23.222(0.369)         & 3560.982(232.572)           & 2660.367(184.853)       & 23.109(0.346)         & 3790.041(204.149)          & 2454.610(216.033)       & 23.335(0.372)           & 3570.188(234.718)            \\
& Linear QR & 5331.374(421.982)       & 58.032(3.025)         & 9186.867(200.708)           & 5613.486(366.482)       & 57.059(2.377)         & 9395.197(209.270)          & 5324.493(423.213)       & 58.261(3.097)           & 9206.238(201.860)            \\ \midrule
\multirow{3}{*}{$\tau=0.95$}
& DQR & \textbf{8.325(14.732)}  & \textbf{9.609(7.311)} & \textbf{1056.564(1168.020)} & \textbf{24.292(23.253)} & \textbf{7.857(5.090)} & \textbf{611.297(644.296)}  & \textbf{10.020(13.874)} & \textbf{15.631(17.196)} & \textbf{2527.314(4691.837)}  \\
& Kernel QR & 2046.284(295.639)       & 75.767(10.299)        & 7536.021(1353.733)          & 2287.394(249.019)       & 71.394(7.767)         & 7009.876(943.809)          & 2051.629(293.556)       & 76.384(10.065)          & 7623.915(1347.725)           \\
& Linear QR & 1555.293(343.043)       & 176.002(17.312)       & 37711.316(6766.304)         & 1830.525(271.066)       & 164.519(13.110)       & 33222.203(4909.937)        & 1555.301(336.023)       & 176.830(17.026)         & 38014.789(6691.186)          \\ \bottomrule
		\end{tabular}%
}
\end{table}
\end{landscape}

\begin{landscape}
	\begin{table}[H]
		\setlength{\tabcolsep}{2pt} 
		\renewcommand{\arraystretch}{1.3} 
		\centering
		\caption{\footnotesize Data is generated from additive model with training sample size $n=512$ and the number of replications $R = 10$. The averaged excess risks, the $L_1$ and the $L_2^2$ test errors with the corresponding standard deviations (in parentheses) are reported for the estimators trained by different methods.}
		\label{tab:add1A}
		\resizebox{1.4\textwidth}{!}{%
			\begin{tabular}{@{}cc|ccc|ccc|ccc@{}}
				\toprule
				\multicolumn{2}{c|}{$n=512$}              & \multicolumn{3}{c|}{$t(3)$}                                            & \multicolumn{3}{c|}{\textit{Sine}}                                           & \multicolumn{3}{c}{\textit{Exp}}                                            \\
				Quantile                     & Method    & Excess risk           & $L_1$                 & $L_2^2$               & Excess risk           & $L_1$                 & $L_2^2$               & Excess risk           & $L_1$                 & $L_2^2$               \\ \midrule
				\multirow{3}{*}{$\tau=0.05$}
				& DQR       & \textbf{0.205(0.022)} & \textbf{0.983(0.045)} & \textbf{1.471(0.142)} & \textbf{0.101(0.033)} & \textbf{0.555(0.141)} & \textbf{0.521(0.274)} & \textbf{0.072(0.017)} & \textbf{0.536(0.084)} & \textbf{0.839(0.151)} \\
				& Kernel QR & 1.556(0.039)          & 3.791(0.052)          & 22.082(0.777)         & 1.487(0.065)          & 4.001(0.115)          & 24.185(1.344)         & 1.527(0.035)          & 4.297(0.090)          & 27.396(0.949)         \\
				& Linear QR & 3.778(0.194)          & 6.413(0.212)          & 60.737(3.483)         & 3.711(0.162)          & 6.796(0.192)          & 66.570(3.205)         & 3.743(0.186)          & 7.058(0.220)          & 70.985(3.841)         \\ \midrule
				\multirow{3}{*}{$\tau=0.25$}
				& DQR       & \textbf{0.283(0.038)} & \textbf{0.751(0.045)} & \textbf{0.983(0.165)} & \textbf{0.141(0.033)} & \textbf{0.440(0.054)} & \textbf{0.344(0.089)} & \textbf{0.117(0.027)} & \textbf{0.443(0.091)} & \textbf{0.439(0.124)} \\
				& Kernel QR & 4.476(0.295)          & 3.291(0.060)          & 16.293(0.941)         & 4.218(0.225)          & 3.249(0.074)          & 15.863(0.806)         & 4.547(0.272)          & 3.398(0.073)          & 17.589(1.014)         \\
				& Linear QR & 9.198(0.775)          & 4.787(0.168)          & 34.673(2.557)         & 8.969(0.399)          & 4.792(0.093)          & 34.818(1.421)         & 9.506(0.800)          & 4.961(0.191)          & 37.463(2.899)         \\ \midrule
				\multirow{4}{*}{$\tau=0.5$}  & DLS       & 0.926(0.149)          & 0.721(0.046)          & 0.924(0.148)          & 0.279(0.041)          & \textbf{0.400(0.027)} & \textbf{0.280(0.041)} & 0.261(0.072)          & \textbf{0.346(0.050)} & \textbf{0.263(0.071)} \\
				& DQR       & \textbf{0.349(0.077)} & \textbf{0.716(0.062)} & \textbf{0.894(0.194)} & \textbf{0.161(0.025)} & 0.449(0.032)          & 0.353(0.058)          & \textbf{0.157(0.047)} & 0.395(0.069)          & 0.329(0.113)          \\
				& Kernel QR & 3.627(0.487)          & 2.903(0.035)          & 11.458(0.325)         & 3.207(0.393)          & 2.846(0.036)          & 10.992(0.304)         & 3.555(0.470)          & 2.879(0.039)          & 11.251(0.344)         \\
				& Linear QR & 7.035(0.777)          & 4.033(0.036)          & 23.116(0.621)         & 6.524(0.627)          & 4.041(0.027)          & 23.478(0.546)         & 7.166(0.735)          & 4.037(0.033)          & 23.279(0.511)         \\ \midrule
				\multirow{3}{*}{$\tau=0.75$}
				& DQR       & \textbf{0.448(0.084)} & \textbf{0.795(0.055)} & \textbf{1.073(0.200)} & \textbf{0.180(0.039)} & \textbf{0.496(0.047)} & \textbf{0.418(0.077)} & \textbf{0.178(0.035)} & \textbf{0.405(0.068)} & \textbf{0.387(0.094)} \\
				& Kernel QR & 1.584(0.258)          & 3.211(0.074)          & 15.307(0.835)         & 1.413(0.157)          & 3.302(0.086)          & 16.584(0.910)         & 1.632(0.233)          & 3.309(0.091)          & 16.540(1.051)         \\
				& Linear QR & 2.474(0.281)          & 4.692(0.112)          & 33.280(1.672)         & 2.437(0.320)          & 4.882(0.177)          & 36.300(2.648)         & 2.555(0.267)          & 4.840(0.113)          & 35.596(1.674)         \\ \midrule
				\multirow{3}{*}{$\tau=0.95$}
				& DQR       & 0.479(0.113)          & \textbf{1.137(0.117)} & \textbf{1.908(0.427)} & \textbf{0.149(0.050)} & \textbf{0.584(0.046)} & \textbf{0.554(0.097)} & \textbf{0.184(0.051)} & \textbf{0.531(0.119)} & \textbf{0.856(0.185)} \\
				& Kernel QR & 0.378(0.118)          & 3.787(0.100)          & 21.959(1.117)         & 0.375(0.070)          & 4.120(0.108)          & 25.568(1.178)         & 0.471(0.100)          & 4.290(0.062)          & 27.340(0.811)         \\
				& Linear QR & \textbf{0.366(0.031)} & 6.393(0.137)          & 60.086(2.248)         & 0.432(0.027)          & 6.900(0.201)          & 67.984(3.290)         & 0.474(0.021)          & 7.014(0.139)          & 69.932(2.296)         \\ \bottomrule
			\end{tabular}%
		}
	\end{table}
\end{landscape}

\begin{landscape}
	\begin{table}[H]
		\setlength{\tabcolsep}{2pt} 
		\renewcommand{\arraystretch}{1.3} 
		\centering
		\caption{\footnotesize Data is generated from additive model with training sample size $n=1024$ and the number of replications $R = 10$. The averaged excess risks, the $L_1$ and the $L_2^2$ test errors and the corresponding standard deviations (in parentheses) are reported for the estimators trained by different methods.}
		\label{tab:add2A}
		\resizebox{1.4\textwidth}{!}{%
			\begin{tabular}{@{}ll|ccc|ccc|ccc@{}}
				\toprule
				\multicolumn{2}{c|}{$n=1024$}             & \multicolumn{3}{c|}{$t(3)$}                                            & \multicolumn{3}{c|}{\textit{Sine}}                                           & \multicolumn{3}{c}{\textit{Exp}}                                            \\
				Quantile                     & Method    & Excess risk           & $L_1$                 & $L_2^2$               & Excess risk           & $L_1$                 & $L_2^2$               & Excess risk           & $L_1$                 & $L_2^2$               \\ \midrule
				\multirow{3}{*}{$\tau=0.05$}
				& DQR & \textbf{0.173(0.047)} & \textbf{0.915(0.116)} & \textbf{1.272(0.277)} & \textbf{0.076(0.021)} & \textbf{0.499(0.087)} & \textbf{0.411(0.156)} & \textbf{0.046(0.010)} & \textbf{0.562(0.134)} & \textbf{0.877(0.205)} \\
				& Kernel QR & 0.765(0.036)          & 2.526(0.070)          & 9.804(0.545)          & 0.746(0.030)          & 2.764(0.063)          & 11.633(0.525)         & 0.734(0.035)          & 2.980(0.071)          & 13.220(0.619)         \\
				& Linear QR & 3.735(0.116)          & 6.378(0.123)          & 60.045(2.097)         & 3.704(0.076)          & 6.788(0.091)          & 66.504(1.507)         & 3.698(0.108)          & 7.014(0.132)          & 70.101(2.202)         \\ \midrule
				\multirow{3}{*}{$\tau=0.25$}
				& DQR & \textbf{0.257(0.034)} & \textbf{0.704(0.061)} & \textbf{0.857(0.134)} & \textbf{0.104(0.011)} & \textbf{0.384(0.041)} & \textbf{0.251(0.046)} & \textbf{0.085(0.019)} & \textbf{0.362(0.064)} & \textbf{0.358(0.080)} \\
				& Kernel QR & 2.216(0.125)          & 2.265(0.059)          & 7.582(0.441)          & 2.173(0.118)          & 2.273(0.059)          & 7.733(0.445)          & 2.203(0.119)          & 2.332(0.053)          & 8.219(0.415)          \\
				& Linear QR & 9.038(0.378)          & 4.745(0.084)          & 33.998(1.251)         & 9.101(0.417)          & 4.818(0.096)          & 35.165(1.438)         & 9.064(0.415)          & 4.850(0.100)          & 35.712(1.488)         \\ \midrule
				\multirow{4}{*}{$\tau=0.5$}  & DLS       & 0.861(0.127)          & 0.705(0.050)          & 0.858(0.126)          & 0.214(0.058)          & \textbf{0.354(0.051)} & \textbf{0.214(0.058)} & 0.233(0.049)          & \textbf{0.282(0.027)} & 0.233(0.048)          \\
				& DQR & \textbf{0.321(0.067)} & \textbf{0.685(0.065)} & \textbf{0.812(0.155)} & \textbf{0.110(0.017)} & 0.377(0.036)          & 0.236(0.039)          & \textbf{0.103(0.021)} & 0.296(0.055)          & \textbf{0.232(0.053)} \\
				& Kernel QR & 1.662(0.155)          & 2.020(0.038)          & 5.375(0.206)          & 1.630(0.146)          & 1.969(0.044)          & 5.072(0.217)          & 1.684(0.163)          & 1.974(0.040)          & 5.107(0.201)          \\
				& Linear QR & 6.144(0.404)          & 4.004(0.011)          & 22.629(0.463)         & 6.234(0.375)          & 4.004(0.011)          & 22.707(0.468)         & 6.290(0.412)          & 4.004(0.012)          & 22.702(0.462)         \\ \midrule
				\multirow{3}{*}{$\tau=0.75$}
				& DQR & \textbf{0.360(0.095)} & \textbf{0.731(0.059)} & \textbf{0.896(0.135)} & \textbf{0.112(0.020)} & \textbf{0.422(0.037)} & \textbf{0.286(0.043)} & \textbf{0.136(0.038)} & \textbf{0.299(0.030)} & \textbf{0.298(0.069)} \\
				& Kernel QR & 0.703(0.053)          & 2.284(0.036)          & 7.744(0.273)          & 0.729(0.045)          & 2.294(0.046)          & 7.928(0.355)          & 0.780(0.052)          & 2.355(0.053)          & 8.423(0.423)          \\
				& Linear QR & 2.235(0.224)          & 4.808(0.114)          & 35.111(1.710)         & 2.299(0.189)          & 4.878(0.116)          & 36.159(1.787)         & 2.350(0.170)          & 4.936(0.116)          & 37.083(1.788)         \\ \midrule
				\multirow{3}{*}{$\tau=0.95$}
				& DQR & 0.336(0.108)          & \textbf{1.002(0.131)} & \textbf{1.477(0.294)} & \textbf{0.100(0.015)} & \textbf{0.528(0.036)} & \textbf{0.430(0.046)} & \textbf{0.144(0.054)} & \textbf{0.440(0.084)} & \textbf{0.718(0.209)} \\
				& Kernel QR & \textbf{0.132(0.017)} & 2.564(0.065)          & 10.104(0.546)         & 0.177(0.013)          & 2.759(0.056)          & 11.613(0.496)         & 0.220(0.026)          & 2.971(0.072)          & 13.162(0.628)         \\
				& Linear QR & 0.338(0.014)          & 6.409(0.121)          & 60.352(1.964)         & 0.422(0.015)          & 6.873(0.127)          & 67.675(2.031)         & 0.450(0.014)          & 7.134(0.105)          & 71.839(1.767)         \\ \bottomrule
			\end{tabular}%
		}
	\end{table}
\end{landscape}

\end{document}